\documentclass[12pt]{article}

\usepackage[english]{babel}
\usepackage{amsfonts}
\usepackage{amsmath}
\usepackage{amssymb, latexsym}
\usepackage{graphicx}
\usepackage[usenames]{color}

\usepackage{fullpage}
\usepackage{amsmath}

\usepackage{mathrsfs}
\usepackage{ulem}

\newtheorem{theorem}{Theorem}[section]
\newtheorem{proposition}[theorem]{Proposition}

\newtheorem{lemma}[theorem]{Lemma}
\newtheorem{remark}[theorem]{Remark}
\newtheorem{definition}[theorem]{Definition}

\numberwithin{equation}{section} 

\def\calx{\mathcal{X}}
\def\caly{\mathcal{Y}}
\def\bbc{{\mathbb{C}}}
\def\bbe{{\mathbb{E}}}
\def\bbr{{\mathbb{R}}}
\def\bbp{{\mathbb{P}}}

\def\bbx{{\mathbb{X}}}
\def\bby{{\mathbb{Y}}}
\def\bbw{{\mathbb{W}}}
\renewcommand{\Re}{{\rm Re} \,}

\def\en t{{{\rm Z}\mkern-5.5mu{\rm Z}}}

\def\<{\left<}
\def\>{\right>}
\def\({\left(}
\def\){\right)}

\def\9{{\infty}}
\def\barr{\begin{array}}
\def\earr{\end{array}}

\def\wt{\widetilde}
\def\wh{\widehat}
\def\ol{\overline}

\def\vf{{\varphi}}
\def\lbb{{\lambda}}
\def\g{{\gamma}}
\def\a{{\alpha}}

\def\n{\noindent }

\def\3{\subset }
\def\na{{\nabla}}
\def\sk{\smallskip }

\def\bk{\bigskip }

\def\ve{{\varepsilon}}
\def\p{{\partial}}

\begin{document}

\begin{center}
{\Large{\bf Stochastic nonlinear Schr\"odinger equations
in the defocusing mass and energy critical cases}}
\bigskip\bk

{\large{\bf Deng Zhang}}
\footnote[2]{Department of Mathematics,
Shanghai Jiao Tong University, 200240 Shanghai, China. \\
Email: dzhang@sjtu.edu.cn }
\end{center}

\bk\bk\bk

\begin{quote}
\n{\small{\bf Abstract.}
We study the stochastic nonlinear Schr\"odinger equations with linear multiplicative noise,
particularly in the defocusing mass-critical and energy-critical cases.
For general initial data,
we prove the global existence and uniqueness of solutions in both cases.
When the quadratic variation of noise is globally bounded,
we also obtain the rescaled scattering behavior of solutions in the
spaces $L^2$, $H^1$ as well as  the pseudo-conformal space.
Moreover,
the Stroock-Varadhan type theorem for the topological support of
solutions are also obtained
in the Strichartz and local smoothing spaces.  }


\end{quote}

\vfill

\section{Introduction} \label{Sec-Intro}

This work is devoted to stochastic nonlinear Schr\"odinger equations with linear multiplicative noise
in the defocusing mass-critical and energy-critical cases,
it is a continuation of a series of work \cite{BRZ14,BRZ16,HRZ18}.
Precisely, we consider
\begin{equation} \label{equa-x}
\begin{split}
   idX&= \Delta Xdt   + \lbb F(X) dt  -i\mu Xdt + i\sum\limits_{k=1}^N  X  G_k d\beta_k(t)
,  \\
   X(0)&=X_0.
 \end{split}
\end{equation}
Here,
the nonlinearity $F(X)= |X|^{\a-1}X$,
$\a>1$,
$\lbb =-1$ (resp. $\lbb=1$) corresponds to the defocusing (resp. focusing) case,
$\beta_k$ are standard real-valued
Brownian motions on a probability space $(\Omega, \mathscr{F},
\mathbb{P})$ with normal filtration
$(\mathscr{F}_t)_{t\geq 0}$,
and $G_k(t,x)=g_k(t)\phi_k(x)$,
where $g_k$ are real-valued predictable processes satisfying
$g_k \in L^2_{loc}(\bbr^+; \bbr)$ $\bbp$-a.s.,
and $\phi_k\in C^\9(\bbr^d; \mathbb{C})$,
$d\geq 1$.
For simplicity, we assume $N<\9$,
but the arguments in this paper extend also to the case where
$N=\9$ under appropriate summable conditions on $G_k$.

Moreover,
the term
$\mu$  is of the form
$$\mu(t,x) = \frac12 \sum_{k=1}^N |G_k(t,x)|^2. $$
In particular,
in the {\it conservative case} where $\Re G_k = 0$,
$1\leq k\leq N$,
$-i\mu Xdt + i\sum_{k=1}^N  X  G_k d\beta_k(t)$
is indeed the Stratonovitch differential,
and, via It\^o's formula,
the mass is pathwisely conserved $|X(t)|_2^2 =|X_0|_2^2$, $t\in [0,T]$.
Hence, for the normalized initial state $|X_0|_2=1$,
the quantum system evolves on the unit bass of $L^2$
and verifies the conservation of probability.
See, e.g., \cite{BCIR94, BCIRG95} for  applications
in molecular aggregates with thermal fluctuations.

The {\it non-conservative case} (i.e.,
$\Re G_k \not = 0$ for some $1\leq k\leq N$)
plays an important role in the application to open quantum systems \cite{BG09},
one of the main features is that
$t \mapsto |X(t)|_2^2$
is a continuous martingale.
This fact
implies the mean norm square conservation $\bbe |X(t)|_2^2$, $t\in [0,T]$,
and enables one to define the ``physical'' probability law
\begin{align*}
    \wh{\mathbb{P}}^T_{X_0} (d \omega)
    := (\mathbb{E}_{\mathbb{P}} [|X_0|_2^2])^{-1} |X(T,\omega)|_2^2\ \mathbb{P} (d \omega)
\end{align*}
of the events occurring in $[0,T]$.
We refer to \cite{BG09} for more details.
For more physical applications, e.g. nonlinear optics,
Bose-Einstein condensation and
the Gross-Pitaevskii equation,
we refer to \cite{SS99}.

For stochastic nonlinear Schr\"odinger equations,
most results in literature center around the subcritical case.
The first global well-posedness results
were  proved by de Bouard and Debussche \cite{BD99,BD03},
by using the theory of radonifying operators.
Later,
the compact manifold case was studied by
Brze\'{z}niak and Millet  \cite{BM14},
where  more general stochastic Strichartz estimates were proved.
See also \cite{BHW17, BHW18}.
Recently, the global well-posedness of \eqref{equa-x}
for the full subcritical exponents was proved in \cite{BRZ14,BRZ16},
the new method introduced is the rescaling approach,
which can be viewed as Doss-Dussman type transformations in Hilbert spaces.
We also refer to  \cite{H16} for the global well-posedness
in the full mass-subcritical case with quite general multiplicative noise.
See also \cite{BRZ16, BRZ18,CG15,Z17}.

The critical case is much more subtle
and, to the best of our knowledge,
quite few results are known  for  stochastic nonlinear Schr\"odinger equations.
See, e.g., \cite{BRZ14,BRZ16,BD99,BD03,H16} for the local well-posedenss results.
The main difference between the subcritical and critical cases is,
that the maximal existing time of solutions depends only on the $L^2$- or $H^1$-norm of initial data
in the subcritical case,
while on the whole profile  in the critical case.
Hence,
the standard energy method works well for the global well-posedness in the subscritical case,
however, it fails in the critical case.

In contrast,
the critical case in the deterministic case
has been extensively studied in literature.
In the defocusing mass- and energy-critical cases,
it was conjectured that deterministic solutions exist globally and even scatter at infinity.
This conjecture was first proved,
via the energy induction method,
by Bourgain  in the seminal work \cite{B99}
for the energy-critical case with radial initial data in dimensions three and four.
Later, for general initial data and dimensions,
it was proved
by the I-team \cite{CKSTT08}, Ryckman and Visan \cite{RV07}
and Visan \cite{V07},
based on the energy induction method and
interaction Morawetz estimates.
See also the concentration compactness method introduced in
\cite{KM06}.
Recently,
this conjecture in the mass-critical case
was proved by Dodson \cite{D12,D16.1,D16.2} for general initial data,
where the key ingredients are long-time Strichartz estimates.
See also \cite{KTV14, T06}.

However,
it is quite hard to obtain these estimates in the stochastic case.
Actually,
the presence of noise in \eqref{equa-x}
destroys the symmetries of  equation and the conservation laws (e.g. of the mass and the Hamiltonian),
the corresponding It\^o formulas
consist of several stochastic integrals,
with which it is very difficult to obtain sharp estimates
as in the deterministic  case.
Moreover,
even a Banach space $\mathscr{X}$ is compactly embedded into another one $\mathscr{Y}$,
one does not generally have the compact embedding of $L^p(\Omega; \mathscr{X})$
into $L^p(\Omega; \mathscr{Y})$, $p\geq 1$.

Hence, the global existence of stochastic solutions in the mass- and energy-critical cases
with general initial data
has been an open problem.
See the recent progress \cite{FX18.1,FX18.2} 
for the global well-posedness
in the conservative mass-critical case for dimension $d=1$.
(See also Remark \ref{Rem-Com-FANXU} below.)

In this paper, we prove  the global well-posedness
of \eqref{equa-x}
in the mass-critical case
for all dimensions $d\geq 1$.
Moreover,  in the energy-critical case,
we prove the global well-posedness
for dimensions $3\leq d\leq 6$,
and we are also able to prove conditional global well-posedness results
for high dimensions $d>6$,
assuming an {\it a-priori} bound of the energy.

Thus, together with the previous work \cite{BRZ14,BRZ16}
and  the {\it a-priori} bound of the energy in the energy-critical case with $d>6$,
the global existence and uniqueness of solutions to \eqref{equa-x}
are obtained
for the full subcritical and critical exponents of the nonlinearity in the defocusing case.
We  would like also to mention that,
these results also apply to the non-conservative case,
which is important in the physical context \cite{BG09}.

The proof presented below is
different from that in \cite{FX18.1,FX18.2}
and is
based on a new application of the rescaling approach.
It is also based on the work \cite{CKSTT08, D12,D16.1,D16.2, RV07,V07} mentioned above
and on the stability results for
nonlinear Schr\"odinger equations with lower order perturbations
(see Theorems \ref{Thm-Sta-L2}, \ref{Thm-Sta-H1-dlow} and \ref{Thm-Sta-H1} below),
which are also of independent interest.

Another main interest of this paper lies in the large time behavior of global solutions to \eqref{equa-x}.
As mentioned above,
in the defocusing case
deterministic global solutions scatter at infinity,
namely, behave asymptotically like linear solutions.
However, the situation becomes quite difficult in the stochastic case
because of  rapid fluctuations of noise at large time.
Very recently, in \cite{HRZ18},
the rescaled scattering behavior of global  solutions to \eqref{equa-x}
is proved for the energy-subscritical exponents
$\a \in [\max\{2, 1+\frac 4d\}, 1+\frac{4}{d-2})$,
$3\leq d\leq 6$,
and
it is also proved that
the non-conservative noise  has the effect to improve scattering with high probability,
even in the regime where deterministic solutions fail to scatter.
The energy-critical case is also studied there,
however, relying on the {\it a-priori} assumption of global existence of solutions,
which is another motivation for the present work.

When the quadratic variation of noise  is globally bounded,
on the basis of \cite{HRZ18},
we   prove the rescaled scattering behavior of
global solutions
to \eqref{equa-x}
in the spaces $L^2$, $H^1$ as well as the pseudo-conformal space, respectively.
These results are new in the $L^2$ case
and also improve those of \cite{HRZ18} in the $H^1$ and  pseudo-conformal spaces.

At last,
we give a characterization of the support of the law of
global solutions to \eqref{equa-x},
in both mass-critical and energy-critical cases.
We prove that
the law of stochastic solutions is supported on the closure of all deterministic controlled trajectories
in the Strichartz and local smoothing spaces (see Theorem \ref{Thm-Sca} below).

We would like to  mention that,
for each deterministic controlled trajectory,
the global well-posedness
can be proved by stability results as in \cite{TV05,TVZ07}.
So, if the support theorem is {\it a-priori} assumed to hold,
the related stochastic trajectories should also exist globally.
This, actually, gives an intuitive point of view for the global well-posedness
of stochastic solutions
at the beginning of this work. \\

{\bf Notation.}
For $z\in \bbc$,
we set $F(z):=|z|^{\a-1}z$
with $\a=1+\frac 4d$ or $\a=1+\frac{4}{d-2}$
in the mass-critical or energy-critical case, respectively.
We denote by $F_z$ and $F_{\ol{z}}$
the usual complex derivatives
$F_z= \frac 12 (\frac{\p F}{\p x} - i \frac{\p F}{\p y})$,
$F_{\ol z}= \frac 12 (\frac{\p F}{\p x} + i \frac{\p F}{\p y})$.
For any $x=(x_1,\cdots,x_d) \in \bbr^d$
and multi-index $\a=(\a_1,\cdots, \a_d)$,
we use the notations
$|\a|= \sum_{j=1}^d \a_j$,
$\<x\>=(1+|x|^2)^{1/2}$,
$\p_j = \frac{\p}{\p{x_j}}$,
$\partial_x^\a=\partial_{x_1}^{\a_1}\cdots \partial_{x_d}^{\a_d}$,
$\<\na\>=(I-\Delta)^{1/2}$.

Let $\mathscr{S}$ denote the space of rapid decreasing functions
and $\mathscr{S}'$ be the dual space of $\mathscr{S}$.
For any $f\in \mathscr{S}$, $\mathscr{F}(f)$
is the Fourier transform of $f$,
i.e. $\mathscr{F}(f)(\xi) = \int e^{-ix\cdot \xi} f(x)dx$.
Given $1\leq p \leq \9$, $s\geq 0$,
$L^p = L^p(\bbr^d)$ is the space of $p$-integrable complex functions with the norm $|\cdot|_{L^p}$,
$W^{s,p}= \<D\>^{-s}L^p(\bbr^d)$ is the usual Sobolev space
with the norm $\|\cdot\|_{W^{s,p}}$.
In particular, we write  $|\cdot|_2 = |\cdot|_{L^2}$,
$|\cdot|_{H^1}= \|\cdot\|_{W^{1,2}}$.

For any Banach space $\calx$
and any interval $I \subseteq \bbr^+$,
$L^p(I; \calx)$ is the space of $p$-integrable $\calx$-valued functions
with the norm $\|\cdot\|_{L^p(0,T; \calx)}$,
and $C(I;\calx)$ is the space of continuous
$\calx$-valued functions with the super norm in $t$.
Moreover,
for two Banach spaces $\calx, \caly$,
the norm of $\calx \cap \caly$ is
$\|\cdot\|_{\calx} + \|\cdot\|_{\caly}$,
and
$\calx + \caly$ is equipped with norm
$\|u\|_{\calx + \caly} =\inf\{\|u_1\|_\calx + \|u_2\|_\caly: u=u_1 + u_2, u_1 \in \calx, u_2 \in \caly\}$.

A pair $(p,q)$ is called a Strichartz pair,
if $\frac 2 q = d(\frac 12 - \frac 1p)$,
$(p,q)\in [2,\9]\times[2,\9]$
and $(d,p,q)\not = (2,\9,2)$.
For any interval $I\subseteq\bbr^+$,
define the Strichartz spaces by
\begin{align*}
    S^0(I):= \bigcap\limits_{(p,q):Strichartz\ pair} L^q(I; L^p), \ \
    N^0(I):= \bigcup\limits_{(p,q):Strichartz\ pair} L^{q'}(I; L^{p'}).
\end{align*}
Similarly, let
$S^1(I) = \{u\in \mathscr{S}': \|u\|_{S^0(I)} + \|\na u\|_{S^0(I)} <\9\}$,
and
$N^1(I) = \{u\in \mathscr{S}': \|u\|_{N^0(I)} + \|\na u\|_{N^0(I)} <\9\}$.
In particular,
the Strichartz spaces $V(I)= L^{2+\frac 4d}(I; L^{2+\frac 4d})$,
$W(I) = L^{\frac{2(d+2)}{d-2}}(I; L^{\frac{2d(d+2)}{d^2+4}})$
and $\bbw(I) = L^{\frac{2(d+2)}{d-2}}(I; W^{1, \frac{2d(d+2)}{d^2+4}})$
will be frequently used
in the mass and energy critical spaces.

We use the exotic Strichartz spaces
$X^0(I)$, $\mathbb{X}(I)$ and $\bby(I)$
with the norms
\begin{align*}
   \|u\|_{X^0(I)}
   = \|u\|_{L^{\frac{d(d+2)}{2(d-2)}}(I;L^{\frac{2d^2(d+2)}{(d+4)(d-2)^2}})},&  \ \
   \|u\|_{\bbx(I)}
   := \| \<\na \>^{\frac{4}{d+2}} u\|_{L^{\frac{d(d+2)}{2(d-2)}}(I;L^{\frac{2d^2(d+2)}{d^3-4d+16}})}, \\
 \|u\|_{\bby(I)}
   := &\||\<\na \>^{\frac{4}{d+2}} u \|_{L^{{\frac{d}{2}}}(I; L^{\frac{2d^2(d+2)}{d^3+4d^2+4d-16}})},
\end{align*}
which are the inhomogeneous versions of exotic Strichartz spaces in \cite{KV13}.

We also use the local smoothing spaces defined by,
for $\a, \beta \in \mathbb{R}$,
$$L^2(I;H^\a_{\beta})=\{u\in \mathscr{S}': \int_{I} \int \<x\>^{2\beta}|\<\na\>^{\a} u(t,x)|^2  dxdt <\9 \}.$$

Throughout this paper,
we use $C(\cdots)$ for various constants that may
change from line to line.

\section{Formulations of main results}
Let us start with the definition of solutions to \eqref{equa-x}.

\begin{definition}\label{def-x}
Fix $T>0$. An $L^2$-(resp., $H^1$-)solution to \eqref{equa-x} is an $L^2$-(resp., $H^1$-)valued continuous
$(\mathscr{F}_t)$-adapted process $X=X(t)$, $t\in[0,T],$ such that $|X|^\a\in L^1([0,T], H^{-1})$
and it satisfies $\bbp$-a.s.,
\begin{align}\label{equa-x'}
   X(t) =&  X_0 -  \int^t_0 (i\Delta X(s) +\mu X(s) + \lambda i F(X(s)) )ds \nonumber \\
         & + \sum\limits_{k=1}^N\int ^t_0 X(s)G_k(s) d \beta_k(s),\ \forall t\in [0,T],
\end{align}
as an It\^o equation in $H^{-2}$ (resp. $H^{-1}$).
\end{definition}

We assume the asymptotically flat condition
as in \cite{BRZ14,BRZ16,HRZ18}.
\begin{itemize}
  \item[{\rm}(H0)]
  For each $1\leq k\leq N$,
  $G_k(t,x)= g_k(t)\phi_k(x)$,
  $g_k$ are real-valued predictable processes,
  $g_k \in L^\9(\Omega\times [0,T])$,
  $0<T<\9$,
  and
  $\phi_k \in C^\9(\bbr^d,\mathbb{C})$ satisfying that for any muti-index $\g$,
  $\g\not =0$,
  \begin{align} \label{asymflat}
     \limsup\limits_{|x|\to \9} |x|^2 |\partial_x^\g \phi_k(x)| = 0.
  \end{align}
\end{itemize}

\begin{remark}
The condition \eqref{asymflat}  is
slightly stronger than $(1.3)$ in \cite{HRZ18},
mainly for the convenience to perform pseudo-differential calculus.
Moreover, one can weaken the smoothness condition to that $\phi_k \in C^n$
for $n$ large enough.
\end{remark}

We have the local well-posedness results
in the mass- and energy-critical cases.

\begin{theorem} \label{Thm-LWP} ({\it Local well-posedness})

Consider \eqref{equa-x} in the mass-(resp., energy-)critical case,
i.e., $\a=1+4/d$, $d\geq 1$
(resp., $\a=1+4/(d-2)$, $d\geq 3$).
Assume $(H0)$.
Then, for each $X_0 \in L^2$ (resp. $X_0\in H^1$),
there exits a unique $L^2$-(resp., $H^1$-)solution $X$ to \eqref{equa-x} on $[0,\tau^*)$,
where the maximal existing time  $\tau^*$ is an $\{\mathcal{F}_t\}$-stopping time,
such that $\mathbb{P}$-a.s. for any $t\in (0,\tau^*)$
and any Strichartz pair $(\rho, \g)$,
\begin{align}
      X|_{[0,t]} & \in C([0,t]; L^2) \cap L^\g(0,t;L^\rho)  \label{esti-X-LpLq}\\
  (resp.,\  X|_{[0,t]} &\in C([0,t];H^1) \cap L^\g(0,t;W^{1,\rho}) ).  \label{esti-X-LpW1q}
\end{align}
Moreover,
$X$ exists globally $\bbp$-a.s. if for any $0<T<\9$,
\begin{align}
     & \|X(\omega)\|_{L^{2+\frac 4 d}(0,\tau^* \wedge T; L^{2+\frac 4d})} < \9,\ \ \bbp-a.s.  \label{gloexist-L2} \\
    (resp.,\  & \|X(\omega)\|_{L^{\frac{2(d+2)}{d-2}}(0,\tau^*\wedge T; L^{\frac{2(d+2)}{d-2}})} < \9,\ \ \bbp-a.s.). \label{gloexist-H1}
\end{align}
\end{theorem}

The proofs are similar to
those of \cite[Proposition $5.1$]{BRZ14} and \cite[Theorem $2.1$]{BRZ16},
and the last assertion concerning the global existence follows from the
blow-up alternative results as in \cite{BRZ14,BRZ16}.

The main result of this paper is formulated below,
concerning the global existence and uniqueness of solutions to \eqref{equa-x} in the critical cases.

\begin{theorem} \label{Thm-GWP} (Global Well-Posedness)

$(i)$ Consider \eqref{equa-x} in the defocusing mass-critical case,
i.e., $\lbb=-1$, $\a = 1+ 4/d$, $d\geq 1$. Assume $(H0)$.
Then,
for each $X_0\in L^2$ and  $0<T<\9$,
there exists a unique $L^2$-solution to \eqref{equa-x} on $[0,T]$,
satisfying that for any $p\geq 1$,
\begin{align} \label{thm-L2-L2}
   \bbe \|X\|^p_{C([0,T];L^2)} \leq C(p,T) <\9,
\end{align}
and for any Strichartz pair $(\g,\rho)$,
\begin{align} \label{thm-L2-Lpq}
    X\in L^\g(0,T; L^\rho) \cap L^2(0,T; H^\frac12_{-1}),\ \ \mathbb{P}-a.s..
\end{align}

$(ii)$ Consider \eqref{equa-x} in the defocusing energy-critical case,
i.e., $\lbb=-1$, $\a = 1+ 4/(d-2)$, $d\geq 3$. Assume $(H0)$.
In the high dimensional case where $d>6$,
assume additionally that  for each $0<t<T$,
\begin{align} \label{bdd-E-assum-d6}
    E_T:= \sup\limits_{0\leq t<\tau^* \wedge T} |X(t)|_{H^1} \leq C(T) <\9,\ \ \bbp-a.s..
\end{align}

Then,
for each $X_0\in H^1$ and  $0<T<\9$,
there exists a unique $H^1$-solution to \eqref{equa-x} on $[0,T]$,
satisfying that for any $p\geq 1$,
\begin{align} \label{thm-H1-H1}
   \bbe \|X\|^p_{C([0,T];H^1)} + \bbe \|X\|^p_{L^{\frac{2d}{d-2}}(0,T;L^{\frac{2d}{d-2}})} \leq C(p,T)<\9,
\end{align}
and for any Strichartz pair $(\g,\rho)$,
\begin{align}  \label{thm-H1-LpWq}
    X\in L^\g(0,T; W^{1,\rho}) \cap L^2(0,T; H^\frac 32_{-1}),\ \ \mathbb{P}-a.s..
\end{align}
\end{theorem}

\begin{remark}
We also have the stability results  in both mass- and energy-critical cases,
see Theorems \ref{Thm-Sta-L2}, \ref{Thm-Sta-H1-dlow} and \ref{Thm-Sta-H1} below.
\end{remark}

\begin{remark} \label{Rem-globdd-E-d6}
It is possible to obtain the pathwise bound \eqref{bdd-E-assum-d6}
by using the It\^o formula of  Hamiltonian \eqref{Ito-H} below.
The formula \eqref{Ito-H} can be derived directly by a formal computation,
however, the rigorous proof in the high dimensional case where $d>6$ is technically unclear.
See also Remark \ref{Rem-globdd-E-d6-proof} below.
\end{remark}

\begin{remark}  \label{Rem-Com-FANXU}
We would like to mention that,
the global well-posedness of  stochastic nonlinear Schr\"odinger equations
has been recently proved in \cite{FX18.1,FX18.2} for the mass-critical case
for dimension $d=1$
in the {\it conservative case},
under a different spatial decay assumption on the noise.
The results in \cite{FX18.1,FX18.2} also hold in the case
where one has a uniform pathwise control of mass 
(see \cite[Remark 1.7]{FX18.1}).
In Theorem \ref{Thm-GWP} above,
we prove the global well-posedenss of \eqref{equa-x}
in the mass-critical case for all dimensions $d\geq 1$.
In addition,
Theorem \ref{Thm-GWP} also applies  to the
{\it non-conservative case},
which is important in the physical context \cite{BG09}.
Furthermore,
Theorem \ref{Thm-GWP} proves the global well-posedness
(resp. conditional global well-posedness)
in the energy-critical case for dimensions $3\leq d\leq 6$
(resp. $d>6$),
which is not discussed in \cite{FX18.1,FX18.2}.
Below we also prove scattering
and the Stroock-Varadhan type support theorem for \eqref{equa-x},
see Theorems \ref{Thm-Sca} and \ref{Thm-Supp}.
\end{remark}

We can also enhance the estimates \eqref{thm-L2-Lpq} and \eqref{thm-H1-LpWq}
to the whole time regime,
provided that $g_k\in L^2(\bbr^+)$, $1\leq k\leq N$, a.s.. Namely, we have

\begin{theorem} \label{Thm-S0S1-Global}
Consider the situations in Theorem \ref{Thm-GWP} $(i)$ (resp. $(ii)$).
Assume additionally that $g_k\in L^2(\bbr^+)$, $1\leq k\leq N$, a.s..
Then, for each $X_0\in L^2$ (resp. $X_0\in H^1$),
the solution $X$ to \eqref{equa-x}
satisfies that  for any Strichartz pair $(\rho, \g)$,
\begin{align}
    &X\in L^\g(\bbr^+; L^{\rho}) \cap L^2(\bbr^+; H^\frac 32_{-1}), \ \ \bbp-a.s.   \label{globbd-L2-Lpq} \\
   (resp.\ & X\in L^\g(\bbr^+; W^{1,\rho}) \cap L^2(\bbr^+; H^\frac 32_{-1})\ \ \bbp-a.s..).  \label{globdd-H1-LpWq}
\end{align}
\end{theorem}

Next, we study the scattering behaviour of global solutions to \eqref{equa-x} at infinity.

Besides in  $L^2$ and $H^1$,
we also work with the pseudo-conformal space, i.e.,
$\Sigma:=\{f\in H^1: |\cdot|f(\cdot) \in L^2\}$,
in which
we assume that the time functions $g_k$ in $(H0)$
have appropriate integrability and
decay speed at infinity as in  \cite{HRZ18}.
\begin{itemize}
\item[{\rm(H1)}] For each $1\leq k\leq N$,
\begin{align} \label{AF-3}
    \limsup\limits_{|x|\to \9} |x|^3 |\partial_x^\g \phi_k(x)| = 0,\ \ 1\leq |\g|\leq 3,
\end{align}
$esssup_{\Omega} \int_0^\9(1+t^4)g_k^2(t) dt <\9$,
$1\leq k\leq N$,
and for $\bbp$-a.e. $\omega \in \Omega$,
\begin{align} \label{ILog}
     \lim\limits_{t\nearrow 1} (1-t)^{-3} \(\int_{\frac{t}{1-t}}^\9 g_k^2 (\omega, s)ds \ln \ln \({\int_{\frac{t}{1-t}}^\9 g_k^2 (\omega, s)ds}\)^{-1} \)^{\frac 12} =0.
\end{align}
\end{itemize}

\begin{remark}
As mentioned in \cite[Remark 1.4]{HRZ18},
the $L^\9(\Omega)$ condition on $\int_0^\9 (1+t^4)g_k^2(t)dt$
can be weakened by some suitable exponential integrability.
\end{remark}

In order to formulate the scattering results,
we shall use the rescaling function
\begin{align} \label{vf*}
   \vf_*(t)=-\sum_{k=1}^N \int_t^\9 G_k(s) d\beta_k(s) + \frac 12 \sum_{k=1}^N \int_t^\9 \(|G_k(s)|^2+G^2_k(s) \)ds,
\end{align}
Note that,
$\vf_* \in C(\bbr^+; W^{1,\9})$
if $g_k\in L^2(\bbr^+)$, $1\leq k\leq N$, a.s..
Then, letting
\begin{align} \label{z*}
   z_*(t):= e^{-\vf_*(t)} X(t),
\end{align}
we have
\begin{align} \label{equa-RNLS-Sca}
     i \partial_t z_*  =   e^{-\vf_*(t)}\Delta (e^{\vf_*(t)}z_*) -  e^{(\a-1) \Re \vf_*(t)} F(z_*),
\end{align}
with $ z_*(0) = X_0$. Here,
\begin{align} \label{Op-A*}
      e^{-\vf_*(t)}\Delta (e^{\vf_*(t)}z_*) =   (\Delta + b_*(t) \cdot \na + c_*(t))z_*
\end{align}
with the coefficients of lower order perturbations
\begin{align}
     b_*(t) =& -2 \sum\limits_{k=1}^N \int_t^\9 \na G_k(s) d\beta_k(s) + 2  \int_t^\9 \na \wh{\mu}(s)ds, \label{b*}\\
     c_*(t) =& \sum\limits_{j=1}^N \( \sum\limits_{k=1}^N  \int_t^\9 \partial_j G_k(s) d\beta_k(s) - \int_t^\9 \partial_j\wh{\mu}(s)ds\)^2 \nonumber \\
               &-\sum\limits_{k=1}^N \int_t^\9 \Delta G_k(s) d\beta_k(s)  + \int_t^\9 \Delta \wh{\mu}(s)ds, \label{c*}
\end{align}
and
\begin{align} \label{mu}
    \wh{\mu}(s,x) = \frac 12 \sum\limits_{k=1}^{N} (|G_k(s,x)|^2 + G_k(s,x)^2)
             = \sum\limits_{k=1}^N (\Re G_k) G_k (s,x).
\end{align}

It is also convenient to use the notations $U_*(t,s)$ (resp. $U(t,s)$),
$s,t\geq 0$, for  the evolution operators
corresponding to the random equation \eqref{equa-RNLS-Sca}
(resp. \eqref{equa-RNLS} with $\sigma=0$)
in the homogeneous case $F\equiv 0$.

We are now ready to state the scattering result.
\begin{theorem} \label{Thm-Sca} ({\it Scattering})

$(i)$
Consider the defocusing mass-critical case,
i.e., $\lbb=-1$, $\a = 1+ 4/d$, $d\geq 1$.
Assume $(H0)$ and that $g_k\in L^2(\bbr^+)$, $1\leq k\leq N$, a.s..
Then, for each $X_0\in L^2$,
the global $L^2$-solution $X$ to \eqref{equa-x} scatters at infinity, i.e.,
$\bbp$-a.s. there exist  $v_+, u_+\in L^2$ such that
\begin{align} \label{Sca-L2.1}
      e^{it\Delta}e^{-\vf_*(t)} X(t) \to  v_+,\ \ in\ L^2,\ as\ t\to \9,
\end{align}
and
\begin{align} \label{Sca-L2.2}
   U_*(0,t) e^{-\vf_*(t)} X(t) \to u_+, \ \ in\ L^2, \ as\ t\to \9.
\end{align}

$(ii)$
Consider the defocusing energy-critical case,
i.e., $\lbb=-1$, $\a = 1+ 4/(d-2)$, $d\geq 3$.
Assume $(H0)$ and that $g_k\in L^2(\bbr^+)$, $1\leq k\leq N$, a.s..
In the high dimensional case where $d>6$,
assume additionally  that
\begin{align} \label{globdd-E-assum-d6}
    E_\9: = \sup\limits_{0\leq t<\9}  |X(t)|_{H^1} \leq C <\9,\ \ a.s..
\end{align}

Then, for each $X_0\in H^1$,
the global $H^1$-solution satisfies the asymptotics \eqref{Sca-L2.1} and \eqref{Sca-L2.2}
with $H^1$ replacing $L^2$.

$(iii)$
Consider the situations as in the defocusing energy-critical case in $(ii)$,
$d\geq 3$.
Then, for each $X_0\in \Sigma$,
the asymptotic \eqref{Sca-L2.1} holds with $\Sigma$ replacing $L^2$.
\end{theorem}

\begin{remark}
Unlike in the deterministic case,
the scattering behavior of stochastic solutions to \eqref{equa-x}
is closely related to the rescaling function $e^{-\vf_*}$,
which, actually, encodes the information of noise in \eqref{equa-x}.
We would like also to mention that,
the rescaling function here is different from that in the proof of
global well-posedness in Theorem \ref{Thm-GWP}
(see \eqref{vf} below).
\end{remark}

\begin{remark}
As in the $H^1$ case in Theorem \ref{Thm-GWP} $(ii)$,
it is possible to obtain the global pathwise bound \eqref{globdd-E-assum-d6}
from the Hamiltonian \eqref{Ito-H} below,
by using similar arguments
as in the proof of \cite[(1.7)]{HRZ18}.
However, the rigorous derivation of \eqref{Ito-H} is technically unclear.
\end{remark}

Next,
we characterize the topological support of the law of global solutions to \eqref{equa-x},
in both mass-critical and energy-critical cases.

Support theorem for diffusions was initiated in the seminal papers
\cite{SV72, SV72.2}
and has been extensively studied in literature.
We refer to \cite{G05}
and \cite{G07} for stochastic nonlinear Schr\"odinger equations
with additive noise and with fractional noise, respectively.
See also \cite{MS94,MS94.2} and references therein.

Let $\mathscr{H}$ denote the Cameron-Martin space associated with
the Brownian motions $\beta=(\beta_1,\cdots,\beta_N)$,
i.e., $\mathscr{H}=\{h \in H^1(0,T; \bbr^N): h(0)=0\}$.
For any $h=(h_1,\cdots,h_N)\in \mathscr{H}$,
let $X(\beta+h)$ be the solution to \eqref{equa-x} with the driven process $\beta+h$ replacing
the Brownian motion $\beta$.
Moreover, let $S(h)$ denote the solution to the controlled equation below
\begin{align} \label{equa-Sh}
   i dS(h) &= \Delta S(h) dt + \lbb  F(S(h)) dt - i \wh{\mu} S(h) dt + iS(h) G_k \dot{h}_k dt, \\
    S(h)(0) &= X_0,  \nonumber
\end{align}
where $\wh{\mu}$ is as in \eqref{mu},
and $\dot{h}_k$ is the derivative of $h_k$.
We also use the notation ${\rm supp} (\bbp \circ X^{-1})$
for the topological support of  law of solutions to \eqref{equa-x}.

\begin{theorem} \label{Thm-Supp} (Support Theorem)

$(i)$ Consider the defocusing mass-critical case,
i.e., $\lbb=-1$, $\a = 1+ 4/d$, $d\geq 1$.
Assume $(H0)$ and that
$g_k$ are deterministic and continuous, $1\leq k\leq N$.
Let $X$ be the global $L^2$-solution to \eqref{equa-x}
corresponding to $X(0)=X_0\in L^2$.

Then,
the support  ${\rm supp} (\bbp \circ X^{-1})$
in the spaces $S^0(0,T)$ and $L^2(0,T; H^{\frac 12}_{-1})$
is the closure of the set $\{S(h), h\in \mathscr{H}\}$.

$(ii)$ Consider the defocusing energy-critical case,
i.e., $\lbb=-1$, $\a = 1+ 4/(d-2)$, $3\leq d\leq 6$.
Assume $(H0)$ and that $g_k$ are deterministic and continuous, $1\leq k\leq N$.
Let $X$ be the global $H^1$-solution to \eqref{equa-x}
with  $X(0)=X_0\in H^1$.

Then,
the support  ${\rm supp} (\bbp \circ X^{-1})$
in the spaces $S^1(0,T)$ and $L^2(0,T; H^{\frac 32}_{-1})$
is the closure of  the set $\{S(h), h\in \mathscr{H}\}$.
\end{theorem}

\begin{remark}
Theorem \ref{Thm-Supp} applies in particular to the stochastic nonlinear Schr\"odinger equations
in \cite{BRZ14,BRZ16}, where $g_k \equiv 1$, $1\leq k\leq N$.
\end{remark}

\begin{remark}
We also expect the support theorem to hold in high dimensions $d>6$ in the energy-critical case,
yet the present stability result Theorem \ref{Thm-Sta-H1} can not help us,
due to the smallness condition on the time function $g$ in \eqref{Sta-H1-ve} below.
\end{remark}

\begin{remark}
Equation \eqref{equa-Sh} can be viewed as
a subcritical (linear) perturbation of the  nonlinear Schr\"odinger equation \eqref{equa-NLS} below.
This observation helps to obtain the global well-posedness of \eqref{equa-Sh}
by using the stability results in \cite{KV13, TV05, TVZ07}.
So,
if the support theorem is assumed {\it a-priori} to hold,
then, intuitively, the stochastic solution $X$ itself should also exist globally.
This viewpoint, actually, offers an intuition for the
global well-posedness of  \eqref{equa-x} in critical cases.
\end{remark}

The proof of Theorem \ref{Thm-GWP}
is mainly based on a new application of rescaling approach and the theory of stability.

In order to prove the global existence of solutions to \eqref{equa-x},
in view of Theorem \ref{Thm-LWP},
we only need to obtain the global bounds of the
$L^{2+\frac 4 d}(0,\tau^*; L^{2+\frac 4d})$- and
$L^{\frac{2(d+2)}{d-2}}(0,\tau^*; L^{\frac{2(d+2)}{d-2}})$-norms
of solutions
in the critical cases.
Such estimates were obtained in the deterministic case by using
the energy induction method
or the concentration-compact method,
combined with
the conservation laws
(e.g., of the mass and  Hamiltonian)
and interaction Morawetz estimates.
However,
the presence of Brownian motions in \eqref{equa-x} destroys the conservation laws,
the related It\^o formulas  actually
consist of several stochastic integrals
(see \eqref{Ito-L2}, \eqref{Ito-H} below),
which make it quite hard to obtain  estimates as in the deterministic case.

Proceeding differently,
we perform a series of rescaling transformations on
a random partition (depending on the growth of  noise) of any bounded time interval.
On each small time piece,
we compare the resulting random equation with
the standard nonlinear Schr\"odinger equation with the same initial data,
by using the stability results
(see Theorems \ref{Thm-Sta-L2},
\ref{Thm-Sta-H1-dlow} and \ref{Thm-Sta-H1} below)
and the work of \cite{CKSTT08,RV07,V07} and \cite{D12,D16.1,D16.2}
in the deterministic defocusing mass- and energy-critical cases, respectively
(see Theorems \ref{Thm-L2GWP-Det} and \ref{Thm-H1GWP-Det} below).
Then,
by virtue of the global pathwise bounds of
mass and energy in the defocusing case,
we are able to put together all  finitely many bounds
in the previous step to
obtain the desirable global bounds.

For the reader's convenience,
let us explain more precisely the procedure above
on a random time interval $[\sigma, \sigma+\tau]$,
where $\sigma$ and $\sigma+\tau$ are $(\mathscr{F}_t)$-stopping times.

We use the rescaling transformation
\begin{align} \label{vsigma}
     v_\sigma(t) := e^{-\vf_\sigma(t)} X(\sigma+t),\ \ t\in[0,\tau],
\end{align}
where
\begin{align} \label{vf}
   \vf_\sigma(t,x)
   := \int_\sigma^{\sigma+ t} G_k(s,x) d\beta_k(s)
     -  \int_\sigma^{\sigma+t} \wh{\mu}(s,x) ds
\end{align}
with $\wh{\mu}$ as in \eqref{mu}.

The rescaling transformation
can be regarded as
a Doss-Sussman type transformation in Hilbert space.
See, e.g.,  \cite{BR15}
for the applications of rescaling approach
to general stochastic partial differential equations with coercive structure.
See \cite{BRZ18} for the application to optimal bilinear control problems,
see also \cite{BRZ17, Z17}
for  other quite general stochastic dispersive equations.

The nice feature is that it reveals the structure of stochastic equation \eqref{equa-x}
by reducing to the random equation with lower order perturbations below
\begin{align} \label{equa-RNLS}
   i\p_t v_\sigma
  =& e^{-\vf_\sigma} \Delta (e^{\vf_\sigma} v_\sigma) - e^{(\a-1)\Re \vf_\sigma}  F(v_\sigma), \\
  v(0)=& X(\sigma), \nonumber
\end{align}
where
\begin{align} \label{Op-A}
     e^{-\vf_\sigma}(e^{\vf_\sigma} v_\sigma)
         =  (\Delta + b_\sigma(t) \cdot \na + c_\sigma(t))v_\sigma,
\end{align}
and the coefficients
\begin{align}
   b_\sigma(t)
   =& 2 \na \vf_\sigma(t)
   = 2 \sum\limits_{k=1}^N \int_\sigma^{\sigma+t} \na G_k(s) d\beta_k(s)
                  - 2 \int_\sigma^{\sigma+t} \na \wh{\mu}(s) ds, \label{b} \\
   c_\sigma(t)
   =&\Delta \vf_\sigma + \sum\limits_{j=1}^N (\p_j \vf_\sigma)^2 \nonumber \\
   =& \sum\limits_{j=1}^N
           \( \sum\limits_{k=1}^N \int_\sigma^{\sigma+t} \p_j G_k(s) d\beta_k(s)
             -\int_\sigma^{\sigma+t} \p_j \wh{\mu}(s) ds \)^2 \nonumber \\
         & + \sum\limits_{k=1}^N \int_\sigma^{\sigma+t} \Delta G_k(s) d\beta_k(s)
           - \int_\sigma^{\sigma+t} \Delta \wh{\mu}(s)ds.  \label{c}
\end{align}

The result below connects equations \eqref{equa-x} and \eqref{equa-RNLS},
which generalizes the case where $\sigma\equiv 0$ in \cite{BRZ14}-\cite{BRZ18}.
The proof is postponed to the Appendix.

\begin{theorem} \label{Thm-Rescale-sigma}
Consider the situations in Theorem \ref{Thm-LWP}.
Let $X$ be the $L^2$-(resp. $H^1$-)solution to \eqref{equa-x} on $[0,\tau^*)$  with $X(0)=X_0 \in L^2$,
where $\tau^*$ is the maximal existing time.
Let $v_\sigma$ be as in \eqref{vsigma},
where
$\sigma$ is any $(\mathscr{F}_t)$-stopping time
satisfying  $0\leq \sigma <\tau^*$.
Then,
$v_\sigma$ satisfies \eqref{equa-RNLS}
on $[0,\tau^*-\sigma)$ in the space $H^{-2}$ (resp, $H^{-1}$) almost surely.
\end{theorem}

The key obervation here is,
that the amplitude of lower order perturbations
depends only on the trajectories of noises $\{\beta_k(t), \sigma\leq t\leq \sigma+\tau\}$.
This fact inspires us to view Equation \eqref{equa-RNLS},
if the random interval is short enough,
as a small perturbation of the  nonlinear Schr\"odinger equation
\begin{align} \label{equa-NLS}
  i\p_t \wt{u}
  =&  \Delta \wt{u} - F(\wt{u}), \\
  \wt{u}(0)=& v_\sigma(0)= X(\sigma). \nonumber
\end{align}

Now,
it becomes clear that a stability-type result will fulfill the comparison procedure above.
It should be mentioned that,
because of the lower order perturbations,
we need to prove stability results for the equation of similar form as in \eqref{equa-RNLS}.
For this reason, we reformulate  \eqref{equa-NLS}  as follows
\begin{align} \label{equa-NLS*}
    i \p_t \wt{u}
    =  e^{-\vf_\sigma}\Delta(e^{\vf_\sigma}\wt{u})
      - e^{(\a-1)\Re \vf_{\sigma}}  F(\wt{u}) + e
\end{align}
with the error term
\begin{align} \label{Error-NLS}
   e = -(b_\sigma(t) \cdot \na  + c_\sigma(t))\wt{u}
       - (1-e^{(\a-1)\Re \vf_\sigma}) F(\wt{u}),
\end{align}
where the coefficients $b_\sigma$, $c_\sigma$ are as in \eqref{b} and \eqref{c}, respectively.

The proof of  stability results in Theorems \ref{Thm-Sta-L2},
\ref{Thm-Sta-H1-dlow} and \ref{Thm-Sta-H1} below
is mainly inspired by the work \cite{KV13,TV05,TVZ07}.
However,
it relies heavily on Strichartz estimates for the Laplacian with lower order perturbations
(see Theorem \ref{Thm-Stri}).
Moreover, another important role here is played by the local smoothing spaces,
which enable us to control the lower order perturbations arising from the operator
$e^{-\vf_\sigma}\Delta(e^{\vf} \cdot)$,
for which the pseudo-differential calculus is performed.

Finally,
the proof of scattering in Theorem \ref{Thm-Sca} is based on the
very recent work \cite{HRZ18},
and the proof of the Stroock-Varadhan type support theorem
(i.e., Theorem \ref{Thm-Supp})
is inspired by the work \cite{MS94}.
In both cases,
we shall construct appropriate rescaling transformations,
related to the structure of our problems.
See, e.g., \eqref{vf*},
and \eqref{res-zn-Sbetan}, \eqref{res-yn-Xbetan} below.
We also emphasize that,
again the key ingredients are the
stability results in Theorems \ref{Thm-Sta-L2},
\ref{Thm-Sta-H1-dlow} and \ref{Thm-Sta-H1}
for the Laplacian with lower order perturbations. \\

The remainder of this paper is structured as follows.
In Section \ref{Sec-Pre},
we present the preliminaries used in this paper,
including the pseudo-differential operators,
the Strichartz and local smoothing estimates and
the exotic Strichartz spaces.
Then, we prove the stability results in both mass- and energy-critical cases in Section \ref{Sec-Sta}.
Sections \ref{Sec-GWP}, \ref{Sec-Sca} and \ref{Sec-Supp}
are mainly devoted to the proof of Theorems \ref{Thm-GWP}, \ref{Thm-Sca} and \ref{Thm-Supp}, respectively.
Finally, some technical proofs are postponed to the Appendix, i.e. Section \ref{Sec-App}.

\section{Preliminaries} \label{Sec-Pre}

This section collects some preliminaries used in this paper.

\subsection{Pseudo-differential operators}

We recall some basic facts of pseudo-differential operators.
For more details see \cite{K81,T00,Z17} and references therein.

We say that
$a\in C^\9(\bbr^d \times \bbr^d)$ is a symbol of class $S^m$,
if for any multi-indices $\a,\beta \in \mathbb{N}^d$,
$|\partial^\a_\xi\partial^\beta_x a(x,\xi)| \leq C_{\a,\beta} \<\xi\>^{m-|\a|}$.
The semi-norms $|a|_{S^m}^{(l)}$ are defined by
$$|a|_{S^m}^{(l)} = \max_{|\a+\beta|\leq l} \sup_{\bbr^{2d}}
                        \{ |\partial^\a_\xi \partial^\beta_x a(x,\xi)|\<\xi\>^{-(m-|\a|)}\},\ \ l\in \mathbb{N}. $$
Let $\Psi_a$  denote the pseudo-differential operator
related to the symbol $a(x,\xi)$, i.e.,
\begin{align*}
   \Psi_a v(x)   = (2\pi)^{-d} \int e^{ix \cdot \xi} a(x,\xi) \mathscr{F}(v)(\xi) d\xi,\ \ v\in \mathscr{S}.
\end{align*}
In this case, we write $\Psi_a \in S^m$  when no confusion arises.

\begin{lemma} \label{Lem-Err}
Let $a_i\in S^{m_i}$, $i=1,2$. Then, $\Psi_{a_1} \circ \Psi_{a_2} = \Psi_a \in S^{m_1+m_2}$
with
\begin{align*}
   a(x,\xi)
   = (2\pi)^{-d} \iint e^{-iy\cdot \eta} a_1(x,\xi+\eta) a_2(x+y,\xi) dy d\eta .
\end{align*}
\end{lemma}

Note that, the commutator $i[\Psi_{a}, \Psi_{b}]:= i(\Psi_{a} \Psi_{b} -  \Psi_{b}  \Psi_{a})$
is an operator with symbol in $S^{m_1+m_2-1}$, and the principle symbol is the Poisson bracket
\begin{align*}
    H_{a}b := \{a,b\} = \sum\limits_{j=1}^d \partial_{\xi_j}a \partial_{x_j}b - \partial_{\xi_j}b \partial_{x_j}a.
\end{align*}

One can also expand the composition of two pseudo-differential operators
into any finite order
and estimate the remainder.
See Lemmas $3.1$ and $3.2$ in \cite{Z17}.

\begin{lemma} \label{Lem-L2-Bdd}
Let  $a\in S^0$, $p\in (1,\9)$.
Then, for some $C>0$ and $l \in \mathbb{N}$,
\begin{align} \label{pdo-l2}
    \|\Psi_a \|_{\mathcal{L}(L^p)} \leq C |a|_{S^0}^{(l)} .
\end{align}
\end{lemma}

\subsection{Strichartz and local smoothing estimates} \label{Subsec-Stri}

We first present the Strichartz and local smoothing estimates below.

\begin{theorem}   \label{Thm-Stri}
Let $ I=[t_0,T]\subseteq \bbr^+$.
Consider the equation
\begin{align} \label{equa-stri}
   & i\partial_t u = e^{-\Phi} \Delta(e^\Phi u) + f.
\end{align}
Here,
the function
$\Phi =\Phi(t,x)$ is continuous on $t$ for each $x\in \bbr^d$, $d\geq 1$,
and satisfies that  for each multi-index $\g$,
\begin{align} \label{psi-Stri}
     \sup\limits_{t\in I} |\p_x^\g \Phi(t,x)| \leq C(\g) \sup\limits_{t\in I}g(t)  \<x\>^{-2}
\end{align}
for some positive and continuous function $g$.
Then, for any  $u(t_0)\in L^2$
and $f\in N^0(I) + L^2{(I; H^{-\frac 12}_1)}$,
the solution $u$ to \eqref{equa-stri} satisfies
\begin{align}  \label{L2-Stri}
    \|u\|_{S^0(I) \cap L^2{(I; H^{\frac12}_{-1})}}\leq
    C_T (|u(t_0)|_{2}+\|f\|_{N^0(I) + L^2{(I; H^{-\frac 12}_1)}}).
\end{align}
Moreover, if in addition $u(t_0)\in H^1$, $d\geq 3$,
$f\in N^1(I)+ L^2{(I; H^{\frac 12}_1)}$,
then
\begin{align}  \label{H1-Stri}
    \|u\|_{S^1(I) \cap L^2{(I; H^\frac 32_{-1})}}
    \leq&
    C_T \big(|u(t_0)|_{H^1}+\|f\|_{N^1(I)+ L^2{(I; H^{\frac 12}_1)}}).
\end{align}
\end{theorem}
Below, we use the notation $C_T$ for the
constant in Strichartz estimates  above
throughout the paper.
We may assume $C_T \geq 1$ without lose of generality.

\begin{remark}
Estimates \eqref{L2-Stri} and \eqref{H1-Stri} are the so called {\it local-in-time} estimates,
in that the constant $C_T$ depends on time.
Quantitative estimates and $L^p(\Omega)$-integrability
of $C_T$ have been obtained in \cite{Z17}
for quite general stochastic dispersive equations,
including stochastic Schr\"odinger equations with variable coefficients
as well as the stochastic Airy equation.
See also \cite{MMT08} for more general situations
where  Hamiltonian flows associated to Schr\"odinger operators are trapped.
\end{remark}

{\bf Proof.}
Estimate \eqref{L2-Stri} can be proved similarly as  in \cite[Theorem $2.11$]{Z17}.
See also Remark $2.14$ in \cite{Z17}.
Actually,
the asymptotically flat condition \eqref{psi-Stri} guarantees that the
lower order perturbations arising in the operator $e^{-\Phi}\Delta( e^{\Phi} \cdot)$
can be controlled, via the G{\aa}rding inequality,
by the Poisson bracket $i[\Psi_h,\Delta]$
for some appropriate symbol $h\in S^0$
(see the proof of \cite[Theorem 4.1]{Z17}).
We refer to \cite{Z17} for more details.
See also \cite[Lemma 4.1]{BRZ14} and \cite[Lemma 2.7]{BRZ16} for the
special case where $\Phi$ is as in \eqref{vf} with $\sigma \equiv 0$.

Regarding \eqref{H1-Stri},
Applying the operator $\<\na\>$ to both sides of \eqref{equa-stri}  we get
\begin{align}  \label{equa-nau}
   i \p_t (\<\na\> u) = e^{-\Phi}\Delta(e^\Phi (\<\na\>u))
                        + [\<\na\>, b\cdot \na+c] u + \<\na\>f,
\end{align}
where the coefficients $b=\na \Phi$,
$c= \Delta \Phi+ \sum_{j=1}^d (\p_j\Phi)^2$.
We regard \eqref{equa-nau} as the equation for the unknown $\<\na\>u$.
Then,   \eqref{L2-Stri} yields
\begin{align} \label{H1-stri.0}
    \|u\|_{S^1(I) \cap L^2{(I; H^{\frac 32}_{-1})}}\leq
    C_T& \bigg(|u(t_0)|_{H^1}
         +\|[\<\na\>, b\cdot \na+c] u\|_{L^2(I; H^{-\frac 12}_1)}  \nonumber \\
       &\qquad   +\|f\|_{N^1(I) + L^2{(I; H^{\frac 12}_1)}} \bigg).
\end{align}
Note that,
for the commutator $[\<\na\>, b\cdot \na+c]$,
\begin{align*}
   \<x\>\<\na\>^{-\frac 12} [\<\na\>, b\cdot \na+c] = \Psi_p \<x\>^{-1}\<\na\>^{\frac 12},
\end{align*}
where $\Psi_p:= \<x\>\<\na\>^{-\frac 12} [\<\na\>, b\cdot \na+c] \<\na\>^{-\frac 12} \<x\>$
is a pseudo-differential operator of order $0$
with semi-norms depending on $\sup_{t\in I} g(t)$.
By Lemma \ref{Lem-L2-Bdd} and \eqref{L2-Stri},
\begin{align} \label{esti-bc-Stri}
   \|[\<\na\>, b\cdot \na+c] u\|_{L^2(I; H^{-\frac 12}_1)}
   \leq& C \sup\limits_{t\in I} g(t) \|u\|_{L^2(I; H^{\frac 12}_{-1})} \nonumber \\
   \leq& CC_T \sup\limits_{t\in I} g(t) \|f\|_{N^0(I)+L^2(I; H^{-\frac 12}_{1})}.
\end{align}
Since when $d\geq 3$,
$\<x\>^2$ is a weight of Muckenhoupt class $A_2$
(see, e.g., \cite[Lemma 2.3 (iv)]{FS97}),
by virtue of the boundedness of multiplier $m(\xi)=\<\xi\>^{-1}$
in the weighted space $L^2(\<x\>^2dx)$ (see, e.g., \cite{KW79, K80}),
we have the embedding
$H_1^{\frac 12}  \hookrightarrow H_1^{-\frac 12}$
and so $L^2(I; H^{\frac 12}_{1}) \hookrightarrow L^2(I; H^{-\frac 12}_{1})$.

Therefore, taking into account $N^1(I)  \hookrightarrow N^0(I)$
and plugging \eqref{esti-bc-Stri} into \eqref{H1-stri.0}
we obtain \eqref{H1-Stri}. The proof is complete.
\hfill $\square$ \\

It is known that {\it global-in-time} Strichartz and local smoothing estimates
(i.e., the constant $C_T$ is independent of $T$)
hold for the free Schr\"odinger group $\{e^{-it\Delta}\}$.
See, e.g., \cite{B08, KV13, MMT08} and references therein.
This is also true for
the operator $-i e^{-\Phi}\Delta(e^{\Phi} \cdot)$
when $g$ satisfies some smallness condition,
which is crucial in the study of scattering in Section \ref{Sec-Sca} below.
Precisely, we have

\begin{theorem}  \label{Thm-Stri*}
Consider the situations as in Theorem \ref{Thm-Stri}.
Assume $(H0)$.
Assume additionally that for some $T_*>0$,
$\sup_{t\geq T_*} g(t) \leq \ve$ with $\ve$ sufficiently small.

Then,
the estimates \eqref{L2-Stri} and
\eqref{H1-Stri} also hold
with some constant $C$ independent of $t_0$ and $T$,
and $u(t_0)$ can be replaced by the final datum $u(T)$.
\end{theorem}

\begin{remark}
Similar estimates were proved in \cite[Corollary 5.3]{HRZ18},
with $L^2(I;H^\frac{1}{2}_{-1})$ and $L^2(I;H^{-\frac{1}{2}}_{1})$
replaced by the local smoothing spaces $LS(I)$ and $LS'(I)$ introduced in \cite{MMT08}, respectively.
We refer to \cite{MMT08} for more general situations.
\end{remark}

{\bf Proof.}
Below we mainly consider the $L^2$ case.
The proof is  quite similar to that of \cite[Corollary 5.3 (i)]{HRZ18}.

Actually, we have from Equation \eqref{equa-stri} that
\begin{align*}
   i \p_t u = \Delta u + (b\cdot \na + c) u + f,
\end{align*}
where $b,c$ are as in the proof of Theorem \ref{Thm-Stri}.

We assume $T_*<  T$ without lose of generality.
First, on the time regime $[t_0, T_*]$,
using \eqref{L2-Stri} we have
\begin{align} \label{esti-uS0-Stri.0}
   \|u\|_{S^0(t_0, T_*) \cap L^2(t_0, T_*; H^{\frac 12}_{-1})}
   \leq C_{T_*} |u(t_0)|_2
         + C_{T_*} \|f\|_{N^0(t_0, T_*) + L^2(t_0, T_*; H^{-\frac 12}_{1})}
\end{align}
Moreover, on the regime $[T_*,T]$,
using the global-in-time Stichartz and local smoothing estimates of
the free Schr\"odinger group $\{e^{-it\Delta}\}$ we get
\begin{align} \label{esti-uS0-Stri.1}
        &\|u\|_{S^0(T_*, T) \cap L^2(T_*, T; H^{\frac 12}_{-1})}  \nonumber  \\
   \leq& C |u(T_*)|_2  + C \|(b\cdot \na + c) u\|_{L^2(T_*, T; H^{-\frac 12}_{1})}
        + C \|f\|_{N^0(T_*, T) + L^2(T_*, T; H^{-\frac 12}_{1})}  \\
   \leq& C |u(t_0)|_2  + C \|(b\cdot \na + c) u\|_{L^2(T_*, T; H^{-\frac 12}_{1})}
          + (C_{T_*}+C) \|f\|_{N^0(T_*, T) + L^2(T_*, T; H^{-\frac 12}_{1})}, \nonumber
\end{align}
where $C$ is independent of $t_0$ and $T$,
and in the last step we also used \eqref{esti-uS0-Stri.0} to bound $|u(T_*)|_2$.
Note that, similarly to \eqref{esti-bc-Stri},
\begin{align*}
   \|(b\cdot \na + c) u\|_{L^2(T_*, T; H^{-\frac 12}_{1})}
   \leq C \sup\limits_{t\in I} g(t) \|u\|_{L^2(T_*, T; H^{\frac 12}_{-1})}
   \leq C \ve  \|u\|_{L^2(T_*, T; H^{\frac 12}_{-1})}.
\end{align*}
Plugging this into \eqref{esti-uS0-Stri.1}
we obtain that for $\ve$ small enough,
\begin{align} \label{esti-uS0-Stri.2}
   \|u\|_{S^0(T_*, T) \cap L^2(T_*, T; H^{\frac 12}_{-1})}
   \leq  C |u(t_0)|_2
         + C \|f\|_{N^0(T_*, T) + L^2(T_*, T; H^{-\frac 12}_{1})},
\end{align}
Now, combining \eqref{esti-uS0-Stri.0} and \eqref{esti-uS0-Stri.2} together
we obtain \eqref{L2-Stri} with the constant uniformly bounded on the whole time regime $\bbr^+$.

The $H^1$ case can be proved similarly.
Moreover, one can use similar arguments as in the proof of \cite[Corollary 5.3 (iii)]{HRZ18}
to replace $u(t_0)$ in \eqref{L2-Stri} and \eqref{H1-Stri} with the final datum $u(T)$.
The proof is complete.
\hfill $\square$ \\

In the end of this subsection,
we collect some estimates in the Strichartz space $V(I)$, $W(I)$ and $\bbw(I)$,
where $I$ is any interval in $\bbr^+$.
These estimates will be frequently used  throughout this paper.
Precisely, we have
\begin{align}  \label{ineq-V}
   & \| |u|^{\frac 4d} v\|_{L^{\frac{2+4}{d+4}}(I\times \bbr^d)}
   \leq \|u\|_{V(I)}^{\frac 4d} \|v\|_{V(I)}, \\
   & \||u|^{\frac{4}{d-2}}  v\|_{L^2(I; L^{\frac{2d}{d+2}})}
    \leq \|u\|^{\frac{4}{d-2}}_{L^{\frac{2(d+2)}{d-2}}(I; L^{\frac{2(d+2)}{d-2}})}
          \|v\|_{W(I)}
    \leq   \|u\|^{\frac{4}{d-2}}_{\bbw(I)}
          \|v\|_{W(I)},  \label{ineq-W.2}
\end{align}
and if $3\leq d\leq 6$,
\begin{align}
   \||u|^{\frac{6-d}{d-2}}  v \na w\|_{L^2(I; L^{\frac{2d}{d+2}})}
   \leq \|u\|^{\frac{6-d}{d-2}}_{\bbw(I)} \|v\|_\bbw
          \|w\|_{\bbw(I)}.   \label{ineq-W.3}
\end{align}
Estimates \eqref{ineq-V}-\eqref{ineq-W.3} can be proved by using the H\"older inequality
and the Sobolev embedding
\begin{align} \label{Sob-W1-Lpq}
   \bbw{(I)} \hookrightarrow L^{\frac{2(d+2)}{d-2}}(I; L^{\frac{2(d+2)}{d-2}}).
\end{align}

\subsection{Exotic Strichartz estimates} \label{Subsec-Exotic}

The exotic Strichartz spaces are introduced primarily
to treat the non-Lipschitzness of the derivatives of nonlinearity,
particularly for dimensions larger than six.
Actually,
for the nonlinearity $F(u)=|u|^{\frac{4}{d-2}} u$, $u\in \bbc$,
we have (see \cite[(1.3), (1.4)]{TVZ07})
\begin{align}
  & |F_z(u)| + |F_{\ol{z}}(u)|
  \leq C |u|^{\frac{4}{d-2}}, \label{Fz.1}\\
  & |F_z(u) - F_z(v)| + |F_{\ol{z}}(u) - F_{\ol{z}}(v)|
  \leq \left\{
         \begin{array}{ll}
           C|u-v|^{\frac{4}{d-2}} , & \hbox{if $d>6$;} \\
           C|u-v|(|u|^{\frac{6-d}{d-2}} + |v|^{\frac{6-d}{d-2}}), & \hbox{if $3\leq d\leq 6$.}
         \end{array}
       \right. \label{Fz.2}
\end{align}
The space $\bbx(0,\tau)$ allows to take $\frac{4}{d+2}$-derivatives of the nonlinearity,
instead of taking the full derivative.

Below We recall some important estimates in the exotic Strichartz spaces when $d\geq 3$,
which are mainly proved in \cite{KV13} in the homogenous case.
The arguments there apply also the inhomogenous case
considered in this paper.

\begin{lemma} \label{Lem-X0-bbX-S1}
For any compact time interval $I \subseteq \bbr^+$,
\begin{align}
    & \|u\|_{X^0(I)} \leq C \|u\|_{\bbx(I)} \leq C \|u\|_{S^1(I)}.  \label{X0-bbX-S1}\\
    & \|u\|_{\bbx(I)} \leq C \|u\|^{\frac{1}{d+2}}_{L^\frac{2(d+2)}{(d-2)}(I \times \bbr^d)}
                             \| u\|_{S^1(I)}^{\frac{d+1}{d+2}}
                      \leq C\|u\|^{\frac{1}{d+2}}_{\bbw(I)}
                             \| u\|_{S^1(I)}^{\frac{d+1}{d+2}}.  \label{bbX-LpS0}
\end{align}
and for some $0<c\leq 1$,
\begin{align}  \label{Lp-bbXS1}
     \|u\|_{L^\frac{2(d+2)}{(d-2)}(I\times \bbr^d)}
       \leq \|u\|^c_{\bbx} \|u\|_{S^1(I)}^{1-c}.
\end{align}
\end{lemma}
The proof is similar to that of \cite[Lemma 3.11]{KV13}.

\begin{lemma} \label{Lem-bbX-bbYLSN1}
Let $I=[t_0,T]$ be any compact interval in $\bbr^+$. We have
\begin{align}
  \|e^{-i(\cdot -t_0) \Delta }u_0\|_{\bbx(I)} \leq&  C |u_0|_{H^1}, \label{bbX-H1} \\
 \bigg \|\int_{t_0}^\cdot e^{-i(\cdot-t_0)\Delta} f(s) ds \bigg \|_{\bbx(I)}
     \leq &  C \|f\|_{\bby(I)+L^2(I;H^\frac 12_1)+N^1(I)}. \label{bbX-bbYLSN1}
\end{align}
\end{lemma}

{\bf Proof.}
Estimate \eqref{bbX-H1} follows from \eqref{X0-bbX-S1}
and the homogenous Strichartz estimates.
For \eqref{bbX-bbYLSN1},
similar arguments as in the proof of \cite[Lemma 3.10]{KV13} yield
\begin{align*}
  \bigg \|\int_{t_0}^\cdot e^{-i(\cdot-t_0)\Delta} f(s) ds \bigg \|_{\bbx(I)}
     \leq    C \|f\|_{\bby(I)}.
\end{align*}
Moreover, using \eqref{X0-bbX-S1} and Strichartz estimates we have
\begin{align*}
   \bigg \|\int_{t_0}^\cdot e^{-i(\cdot-t_0)\Delta} f(s) ds \bigg \|_{\bbx(I)}
   \leq C \bigg \|\int_{t_0}^\cdot e^{-i(\cdot-t_0)\Delta} f(s) ds \bigg \|_{S^1(I)}
   \leq C  \|f\|_{ L^2(I;H^\frac 12_1)+N^1(I)}.
\end{align*}
Combining the estimates above together we prove \eqref{bbX-bbYLSN1}.
\hfill $\square$

\begin{lemma}  \label{Lem-F-bbY-bbX}
For any compact time interval $I \subseteq \bbr^+$,
\begin{align}
       \|F(u)\|_{\bby(I)} \leq&  C \|u\|_{\bbx(I)}^{\frac{d+2}{d-2}}. \label{F-bbY-bbX}
\end{align}
Moreover,
\begin{align} \label{Fz-bbY-S1bbX}
       \|F_z(u+v)w\|_{\bby(I)}
        \leq &  C(\|u\|_{\bbx(I)}^\frac{8}{d^2-4} \|u\|_{S^1(I)}^\frac{4d}{d^2-4}
               + \|v\|_{\bbx(I)}^\frac{8}{d^2-4} \|v\|_{S^1(I)}^\frac{4d}{d^2-4} )
             \|w\|_{\bbx(I)},
\end{align}
and similar estimate also holds for $\|F_{\ol{z}}(u+v)w\|_{\bby(I)} $.
\end{lemma}
The proof is similar to that of \cite[Lemma 3.12]{KV13}.

\section{Stability} \label{Sec-Sta}

This section is devoted to the stability results in the mass and energy critical cases,
which are crucial in the proof of global well-posedness in the next section.

To begin with,
let us start with the easier mass-critical case.

\subsection{Mass-critical case} \label{Subsec-Sta-L2}
The main result of this subsection is formulated below.

\begin{theorem}  \label{Thm-Sta-L2} ({\it Mass-Critical Stability Result}).
Fix $I=[t_0,T]\subseteq \bbr^+$.
Let $v$ be the solution to
\begin{align} \label{equa-v-p}
     i\p_t v = e^{-\Phi}\Delta(e^{\Phi} v) - e^{\frac 4d \Re \Phi}F(v),
\end{align}
where $\Phi$ satisfies \eqref{psi-Stri},  $d\geq 1$,
and $\wt{v}$ solve the perturbed equation
\begin{align} \label{equa-wtv-p}
   i\p_t \wt{v} = e^{-\Phi}\Delta(e^{\Phi} \wt{v}) - e^{\frac 4d \Re \Phi}F(\wt{v}) + e
\end{align}
for some function $e$.
Assume that
\begin{align} \label{Sta-L2-V}
     \|\wt{v}\|_{C(I; L^2)} \leq M,  \ \
    |v(t_0) - \wt{v}(t_0)|_{2} \leq M', \ \
      \|\wt{v}\|_{V(I)} \leq L
\end{align}
for some positive constants $M,M'$ and $L$.
Assume also the smallness conditions
\begin{align} \label{Sta-L2-ve}
    \| U(\cdot, t_0)(v(t_0)-\wt{v}(t_0)) \|_{V(I)} \leq \ve, \ \
     \| e \|_{L^2(I; H^{-\frac 12}_1) + {N}^0(I)} \leq \ve
\end{align}
for some $0<\ve\leq \ve_*$,
where $\ve_* = \ve_*(C_T, D_T, M,M',L)>0$ is a small constant,
$C_T$ is the Strichartz constant in Theorem \ref{Thm-Stri},
$D_T = \|e^{\frac 4d \Re\Phi}\|_{C(I;L^\9)}$.
Then,
\begin{align}
   & \|v-\wt{v}\|_{V(I)} \leq C(C_T, D_T, M,M',L) \ve , \label{Sta-L2.1}\\
   & \|v-\wt{v}\|_{S^0(I) \cap L^2(I; H^\frac 12_{-1})} \leq  C(C_T, D_T, M,M',L)M', \label{Sta-L2.2} \\
   & \|v\|_{S^0(I)\cap L^2(I; H^\frac 12_{-1})} \leq C(C_T, D_T, M, M', L).  \label{Sta-L2.3}
\end{align}
We can take $\ve_*(C_T,D_T, M, M', L)$
(resp. $C(C_T,D_T, M, M', L)$)
to be decreasing (resp. nondecreasing) with respect to each argument.
\end{theorem}

Theorem \ref{Thm-Sta-L2} states that,
if the difference between two initial data
and the error term are small enough in
appropriate spaces,
the two corresponding solutions will also
stay very close to each other in the mass-critical space $V(t_0,T)$.

\begin{remark} \label{Rem-glob-Sta-L2}
Theorem \ref{Thm-Sta-L2} also holds if $\Phi$
is replaced by $\vf_*$  as in \eqref{vf*}.
In this case,
since the Strichartz constants are independent of time and $\vf_*\in L^\9(\bbr^+; L^\9)$,
the constants in \eqref{Sta-L2.1}-\eqref{Sta-L2.3}
are independent of time,
i.e., depend only $M'$, $M$ and $L$.
This fact will be important in the study of scattering in Section \ref{Sec-Sca} later.
\end{remark}

In order to prove Theorem \ref{Thm-Sta-L2},
we first prove the short-time perturbation result below.

\begin{proposition}  \label{Pro-ShortP-L2} ({\it Mass-Critical Short-time Perturbation}).
Let $I=[t_0, T]\subseteq\bbr^+$
and
$v$, $\wt{v}$ be the solutions to Equations
\eqref{equa-v-p} and \eqref{equa-wtv-p}, respectively.
Assume that,
\begin{align} \label{Short-L2-l2}
     \|\wt{v}\|_{C(I; L^2)} \leq M,  \ \
     |v(t_0) - \wt{v}(t_0)|_{2} \leq M'
\end{align}
for some positive constants $M,M'$.
Assume also the smallness conditions
\begin{align} \label{Short-L2-ve}
    \|\wt{v}\|_{V(I)} \leq \delta, \ \
    \| U(\cdot, 0)(v(t_0)-\wt{v}(t_0)) \|_{V(I)} \leq \ve, \ \
     \| e \|_{L^2(I; H^{- \frac 12}_{1}) + {N}^0(I)} \leq \ve
\end{align}
for some $0<\ve\leq \delta$
where $\delta = \delta(C_T, D_T, M,M')>0$ is a small constant,
and $C_T, D_T$ are as in Theorem \ref{Thm-Sta-L2}.
Then, we have
\begin{align}
   & \|v-\wt{v}\|_{V(I)} \leq C(C_T, D_T) \ve, \label{Short-L2.1} \\
   & \|v-\wt{v}\|_{S^0(I)\cap L^2(I; H^\frac 12_{-1})} \leq C(C_T, D_T) M', \label{Short-L2.2} \\
   & \|v\|_{S^0(I)\cap L^2(I; H^\frac 12_{-1})} \leq C(C_T, D_T) (M+M'), \label{Short-L2.3} \\
   & \|e^{\frac 4d \Re \Phi}(F(v) - F(\wt{v})) \|_{{N}^0(I)}
     \leq C(C_T, D_T) \ve.  \label{Short-L2.4}
\end{align}
\end{proposition}

{\bf Proof.}
The proof is similar to that of \cite[Lemma $3.4$]{TVZ07},
however,
based on Theorem \ref{Thm-Stri},
i.e.,
the Strichartz estimates for the Laplacian with lower order perturbations.

Let $z:= v-\wt{v}$.
In view of the equations \eqref{equa-v-p} and \eqref{equa-wtv-p},  we have
\begin{align} \label{equa-z-p}
    i\p_t z =& e^{-\Phi}\Delta(e^{\Phi} z) - e^{\frac 4d \Re \Phi} (F(z+\wt{v}) - F(\wt{v})) - e,  \\
    z(t_0)=& v(t_0) - \wt{v}(t_0),   \nonumber
\end{align}
or equivalently,
\begin{align} \label{equa-z-p*}
   z(t) = U(t,t_0)z(t_0)
          + \int_{t_0}^t U(t,s) (ie^{\frac 4d \Re \Phi} (F(z+\wt{v}) - F(\wt{v})) + i e) ds.
\end{align}

Set $S(I):= \|e^{\frac 4d \Re \Phi} (F(z+\wt{v}) - F(\wt{v})\|_{ N^0(I)}$.
By  \eqref{ineq-V} and  \eqref{Short-L2-ve},
\begin{align} \label{esti-S-L2}
   S(I) \le& \|e^{\frac 4d \Re \Phi} (F(z+\wt{v}) - F(\wt{v}) \|_{L^{\frac{2(d+2)}{d+4}}(I \times \bbr^d)}  \nonumber \\
        \leq& CD_T (\|\wt{v}\|_{V(I)}^{\frac 4d} \|z\|_{V(I)} + \|z\|_{V(I)}^{1+\frac 4d}) \nonumber \\
        \leq& C D_T (\delta^{\frac 4d} \|z\|_{V(I)} + \|z\|_{V(I)}^{1+\frac 4d}).
\end{align}
Moreover,
applying Theorem \ref{Thm-Sta-L2} to \eqref{equa-z-p}
and using \eqref{Short-L2-ve} we have
\begin{align} \label{esti-S-L2*}
     \|z\|_{V(I)}
     \leq& C_T (\|U(\cdot, t_0)z(t_0)\|_{V(I)} + S(I) + \|e\|_{N^0(I) + L^2(I; H^{- \frac 12}_{1})}) \nonumber \\
     \leq& C_T (2\ve + S(I)).
\end{align}
Then, pugging \eqref{esti-S-L2} into \eqref{esti-S-L2*} we obtain
\begin{align*}
      \|z\|_{V(I)}
      \leq C_T (2\ve + C D_T \delta^{\frac 4d} \|z\|_{V(I)}  + CD_T \|z\|_{V(I)}^{1+\frac 4d}).
\end{align*}
Thus, in view of Lemma $6.1$ in \cite{BRZ18},
for $\delta = \delta(C_T,D_T)$ small enough such that
$CC_T D_T \delta^{\frac 4d} \leq \frac 12$
and $4C_T \delta < (1-\frac 1\a)(2\a CC_TD_T)^{-\frac{1}{\a-1}}$
with $\a= 1+ \frac 4d$, we obtain
\begin{align} \label{esti-z-V}
      \|z\|_{V(I)} \leq (d+4) C_T \ve,
\end{align}
which along with \eqref{esti-S-L2} implies \eqref{Short-L2.1} and \eqref{Short-L2.4}.

Now, applying Theorem \ref{Thm-Stri} to \eqref{equa-z-p} again
and using \eqref{Short-L2-l2}, \eqref{Short-L2.4} we have
\begin{align} \label{esti-z-S0.0}
     \|z\|_{S^0(I) \cap L^2(I;H^\frac 12_{-1})}
     \leq& C_T (|z(t_0)|_2 + S(I) + \|e\|_{{N}^0(I) + L^2(I;H^{-\frac 12}_1) } ) \nonumber \\
     \leq& C_T (M' + C(C_T, D_T)\ve +\ve) \nonumber \\
     \leq& 2C_T M',
\end{align}
if $\delta= \delta(C_T,D_T,M')$ is  such that
$(C(C_T, D_T)+1)\delta  \leq M' $.
Thus, \eqref{Short-L2.2} follows.

Similarly,
by Equation \eqref{equa-wtv-p} and conditions
\eqref{Short-L2-l2} and \eqref{Short-L2-ve},
taking a even smaller
$\delta = \delta(C_T, D_T, M,M')$ such that
$D_T \delta^{1+\frac 4d} + \delta \leq M$,
we have
\begin{align} \label{esti-wtv-S0.0}
      \|\wt{v}\|_{S^0(I)\cap L^2(I;H^\frac 12_{-1})}
      \leq& C_T (|\wt{v}({t_0})|_2 + D_T\|\wt{v}\|^{1+\frac 4d}_{V(I)} + \|e\|_{N^0(I)+ L^2(I; H^{- \frac 12}_{1})}) \nonumber \\
      \leq& C_T (M+D_T \delta^{1+\frac 4d} + \delta) \nonumber \\
      \leq& 2 C_T M,
\end{align}
which along with \eqref{esti-z-S0.0} implies that
\begin{align} \label{esti-z-S0}
    \|v\|_{{S}^0(I)\cap L^2(I;H^\frac 12_{-1})}
    \leq \|z\|_{{S}^0(I)}  +  \|\wt{v}\|_{{S}^0(I)}
    \leq 2 C_T (M' + M),
\end{align}
thereby yielding \eqref{Short-L2.3}.
The proof is complete. \hfill $\square$ \\

{\bf Proof of Theorem \ref{Thm-Sta-L2}.}
First fix $\delta = \delta(C_T, D_T, M, 2C_TM')$,
where $\delta$ is as in Proposition \ref{Pro-ShortP-L2}.
We divide $[t_0,T]$ into finitely many
small pieces $I_j= [t_j,t_{j+1}]$,  $0\leq j\leq l$,
such that
$t_{l+1}=T$,
$\|\wt{v}\|_{V(t_j,t_{j+1})} =\delta$,
$0\leq j\leq l-1$,
and
$\|\wt{v}\|_{V(t_l,t_{l+1})} \leq \delta$.
Then, $l\leq (L/\delta)^{2+\frac 4d}<\9$.

Let
$C(0) = C(C_T, D_T)$,
$C(j+1) = \max\{ C(0)C_T^2 (\sum_{k=0}^j C(j)+ 2),  C(0)(1+2C_T)\}$,
$0\leq j\leq l-1$,
where $C(C_T, D_T)$ is the constant in \eqref{Short-L2.1}-\eqref{Short-L2.4}.
Choose $\ve_*= \ve_*(C_T, D_T, M,M',L)$ sufficiently small such that
\begin{align} \label{ve*-L2}
  (\sum\limits_{k=0}^l C(k) +1) \ve_* \leq M',\ \
  C_T^2 (\sum\limits_{k=0}^l C(k) +2) \ve_* \leq \delta.
\end{align}

Below we use inductive arguments to prove for any $0\leq j\leq l$,
\begin{align}
   & \|v-\wt{v}\|_{V(I_j)} \leq C(j)\ve,  \label{Sta-L2.1-proof}\\
   & \|v-\wt{v}\|_{S^0(I_j)\cap L^2(I_j;H^\frac 12_{-1})} \leq C(j) M',\label{Sta-L2.2-proof}\\
   & \|v\|_{S^0(I_j)\cap L^2(I_j;H^\frac 12_{-1})} \leq C(j)(M+M'),\label{Sta-L2.3-proof}\\
   & \|e^{\frac 4d \Re \Phi} (F(v) - F(\wt{v}))\|_{N^0(I_j)} \leq C(j) \ve. \label{Sta-L2.4-proof}
\end{align}

Proposition \ref{Pro-ShortP-L2} yields that
the estimates above hold for $j=0$.
Suppose that \eqref{Sta-L2.1-proof}-\eqref{Sta-L2.4-proof} are also valid for each $0\leq k\leq j<l$.
We shall apply Proposition \ref{Pro-ShortP-L2} to show that they also hold for the case where $j+1$ replaces $j$.

For this purpose,
by  Theorem \ref{Thm-Stri}, \eqref{Sta-L2-V} and the inductive assumptions
\begin{align*}
   |v(t_{j+1} - \wt{v}(t_{j+1})|_2
   \leq& C_T (|v(t_0) - \wt{v}(t_0)|_2 + S(t_0, t_{j+1})
             + \|e\|_{L^2(t_0,t_{j+1}; H^{-\frac 12}_1) + N^0(t_0, t_{j+1})} \\
   \leq& C_T (M' + \sum\limits_{k=0}^jC(k)\ve +\ve )
   \leq  2 C_T M',
\end{align*}
where
$S(t_0, t_{j+1})$ is as in the proof of Proposition \ref{Pro-ShortP-L2},
the last step is due to \eqref{ve*-L2}.

Moreover, by Theorem \ref{Thm-Stri},
\begin{align*}
  \|U(\cdot, t_{j+1})(v(t_{j+1}) - \wt{v}(t_{j+1}))\|_{V(I_{j+1})}
  \leq& C_T |v(t_{j+1}) - \wt{v}(t_{j+1})|_2 \\
  \leq& C_T \|v  - \wt{v} \|_{C([t_j,t_{j+1}]; L^2)}.
\end{align*}
Then, applying Theorem \ref{Thm-Stri} to \eqref{equa-z-p} again we have
\begin{align*}
  &\|U(\cdot, t_{j+1})(v(t_{j+1}) - \wt{v}(t_{j+1}))\|_{V(I_{j+1})}  \\
  \leq& C_T \|U(\cdot, t_0)(v(t_0)-\wt{v}(t_0))\|_{C([0,t_{j+1}])}
        + C^2_T \|e^{\frac 4d \Re \Phi} (F(v) - F(\wt{v})) \|_{N^0([0,t_{j+1}])}   \\
      &  + C_T^2 \|e\|_{N^0(t_0, t_{j+1}) + L^2(t_0,t_{j+1}; H^{-\frac 12}_1)} \\
  \leq&  C_T \ve + C_T^2 (\sum\limits_{k=0}^j C(k) \ve + \ve)
  \leq \delta,
\end{align*}
where the last step is again due to \eqref{ve*-L2}.

Thus, the conditions \eqref{Sta-L2-ve} and \eqref{Sta-L2-ve} of Proposition \ref{Pro-ShortP-L2} are satisfied
with $2C_TM'$ and
$C_T^2 (\sum_{k=0}^j C(k) +2 )\ve$
replacing $M'$ and $\ve$, respectively.
Proposition \ref{Pro-ShortP-L2} yields that
estimates \eqref{Sta-L2.1-proof}-\eqref{Sta-L2.4-proof} are valid with $j+1$ replacing $j$.

Therefore,  inductive arguments
yield that
\eqref{Sta-L2.1-proof}-\eqref{Sta-L2.4-proof} are valid for all $0\leq j\leq l$,
thereby proving Theorem \ref{Thm-Sta-L2}.
The proof is complete.
\hfill $\square$

\subsection{Energy-critical case} \label{Subsec-Sta-H1}
The main results of this subsection are  Theorems \ref{Thm-Sta-H1-dlow}
and \ref{Thm-Sta-H1} below.
The delicate problem here is that
the derivatives of the nonlinearity in \eqref{equa-x}
are Lipschitz when $3\leq d\leq 6$,
however, they are only H\"older continuous in high dimensions
when $d>6$.
In the latter case,
more delicate arguments involving the exotic Strichartz spaces
as well as local smoothing spaces will be used.

To begin with,
we start with  the easier case when $3\leq d\leq 6$.

\subsubsection{The case when $3\leq d\leq 6$}

In this case,
the stability result is quite similar to the previous mass-critical case.

\begin{theorem}  \label{Thm-Sta-H1-dlow} ({\it Energy-critical Stability Result when $3\leq d\leq 6$}).
Consider any bounded compact interval $I=[t_0,T]\subseteq \bbr^+$.
Let $w$ be the solution to
\begin{align} \label{equa-w-p}
     i\p_t w = e^{-\Phi}\Delta(e^{\Phi} w) - e^{\frac{4}{d-2} \Re \Phi} F(w)
\end{align}
with $\Phi$ satisfying \eqref{psi-Stri}, $3\leq d\leq 6$,
and $\wt{w}$ solve the perturbed equation
\begin{align} \label{equa-wtw-p}
   i\p_t \wt{w} = e^{-\Phi}\Delta(e^{\Phi} \wt{w}) - e^{\frac{4}{d-2}\Re\Phi} F(\wt{w})  + e
\end{align}
for some function $e$.
Assume that
\begin{align} \label{Sta-H1-V-dlow}
     \|\wt{w}\|_{C(I;H^1)} \leq E,\ \
      |w(t_0) - \wt{w}(t_0)|_{H^1} \leq E', \ \
      \|\wt{w}\|_{\bbw(I)} \leq L
\end{align}
for some positive constants $E,E'$ and $L$.
Assume also the smallness conditions
\begin{align} \label{Sta-L2-ve-dlow}
    \| U(\cdot, t_0)(w(t_0)-\wt{w}(t_0) \|_{\bbw(I)} \leq \ve, \ \
     \| e \|_{{N}^1(I) + L^2(I; H^{\frac 12}_1)} \leq \ve
\end{align}
for some $0<\ve\leq \ve_*$,
where $\ve_* = \ve_*(C_T, D'_T, E,E',L)>0$ is a small constant,
$C_T$ is the Strichartz constant in Theorem \ref{Thm-Stri}
and $D'_T = \|e^{\frac{4}{d-2} \Re \Phi} \|_{C(I; W^{1,\9})}$.
Then,
\begin{align}
   & \|w-\wt{w}\|_{\bbw(I)} \leq C(C_T, D'_T, E,E',L) \ve , \label{Sta-H1.1-dlow}\\
   & \|w-\wt{w}\|_{S^1(I) \cap L^2(I;H^\frac 32_{-1})} \leq  C(C_T, D'_T, E,E',L)E', \label{Sta-H1.2-dlow} \\
   & \|w\|_{S^1(I)\cap L^2(I;H^\frac 32_{-1})} \leq C(C_T, D'_T, E,E', L).  \label{Sta-H1.3-dlow}
\end{align}
The constants $\ve_*(C_T, D'_T, E,E', L)$ and $C(C_T, D'_T, E,E', L)$
can be taken to be decreasing and nondecreasing with respect to each argument,
respectively.
\end{theorem}

As in the mass-critical case,
Theorem \ref{Thm-Sta-H1-dlow} follows from the
short-time perturbation result below.

\begin{proposition}  \label{Pro-ShortP-H1-dlow} ({\it Energy-Critical Short-time Perturbation when $3\leq d\leq 6$}).
Let $I=[ t_0,T]$, $w$, $\wt{w}$ be as in Theorem \ref{Thm-Sta-H1-dlow},
$3\leq d\leq 6$.
Assume that
\begin{align} \label{Short-H1-h1-dlow}
     \|\wt{w}\|_{C(I; H^1)} \leq E,  \ \
     |w(t_0) - \wt{w}(t_0)|_{H^1} \leq E'
\end{align}
for some positive constants $E,E'$.
Assume also the smallness conditions
\begin{align} \label{Short-H1-ve-dlow}
    \|\wt{w}\|_{\bbw(I)} \leq \delta, \ \
    \| U(\cdot, t_0)(w(t_0)-\wt{w}(t_0)) \|_{\bbw(I)} \leq \ve, \ \
     \| e \|_{{N}^1(I) + L^2(I; H^{ \frac 12}_{1})} \leq \ve
\end{align}
for some $0<\ve\leq \delta$,
where $\delta = \delta(C_T, D'_T, E,E')>0$ is a small constant,
and $C_T, D'_T$ are as in Theorem \ref{Thm-Sta-H1-dlow}.
Then, we have
\begin{align}
   & \|w-\wt{w}\|_{\bbw(I)} \leq C(C_T, D'_T) \ve, \label{Short-H1.1-dlow} \\
   & \|w-\wt{w}\|_{S^1(I)\cap L^2(I;H^\frac 32_{-1})} \leq C(C_T, D'_T) E', \label{Short-H1.2-dlow} \\
   & \|w\|_{S^1(I)\cap L^2(I;H^\frac 32_{-1})} \leq C(C_T, D'_T) (E+E'), \label{Short-H1.3-dlow} \\
   & \|e^{\frac{4}{d-2} \Re\Phi}(F(w)-F(\wt{w})) \|_{{N}^1(I)}
     \leq C(C_T, D'_T) \ve.  \label{Short-H1.4-dlow}
\end{align}
\end{proposition}

{\bf Proof.}
Set $z:= w-\wt{w}$
and $S(I) := \|e^{\frac {4}{d-2} \Re\Phi}(F(\wt{w}+z)-F(\wt{w}))\|_{N^1(I)}$.
Then,
\begin{align*}
    S(I)  \leq
    D'_T (\|F(z+\wt{w}) - F(\wt{w}) \|_{L^2(I; L^{\frac{2d}{d+2}})}
           +  \| \na (F(z+\wt{w}) - F(\wt{w})) \|_{L^2(I; L^{\frac{2d}{d+2}})}).
\end{align*}
By \eqref{ineq-W.2}, \eqref{Fz.1} and \eqref{Short-H1-ve-dlow},
\begin{align*}
    \|F(z+\wt{w}) - F(\wt{w}) \|_{L^2(I; L^{\frac{2d}{d+2}})}
   \leq& C(\|z\|_{\bbw(I)}^{\frac{4}{d-2}} + \|\wt{w}\|_{\bbw(I)}^{\frac{4}{d-2}}) \|z\|_{\bbw(I)} \\
   \leq& C \delta^{\frac{4}{d-2}} \|z\|_{\bbw(I)} + C \|z\|^{\frac{d+2}{d-2}}_{\bbw(I)}.
\end{align*}
Moreover, since by \eqref{Fz.2} we have  (see, e.g., \cite[(3.20)]{HRZ18})
\begin{align*}
    |\na(F(z+\wt{w})-F(\wt{w}))|
   \leq C (|\na \wt{w}| |z|^{\frac{4}{d-2}}
            + |\wt{w}|^{\frac{4}{d-2}} |\na z| + |z|^{\frac{4}{d-2}} |\na z|
            + |\na \wt{w}||\wt{w}|^{\frac{6-d}{d-2}} |z|).
\end{align*}
Taking into account \eqref{ineq-W.2}, \eqref{ineq-W.3} and \eqref{Short-H1-ve-dlow} we get
\begin{align*}
    \|\na(F(z+\wt{w}) - F(\wt{w})) \|_{L^2(I; L^{\frac{2d}{d+2}})}
   \leq& C(\|\wt{w}\|_{\bbw(I)} \|z\|_{\bbw(I)}^{\frac{4}{d-2}}
          + \|\wt{w}\|_{\bbw(I)}^{\frac{4}{d-2}} \|z\|_{\bbw(I)}
          + \|z\|^{\frac{d+2}{d-2}}_{\bbw(I)} ) \\
   \leq& C( \delta \|z\|_{\bbw(I)}^{\frac{4}{d-2}}
          +  \delta ^{\frac{4}{d-2}} \|z\|_{\bbw(I)}
          + \|z\|^{\frac{d+2}{d-2}}_{\bbw(I)} ).
\end{align*}
Thus,
combining the estimates above together we obtain
\begin{align*}
   S(I) \leq CD_T'
            (  \delta ^{\frac{4}{d-2}} \|z\|_{\bbw(I)}
              +  \delta \|z\|_{\bbw(I)}^{\frac{4}{d-2}}
              + \|z\|^{\frac{d+2}{d-2}}_{\bbw(I)}).
\end{align*}
Since $1\leq \frac{4}{d-2} \leq \frac{d+2}{d-2}$ when $3\leq d\leq 6$,
$ \|z\|_{\bbw(I)}^{\frac{4}{d-2}} \leq \|z\|_{\bbw(I)}+  \|z\|^{\frac{d+2}{d-2}}_{\bbw(I)}$,
we come to
\begin{align} \label{esti-S-Short-H1-dlow}
   S(I) \leq 2CD_T'
            (  \delta \|z\|_{\bbw(I)}
              + \|z\|^{\frac{d+2}{d-2}}_{\bbw(I)}).
\end{align}

Moreover,
similarly to \eqref{equa-z-p*}, we have
\begin{align} \label{equa-z-p-H1-dlow}
   z(t) = U(t,0)z(t_0)
          + \int_{t_0}^t U(t,s)(ie^{\frac{4}{d-2}\Re\Phi}(F(z+\wt{w})- F(\wt{w})) +i e )ds
\end{align}
Applying Theorem \ref{Thm-Stri} and using \eqref{Short-H1-ve-dlow}
we have
\begin{align*}
     \|z\|_{\bbw(I)}
\leq  C_T (\|U(\cdot, t_0) z(t_0)\|_{\bbw(I)} + S(I) + \|e\|_{L^2(I; H^{ \frac 12}_{1}) + {N}^1(I)}  )
\leq C_T (2\ve + S(I)).
\end{align*}

Thus, plugging \eqref{esti-S-Short-H1-dlow} into the estimate above we obtain
\begin{align*}
     \|z\|_{\bbw(I)}
\leq&  2C_T (\ve + CD_T'
             \delta \|z\|_{\bbw(I)}
              + CD_T'\|z\|^{\frac{d+2}{d-2}}_{\bbw(I)}) .
\end{align*}
Taking $\delta=\delta(C_T, D'_T)$ very small such that
$ 2C C_T D_T' \delta \leq \frac 12$ we come to
\begin{align*}
     \|z\|_{\bbw(I)}
\leq&  4C_T \ve +  4 CC_TD_T'\|z\|^{\frac{d+2}{d-2}}_{\bbw(I)} .
\end{align*}
Then, by virtue of \cite[Lemma A.1]{BRZ18},
taking $\delta= \delta (C_T,D'_T)$ smaller
such that $ 4C_T \delta < (1-\frac 1 \a)( 4\a CC_TD_T')^{-\frac{1}{\a-1}}$
with $\a=1+\frac{4}{d-2}$
we obtain
\begin{align} \label{esti-z-Short-H1-dlow.1}
    \|z\|_{\bbw(I)} \leq\frac{ 4 \a}{\a-1} C_T \ve,
\end{align}
which together with \eqref{esti-S-Short-H1-dlow} implies \eqref{Short-H1.1-dlow} and \eqref{Short-H1.4-dlow}.

For  \eqref{Short-H1.2-dlow},
applying Theorem \ref{Thm-Stri} to \eqref{equa-z-p-H1-dlow}
and using \eqref{Short-H1-ve-dlow}, \eqref{Short-H1.4-dlow}
we have
\begin{align} \label{esti-wtw-Short-H1-dlow}
   \|z\|_{S^1(I)\cap L^2(I;H^\frac 32_{-1})}
\leq& C_T (|z(t_0)|_{H^1} + S(I) + \|e\|_{{N}^1(I) + L^2(I; H^{ \frac 12}_{1})} ) \nonumber \\
\leq& C_T (E'+ C(C_T, D'_T)\ve + \ve ),
\end{align}
which implies \eqref{Short-H1.2-dlow},
provided $\delta$ is smaller
such that $ C(C_T,D_T')\delta + \delta \leq E'$.

Similarly, by \eqref{equa-wtw-p},
\begin{align} \label{esti-z-S1-Short-H1-dlow}
   \|\wt{w}\|_{S^1(I)\cap L^2(I;H^\frac 32_{-1})}
\leq& C_T (|\wt{w}(t_0)|_{H^1} + \|e^{\frac{4}{d-2} \Re\Phi} F(\wt{w})\|_{N^1(I)} + \|e\|_{{N}^1(I)+L^2(I; H^{ \frac 12}_{1})} ) \nonumber \\
\leq& C_T (|\wt{w}(t_0)|_{H^1} + D'_T \|\wt{w}\|_{\bbw(I)}^{\frac{d+2}{d-2}} + \|e\|_{{N}^1(I)+L^2(I; H^{ \frac 12}_{1})}) \nonumber \\
\leq& C_T (E+ D'_T \delta^{\frac{d+2}{d-2}} + \ve ) \nonumber \\
\leq& 2C_T E,
\end{align}
if we take $\delta$ even smaller such that $D'_T \delta^{\frac{d+2}{d-2}} + \delta \leq E$.

Therefore, we obtain \eqref{Short-H1.3-dlow} from \eqref{Short-H1.2-dlow} and \eqref{esti-z-S1-Short-H1-dlow}
and so finish the proof.
\hfill $\square$

Once Proposition \ref{Pro-ShortP-H1-dlow} obtained,
we can use the partition arguments as in the proof of Theorem $4.1$
to prove Theorem \ref{Thm-Sta-H1-dlow}.
The details are omitted   for simplicity.

\subsubsection{The case when $d>6$}

\begin{theorem}  \label{Thm-Sta-H1} ({\it Energy-Critical Stability Result when $d>6$}).
Consider any bound compact interval $I=[t_0,T]\subseteq \bbr^+$.
Let $w, \wt{w}$ solve the equations \eqref{equa-w-p} and \eqref{equa-wtw-p}, respectively,
and $\Phi$ satisfy \eqref{psi-Stri},
$d>6$.
Assume that,
\begin{align} \label{Sta-H1-W}
   \|\wt{w}\|_{C(I; H^1)} \leq E,  \ \
      \|\wt{w}\|_{\bbw(I) \cap L^2(I;H^\frac 32_{-1})} \leq L
\end{align}
for some positive constants $E$ and $L$.
Assume also the smallness conditions
\begin{align} \label{Sta-H1-ve}
    \|g\|_{C(I;\bbr^+)} \leq \ve,\ \
    |w(t_0) - \wt{w}(t_0)|_{H^1} \leq \ve, \ \
    \| e \|_{ N^1(I)+ L^2(I; H^{\frac 12 }_1)} \leq \ve
\end{align}
for some $0<\ve\leq \ve_*$,
where
$g$ is the time function as in \eqref{psi-Stri},
$\ve_* = \ve_* (C_T, D'_T, E,L)>0$ is a small constant,
and $C_T, D'_T$ are as in Theorem \ref{Thm-Sta-H1-dlow}.
Then,  for some $c = c(C_T,D'_T,E,L)>0$,
\begin{align}
   & \|w-\wt{w}\|_{L^{\frac{2(d+2)}{d-2}}(I\times \bbr^d)} \leq  C(C_T, D'_T, E,L) \ve^c, \label{Sta-H1.1} \\
   & \|w-\wt{w}\|_{{S}^1(I)\cap L^2(I;H^\frac 32_{-1})} \leq  C(C_T, D'_T, E,L)\ve^c,  \label{Sta-H1.2} \\
   & \|w\|_{{S}^1(I)\cap L^2(I;H^\frac 32_{-1})} \leq C(C_T, D'_T, E, L).   \label{Sta-H1.3}
\end{align}
We can take the constants $\ve_*(C_T,D'_T,E,L)$ and $ C(C_T, D'_T, E, L)$
to be decreasing and nondecreasing with respect to each argument, respectively.
\end{theorem}

\begin{remark}
Unlike in the case where $3\leq d\leq 6$,
the smallness condition on $g$
is imposed in \eqref{Sta-H1-ve}   mainly to control the lower order perturbations
(see \eqref{esti-J1-H1-p} below).
One may remove this restriction on $g$,
if the estimate \eqref{bbX-bbYLSN1} still holds
with $L^2(I;H^{\frac 12}_1)$   replaced by
$L^2(0,\tau; H^{-\frac 12  + \frac{4}{d+2}}_1)$,
which, however, is unclear.
\end{remark}

We first prove
the short-time perturbation result below.

\begin{proposition} \label{Pro-ShortP-H1} ({\it Energy-Critical Short-time Perturbations when $d>6$}).
Let $I=[t_0,T]$, $w, \wt{w}$ and $g$ be as in Theorem \ref{Thm-Sta-H1}, $d>6$.
Assume that
\begin{align} \label{Short-wtW-H1}
   \|\wt{w}\|_{C(I; H^1)} \leq E
\end{align}
for some positive constant $E$.
Assume also the smallness conditions
\begin{align}
    &\|\wt{w}\|_{\bbw(I)\cap L^2(I; H^{\frac 32}_{-1})}   \leq \delta, \label{Short-H1-ve.0}\\
  \|g\|_{C([t_0,T];\bbr^+)} \leq \ve,\ \
    &  |w(t_0) - \wt{w}(t_0)|_{H^1} \leq \ve, \ \
    \| e \|_{L^2(I; H^{\frac 12 }_1) + N^1(I)} \leq \ve \label{Short-H1-ve}
\end{align}
for some $0<\ve\leq \delta$,
where
$\delta = \delta (C_T, D'_T, E)>0$ is a small constant
and $C_T, D'_T$  are as in Theorem \ref{Thm-Sta-H1}.
Then, we have
\begin{align}
   & \|w-\wt{w}\|_{{\bbx(I)}} \leq  C(C_T, D'_T, E) \ve , \label{Short-H1.1} \\
   & \|w-\wt{w}\|_{{S}^1(I)\cap L^2(I;H^\frac 32_{-1})} \leq  C(C_T, D'_T, E) \ve^{\frac{4}{d-2}},  \label{Short-H1.2} \\
   & \|w\|_{{S}^1(I)\cap L^2(I;H^\frac 32_{-1})} \leq C(C_T, D'_T, E),  \label{Short-H1.3} \\
   & \|F(w) - F(\wt{w})\|_{\bby(I)}  \leq C(C_T, D'_T, E) \ve, \label{Short-H1.4} \\
   & \| F(w) - F(\wt{w}) \|_{N^1(I)} \leq  C(C_T, D'_T, E) \ve^{\frac{4}{d-2}}.   \label{Short-H1.5}
\end{align}
where $\bbx(I)$, $\bby(I)$ are exotic Strichartz spaces as in Section \ref{Sec-Intro}.
\end{proposition}

\begin{remark}
The exotic Strichartz space $\bbx(I)$ and $\bby(I)$
are used to
deal with the non-Lipschitzness of the derivatives of nonlinearity when $d>6$.
Moreover,
the local smoothing spaces are
introduced primarily
to treat the
lower order perturbations of the Laplacian
arising in the operator $e^{-\Phi}\Delta (e^{\Phi} \cdot)$.
\end{remark}

In order to prove Proposition \ref{Pro-ShortP-H1},
we first prove Lemma \ref{Lem-Short-H1-rough} below .

\begin{lemma} \label{Lem-Short-H1-rough}
Consider the situations in Proposition \ref{Pro-ShortP-H1}.
We have that for $\delta= \delta(C_T, D'_T, E)$ small enough,
\begin{align}
    & \|\wt{w}\|_{S^1(I)\cap L^2(I;H^\frac 32_{-1})} \leq C( C_T, D'_T, E), \label{esti-wtw-S1LS-bdd-Short-H1} \\
    & \|w\|_{\bbw(I) \cap L^2(I; H^\frac 32_{-1})} \leq C(C_T, D'_T, E) \delta, \label{esti-w-bbWLS-da-Short-H1} \\
    & \|w\|_{\bbx(I)} \leq  C(C_T, D'_T, E) \delta^{\frac{1}{d+2}}. \label{esti-w-bbX-da-Short-H1}
\end{align}
\end{lemma}

\begin{remark}
Unlike  in \cite{KV13} (and also \cite{TV05}),
it is more delicate  here
to derive the smallness bound \eqref{esti-w-bbX-da-Short-H1} of $w$ in the exotic Strichartz space $\bbx(I)$,
because of the lower order perturbations in Equation \eqref{equa-w-p}.
Below we first prove the smallness bound \eqref{esti-w-bbWLS-da-Short-H1} of $w$
in the local smoothing space $L^2(I; H^\frac 32_{-1})$,
with which we are able to control the lower order perturbations
and then obtain the estimate \eqref{esti-w-bbX-da-Short-H1}.
\end{remark}

{\bf Proof. }
We first prove \eqref{esti-wtw-S1LS-bdd-Short-H1}.
Applying Theorem \ref{Thm-Stri} to \eqref{equa-wtw-p} and
using \eqref{ineq-W.2}, \eqref{Short-wtW-H1} and \eqref{Short-H1-ve.0} we have
\begin{align} \label{esti-wtw-S1.0}
    \|\wt{w}\|_{S^1(I)}
    \leq& C_T (|\wt{w}(t_0)|_{H^1}
              + \|e^{\frac{4}{d-2} \Re \Phi} F(\wt{w})\|_{L^2(I; W^{1,\frac{2d}{d+2}})}
              + \| e \|_{N^1(I)+L^2(I; H^{\frac 12 }_1)} ) \nonumber \\
    \leq& C_T (E+ CD'_T\|\wt{w}\|_{\bbw(I)}^{\frac{d+2}{d-2}} +\ve) \nonumber \\
    \leq& C_T (E+ C D'_T \delta^{\frac{d+2}{d-2}} +\ve),
\end{align}
which yields  \eqref{esti-wtw-S1LS-bdd-Short-H1}
if $\delta = \delta(D'_T,E)$ is small enough such that
$CD'_T \delta^{\frac{d+2}{d-2}} + \delta \leq E$.

In order to prove \eqref{esti-w-bbWLS-da-Short-H1},
again applying Theorem \ref{Thm-Stri} to \eqref{equa-wtw-p}
and using the H\"older inequality \eqref{ineq-W.2} and \eqref{Short-H1-ve.0} we have
\begin{align*}
  \|U(\cdot,t_0)\wt{w}(t_0)\|_{\bbw(I) \cap L^2(I; H^\frac 32_{-1}) }
  \leq& \|\wt{w}\|_{\bbw(I) \cap L^2(I; H^\frac 32_{-1})}
       + C C_T D'_T \|\wt{w}\|_{\bbw(I)}^{\frac{d+2}{d-2}}  \nonumber \\
      &  + C_T \|e\|_{N^1(I) + L^2(I; H^\frac 12_{-1})} \nonumber \\
  \leq&\delta + CC_T D'_T\delta^{\frac{d+2}{d-2}} + C_T \ve.
\end{align*}
Moreover, by the homogeneous Strichartz estimates and \eqref{Short-H1-ve},
\begin{align*}
    \|U(\cdot, t_0)(w(t_0)-\wt{w}(t_0))\|_{\bbw(I) \cap L^2(I; H^\frac 32_{-1}) }
    \leq C|w(t_0)-\wt{w}(t_0)|_{H^1} \leq C \ve.
\end{align*}

Thus, we obtain
\begin{align} \label{esti-w-homo-bbWLS}
     \|U(\cdot, t_0) w(t_0)\|_{\bbw(I) \cap L^2(I; H^\frac 32_{-1}) }
     \leq C_1(C_T, D'_T)\delta.
\end{align}

Arguing as above
and using \eqref{esti-w-homo-bbWLS}
we deduce from Equation \eqref{equa-w-p} that
\begin{align*}
    \|w\|_{\bbw(I) \cap L^2(I; H^\frac 32_{-1}) }
    \leq& \| U(\cdot, t_0)w(t_0)\|_{\bbw(I) \cap L^2(I; H^\frac 32_{-1}) }
           + C C_T D'_T\|w\|_{\bbw(I)}^{\frac{d+2}{d-2}}  \\
    \leq& C_1(C_T, D'_T)\delta
          + C C_T D'_T\|w\|_{\bbw(I)}^{\frac{d+2}{d-2}}.
\end{align*}
Then, in view of \cite[Lemma A.1]{BRZ18},
taking $\delta=\delta(C_T, D'_T)$ smaller
such that
$C_1(C_T, D'_T)\delta < (1-\frac 1\a)(\a C C_T D'_T)^{-\frac{1}{\a-1}}$,
we obtain \eqref{esti-w-bbWLS-da-Short-H1}.

It remains to prove \eqref{esti-w-bbX-da-Short-H1}.
For this purpose,
we see from \eqref{equa-wtw-p} that
\begin{align*}
  i\p_t \wt{w} = \Delta \wt{w} + (b\cdot \na + c)\wt{w} - e^{\frac{4}{d-2} \Re\psi} F(\wt{w}) +e,
\end{align*}
where $b = 2 \na \Phi$ and $c= \Delta \Phi + \sum_{j=1}^d (\p_j \Phi)^2$.
This yields that
\begin{align} \label{esti-wtw-homo-bbX-p.0}
    \|e^{-i(\cdot-t_0)\Delta} \wt{w}(t_0) &\|_{\bbx(I)}
    \leq  \|\wt{w}\|_{\bbx(I)}
           + \bigg\|\int_0^\cdot e^{-i(\cdot-s)\Delta}( b\cdot \na+ c)\wt{w}ds \bigg\|_{\bbx(I)} \nonumber \\
        &   + \bigg\|\int_0^\cdot e^{-i(\cdot-s)\Delta}e^{\frac{4}{d-2}\Re \Phi} F(w)ds \bigg\|_{\bbx(I)}
           +  \bigg\|\int_0^\cdot e^{-i(\cdot-s)\Delta} e(s) ds \bigg\|_{\bbx(I)} \nonumber \\
    & =:K_0+ K_1 + K_2 +K_3.
\end{align}
Note that, by \eqref{bbX-LpS0},
\begin{align} \label{esti-K0-H1-p}
   K_0
   \leq C \|\wt{w}\|_{\bbw(I)}^{\frac{1}{d+2}} \|\wt{w}\|_{S^1(I)}^{\frac{d+1}{d+2}}.
\end{align}
Moreover, \eqref{bbX-bbYLSN1} yields that
\begin{align*}
   K_1 \leq C\|(b\cdot \na+c) \wt{w}\|_{L^2(I;H^\frac 12_{1})}.
\end{align*}
We see that,
$\|(b\cdot \na+c) \wt{w}\|_{H^\frac 12_{1}}= |\Psi_q \<x\>^{-1} \<\na\>^{\frac 32} \wt{w}|_2$,
where $\Psi_q:=\<x\>\<\na\>^{\frac 12} (b\cdot \na+c) \<\na\>^{-\frac 32} \<x\> \in S^0$.
Then, using Lemma \ref{Lem-L2-Bdd} and \eqref{Short-H1-ve}
we have for some $l\geq 1$,
\begin{align*}
   \|(b\cdot \na+c) \wt{w}\|_{H^\frac 12_{1}}
   \leq C \sup\limits_{t\in I} |(ib(t)\cdot \xi + c(t))|_{S^1}^{(l)} |\<x\>^{-1} \<\na\>^{\frac 32} \wt{w}|_2
   \leq C \delta \|\wt{w}\|_{H^{\frac 32}_{-1}}.
\end{align*}
This yields that
\begin{align} \label{esti-K1-H1-p}
  K_1 \leq C \| (b\cdot \na+ c)\wt{w}\|_{L^2(I;H^\frac 12_{1})}
      \leq C \delta \|\wt{w}\|_{L^2(I;H^\frac 32_{-1})}.
\end{align}
We also deduce  from \eqref{bbX-bbYLSN1} that
\begin{align*}
   K_2 \leq C \|e^{\frac{4}{d-2}\Re \Phi} F(\wt{w})\|_{\bby(I)}.
\end{align*}
The product rule for fractional derivatives (see, e.g., \cite[Chapter 2.1]{T00})
implies
\begin{align*}
     \|e^{\frac{4}{d-2}\Re\Phi} F(\wt{w})\|_{\bby(I)}
     \leq& C\|\<\na\>^{\frac{4}{d+2}}e^{\frac{4}{d-2}\Re \Phi}\|_{C(I;L^\9)}
            \|F(\wt{w})\|_{L^{\frac d2}(I; L^{\frac{2d^2(d+2)}{d^3+4d^2+4d-16}})} \\
         & + C \|e^{\frac{4}{d-2}\Re \Phi}\|_{C(I;L^\9)}
            \|\<\na\>^{\frac{4}{d+2}} F(\wt{w})\|_{L^{\frac d2}(I; L^{\frac{2d^2(d+2)}{d^3+4d^2+4d-16}})} \\
     \leq& C \|e^{\frac{4}{d-2}\Re \Phi}\|_{C(I;W^{1,\9})} \|F(\wt{w})\|_{\bby(I)},
\end{align*}
which along with  \eqref{bbX-LpS0}, \eqref{F-bbY-bbX}
implies that
\begin{align} \label{esti-K2-H1-p}
  K_2  \leq CD'_T \|\wt{w}\|_{\bbx(I)}^{\frac{d+2}{d-2}}
       \leq CD'_T \|\wt{w}\|_{\bbw(I)}^{\frac{1}{d-2}} \|\wt{w}\|_{S^1(I)}^{\frac{d+1}{d-2}}.
\end{align}
Regarding $K_3$, by  \eqref{bbX-bbYLSN1},
\begin{align} \label{esti-K3-H1-p}
   K_3  \leq C \|e\|_{N^1(I) + L^2(I;H^\frac 12 _{1})} .
\end{align}
Thus, plugging \eqref{esti-K0-H1-p}-\eqref{esti-K3-H1-p} into \eqref{esti-wtw-homo-bbX-p.0}
and using \eqref{Short-H1-ve.0}, \eqref{Short-H1-ve} and \eqref{esti-wtw-S1LS-bdd-Short-H1}
yield
\begin{align} \label{esti-wtw-homo-bbX-p}
    \|e^{-i(\cdot-t_0)\Delta} \wt{w}(t_0) \|_{\bbx(I)}
    \leq& C\|\wt{w}\|^{\frac{1}{d+2}}_{\bbw(I)} \|\wt{w}\|^{\frac{d+1}{d+2}}_{S^1(I)}
         + C \delta \|\wt{w}\|_{L^2(I;H^\frac 32_{-1})} \nonumber \\
        & + CD'_T \|\wt{w}\|_{\bbw(I)}^{\frac{1}{d-2}} \|\wt{w}\|_{S^1(I)}^{\frac{d+1}{d-2}}
         + C \|e\|_{N^1(I)+L^2(I;H^\frac 12 _{1})} \nonumber  \\
    \leq& C_2(C_T,D'_T, E)\delta^{\frac{1}{d+2}}.
\end{align}

Moreover,
by \eqref{bbX-H1} and \eqref{Short-H1-ve},
\begin{align*}
   \|e^{-i(\cdot-t_0)\Delta} (w(t_0)- \wt{w}(t_0)) \|_{\bbx(I)}
   \leq C |w(t_0)-\wt{w}(t_0)|_{H^1}
   \leq C \ve.
\end{align*}

Thus, we obtain
\begin{align*}
   \|e^{-i(\cdot-t_0)\Delta} w(t_0) \|_{\bbx(I)}
     \leq C_3(C_T,D'_T, E)\delta^{\frac{1}{d+2}}.
\end{align*}

Now,
similarly as above,
we deduce from Equation \eqref{equa-w-p} that
\begin{align*}
    \|w \|_{\bbx(I)}
    \leq&  \|e^{-i(\cdot-t_0)\Delta} w(t_0)\|_{\bbx(I)}
           + \bigg\|\int_0^\cdot e^{-i(\cdot-s)\Delta} (b\cdot \na+ c)w(s) ds \bigg\|_{\bbx(I)} \\
        &   + \bigg\|\int_0^\cdot e^{-i(\cdot-s)\Delta}e^{\frac{4}{d-2}\Re\Phi} F(w(s))ds \bigg\|_{\bbx(I)} \\
    \leq& C_3(C_T,D'_T,E)\delta^{\frac{1}{d+2}} + C \delta \|w\|_{L^2(I;H^\frac 32_{-1})}
           + C D'_T \|w\|_{\bbx(I)}^{\frac{d+2}{d-2}}
           + C \ve \\
    \leq& C_4(C_T,D'_T,E)\delta^{\frac{1}{d+2}} + CD'_T \|w\|_{\bbx(I)}^{\frac{d+2}{d-2}},
\end{align*}
where the last step is due to \eqref{esti-w-bbWLS-da-Short-H1}.

Therefore,
taking $\delta= \delta(C_T, D'_T, E)$ even smaller such that
$ C_4(C_T,D'_T,E) \delta^{\frac{1}{d+2}} < (1-\frac 1\a)(\a C D'_T)^{-\frac{1}{\a-1}}$
and using \cite[Lemma A.1]{BRZ18} we obtain \eqref{esti-w-bbX-da-Short-H1}.
\hfill $\square$ \\

{\bf Proof of Proposition \ref{Pro-ShortP-H1}.}
We first estimate $\|w-\wt{w}\|_{\bbx(I)}$.

For this purpose, we note that $z:= w-\wt{w}$ satisfies the equation
\begin{align} \label{equa-z-p-H1}
   i\p_t z =& e^{-\Phi}\Delta(e^{\Phi} z) - e^{\frac{4}{d-2}\Re \Phi} (F(z+\wt{w}) - F(\wt{w})) -e \nonumber \\
           =& \Delta z  + (b\cdot \na + c)z   - e^{\frac{4}{d-2}\Re \Phi} (F(z+\wt{w}) - F(\wt{w})) -e.
\end{align}
This yields that
\begin{align*}
    \|z\|_{\bbx(I)}
    \leq& \|e^{-i(\cdot-t_0) \Delta}z(t_0)\|_{\bbx(I)}
         + \bigg \|\int_{t_0}^\cdot e^{-i(\cdot-s)\Delta} (b\cdot \na +c)zds \bigg \|_{\bbx(I)}  \nonumber \\
        & + \bigg \|\int_{t_0}^\cdot e^{-i(\cdot-s)\Delta}e^{\frac{4}{d-2} \Re \Phi} (F(z+\wt{w}) - F(\wt{w})) ds \bigg \|_{\bbx(I)}
         + \bigg \|\int_{t_0}^\cdot e^{-i(\cdot-s)\Delta} e(s) ds \bigg \|_{\bbx(I)}  \\
    =:& J_0 + J_1 + J_2 + J_3.
\end{align*}
First, Theorem \ref{Thm-Stri*}, \eqref{X0-bbX-S1} and \eqref{Short-H1-ve}  yield that
$$  J_0 \leq C\|e^{-i(\cdot-t_0) \Delta}z(t_0)\|_{S^1(I)} \leq C|z(t_0)|_{H^1} \leq C \ve.$$
Moreover, similarly to \eqref{esti-K1-H1-p},
using \eqref{Short-H1-ve.0}, \eqref{Short-H1-ve} and \eqref{esti-w-bbWLS-da-Short-H1} we have
\begin{align} \label{esti-J1-H1-p}
     J_1 \leq C \|(b\cdot \na + c)z \|_{L^2(I;H^\frac 12_{1})}
         \leq C \ve \|z\|_{L^2(I;H^\frac 32_{-1})}
         \leq C_1(C_T, D'_T, E) \ve.
\end{align}
We also use \eqref{bbX-bbYLSN1} and \eqref{Fz-bbY-S1bbX} to get that
\begin{align} \label{esti-J2-H1-p}
   J_2 \leq& C \| e^{\frac{4}{d-2}\Re \Phi} (F(z+\wt{w}) - F(\wt{w}))\|_{\bby(I)} \nonumber \\
       \leq& CD'_T  \|(F(z+\wt{w}) - F(\wt{w}))\|_{\bby(I)} \nonumber \\
       \leq& CD'_T (\|\wt{w}\|_{\bbx(I)}^\frac{8}{d^2-4} \|\wt{w}\|_{S^1(I)}^\frac{4d}{d^2-4}
               + \|z\|_{\bbx(I)}^\frac{8}{d^2-4} \|z\|_{S^1(I)}^\frac{4d}{d^2-4} )
             \|z\|_{\bbx(I)}.
\end{align}
Note that, by \eqref{bbX-LpS0} and \eqref{Short-H1-ve.0},
\begin{align} \label{esti-wtw-bbx-S1}
   \|\wt{w}\|_{\bbx(I)}
   \leq C\|\wt{w}\|_{\bbw(I)}^{\frac {1}{d+2}}
          \|\wt{w}\|_{S^1(I)}^{\frac{d+1}{d-2}}
   \leq C \delta^{\frac{1}{d+2}} \|\wt{w}\|_{S^1(I)}^{\frac{d+1}{d-2}} .
\end{align}
Plugging \eqref{esti-wtw-bbx-S1} into \eqref{esti-J2-H1-p} and
using \eqref{esti-wtw-S1LS-bdd-Short-H1} we obtain
\begin{align*}
     J_2 \leq& C_2(C_T,D_T',E) \delta^{\frac{8}{(d-2)(d+2)^2}} \|z\|_{\bbx(I)}
               +  C D'_T \|z\|_{S^1(I)}^\frac{4d}{d^2-4}  \|z\|^{1+\frac{8}{d^2-4}}_{\bbx(I)}.
\end{align*}
Regarding $J_3$,
similarly to \eqref{esti-K3-H1-p},
by \eqref{bbX-bbYLSN1},
\begin{align*}
     J_3   \leq C\|e\|_{ N^1(I)+ L^2(I;H^\frac 12_1)} \leq  C\ve.
\end{align*}

Thus, combining the estimates of $J_i$ above, $i=0,1,2,3$,  we obtain
that for $\delta = \delta(C_T,D'_T, E)$ small enough
such that
$ C_2(C_T,D_T',E) \delta^{\frac{8}{(d-2)(d+2)^2}}  \leq \frac 12$,
\begin{align} \label{esti-z-bbX-StaH1}
     \|z\|_{\bbx(I)}
    \leq 2 (C_1(C_T,D'_T,E)+2C)\ve
           + 2  CD'_T \|z\|_{S^1(I)}^\frac{4d}{d^2-4}  \|z\|^{1+\frac{8}{d^2-4}}_{\bbx(I)} .
\end{align}

Below we estimate $\|z\|_{S^1(I)}$ and $\|z\|_{L^2(I;H^\frac 32_{-1})}$.
Arguing as in the proof of \eqref{esti-wtw-S1.0},
applying Theorem \ref{Thm-Stri} to \eqref{equa-z-p-H1}
and using   \eqref{ineq-W.2} and \eqref{Short-H1-ve}  we have
\begin{align} \label{esti-z-S1-Short-H1.0}
    \|z\|_{S^1(I) \cap L^2(I; H^{\frac 32 }_{-1})}
   \leq& C_T \|z(t_0)\|_{H^1}
        + C_T  \|e^{\frac{4}{d-2} \Re\Phi} (F(z+\wt{w}) - F(\wt{w}))\|_{N^1(I)} \nonumber \\
       &  + C_T  \| e \|_{N^1(I)+ L^2(I; H^{\frac 12 }_1)}  \nonumber  \\
   \leq& 2C_T \ve + C_T D'_T  \|(F(z+\wt{w}) - F(\wt{w}))\|_{N^1(I)}.
\end{align}

Note that, by H\"older's inequality \eqref{ineq-W.2}, \eqref{Short-H1-ve.0} and \eqref{esti-w-bbWLS-da-Short-H1},
\begin{align} \label{esti-F-Short-H1}
    \|F(z+\wt{w}) - F(\wt{w})\|_{N^0(I)}
   \leq& C (\|\wt{w}\|^{\frac{4}{d-2}}_{\bbw(I)}
            +  \|z\|^{\frac{4}{d-2}}_{\bbw(I)}  )
          \|z\|_{S^1(I)} \nonumber \\
   \leq& C_3(C_T,D'_T,E)\delta^{\frac{4}{d-2}}\|z\|_{S^1(I)}.
\end{align}
Moreover,
arguing as in the proof of  \cite[(3.67)]{KV13}
and using \eqref{Short-H1-ve.0}, \eqref{esti-w-bbX-da-Short-H1} we have
\begin{align} \label{esti-naF-Short-H1}
   \|\na (F(z+\wt{w}) - F(\wt{w}))\|_{N^0(I)}
   \leq& C( \|\na \wt{w}\|_{S^0(I)} \|z\|^{\frac{4}{d-2}}_{X^0(I)}
         + \|z+\wt{w}\|^{\frac{4}{d-2}}_{X^0(I)} \|\na z\|_{S^0(I)} ) \nonumber \\
   \leq& C_4(C_T, D'_T,E) ( \|z\|^{\frac{4}{d-2}}_{\bbx(I)}
         +  \delta^{\frac{4}{d^2-4}}  \|z\|_{S^1(I)}) .
\end{align}

Thus, plugging \eqref{esti-F-Short-H1}, \eqref{esti-naF-Short-H1} into \eqref{esti-z-S1-Short-H1.0}
we get
\begin{align*}
    \|z\|_{S^1(I) \cap L^2(I;H^{\frac 32}_{-1})}
    \leq& 2C_T \ve
          +C_5(C_T,D'_T,E)  (  \|z\|_{\bbx(I)}^{\frac{4}{d-2}}
                         + (\delta^{\frac{4}{d-2}}+\delta^{\frac{4}{d^2-4}}) \|z\|_{S^1(I)} ).
\end{align*}
Taking $\delta = \delta(C_T, D'_T,E)$ small such that
$C_5(C_T,D'_T,E)  (\delta^{\frac{4}{d-2}} + \delta^{\frac{4}{d^2-4}}) \leq \frac 12$
yields
\begin{align} \label{esti-z-S1-Short-H1}
     \|z\|_{S^1(I)\cap L^2(I;H^{\frac 32}_{-1})}
    \leq 4C_T \ve
          + 2C_5(C_T,D'_T,E)  \|z\|^{\frac{4}{d-2}}_{\bbx(I)} .
\end{align}

Now,
plugging \eqref{esti-z-S1-Short-H1} into \eqref{esti-z-bbX-StaH1} we get that
if $c_1:= \frac{8}{d^2-4}$, $c_2:= \frac{24d-16}{(d-2)^2(d+2)}>0$,
\begin{align} \label{esti-z-bbX-StaH1.1}
    \|z\|_{\bbx(I)}
\leq C_6(C_T, D'_T, E) \ve
     +  C_6(C_T, D'_T,E)
         \ve^{\frac{4d}{d^2-4}} (\|z\|_{\bbx(I)}^{1+ c_1}
          + \|z\|_{\bbx(I)}^{1+c_2}).
\end{align}
Since $0<c_1<c_2$,
$\|z\|_{\bbx(I)}^{1+c_1} \leq \|z\|_{\bbx(I)} + \|z\|_{\bbx(I)}^{1+c_2}$,
we have
\begin{align*}
   \|z\|_{\bbx(I)}
   \leq C_6(C_T,D'_T,E)(\ve +  \ve^{\frac{4d}{d^2-4}} \|z\|_{\bbx(I)}
        + 2  \|z\|_{\bbx(I)}^{1+c_2}).
\end{align*}
Then, taking $\delta$ very small such that
$C_6(C_T,D'_T,E) \delta^{\frac{4d}{d^2-4}} \leq \frac 12$,
we come to
\begin{align*}
   \|z\|_{\bbx(I)}
   \leq 2 C_6(C_T,D'_T,E) \ve + 4C_6(C_T,D'_T,E)\|z\|_{\bbx(I)}^{1+c_2}.
\end{align*}
Thus, taking $\delta$ even smaller such that
$2C_6\delta <(1-\frac 1\a)(4\a C_6)^{-\frac{1}{\a-1}}$,
we apply  \cite[Lemma A.1]{BRZ18}
to obtain \eqref{Short-H1.1},
which along with \eqref{esti-z-S1-Short-H1} implies \eqref{Short-H1.2}.

Finally,
\eqref{Short-H1.3} follows from \eqref{Short-H1.2} and \eqref{esti-wtw-S1LS-bdd-Short-H1},
\eqref{Short-H1.4} can be proved by \eqref{Short-H1.1} and similar estimates as in \eqref{esti-J2-H1-p},
and \eqref{Short-H1.5} follow from
\eqref{Short-H1.1}, \eqref{Short-H1.2} and \eqref{esti-naF-Short-H1}.
Therefore, the proof is complete.
\hfill $\square$ \\

{\bf Proof of Theorem \ref{Thm-Sta-H1}.}
Let $\delta=\delta(C_T,D'_T,E)$ be as in Proposition \ref{Pro-ShortP-H1}.
As in the proof of Theorem \ref{Thm-Sta-L2},
since $\|\wt{w}\|_{\bbw(I)} \leq L<\9$,
we can divide $I$ into subintervals $I'_j=[t'_j,t'_{j+1}]$,
such that
$0\leq j\leq l'\leq (\frac{2L}{\delta})^{\frac{2(d+2)}{d-2}}<\9$,
and $\|\wt{w}\|_{\bbw(I'_j)} \leq \frac \delta 2$ for each $0\leq j\leq l'$.
Similarly, since $\|\wt{w}\|_{L^2(I;H^\frac 32_{-1})} \leq L<\9$,
we have another finite partition $I''_{j}=[t''_j, t''_{j+1}]$,
so that
$0\leq j\leq l''\leq (\frac{2L}{\delta})^{2}$
and on each $I''_j$, $\|\wt{w}\|_{L^2(I''_j;H^\frac 32_{-1})} \leq \frac \delta 2$.
Thus, let $\{t_j\} = \{t'_j\} \cup \{t''_j\}$,
we obtain a partition $\{I_j=[t_j,t_{j+1}]\}_{j=0}^l$,
satisfying that
$l\leq (\frac{2L}{\delta})^{\frac{2(d+2)}{d-2}} + (\frac{2L}{\delta})^{2}$
and $\|\wt{w}\|_{\bbw(I_j) \cap L^2(I_j;H^\frac 32_{-1})} \leq \delta$.

Let $C(0): =C(C_T,D'_T,E)$,
$C(j+1) = C(0)(2C_T+C_T D'_T \sum_{k=0}^jC(k))$,
$0\leq j\leq l$,
where $C(C_T,D_T',E)$ is as in Proposition \ref{Pro-ShortP-H1}.
Choose $\ve_*(C_T,D'_T,E,L)$
sufficiently small such that
\begin{align} \label{ve*-H1}
    (2C_T + C_T D'_T \sum\limits_{k=0}^l C(k)) \ve_*^{(\frac{4}{d-2})^{l+1}} \leq \delta.
\end{align}

We claim that on each $I_j$, $0\leq j\leq l$,
estimates \eqref{Short-H1.1}-\eqref{Short-H1.5} hold  with
$I$, $C(C_T,D'_T,E)$, $\ve$
replaced by $I_j$, $C(j)$ and $\ve^{(\frac{4}{d-2})^{j+1}}$, respectively.

Actually,
Proposition \ref{Pro-ShortP-H1} implies that
the claim is true for $j=0$.
Suppose that it is valid for each $0\leq k\leq j<l$.
Then, on the next interval $I_{j+1}$,
applying Theorem \ref{Thm-Stri} to \eqref{equa-z-p-H1}
and using the inductive assumptions and \eqref{ve*-H1} we have
\begin{align*}
    |w(t_{j+1})- \wt{w}(t_{j+1})|_{H^1}
   \leq& C_T |w(t_{0})- \wt{w}(t_{0})|_{H^1}
         + C_T \|e^{\frac{4}{d-2}\Re\Phi} (F(w)- F(\wt{w}))\|_{N^1(t_0,t_{j+1})}  \\
       &  + C_T \|e\|_{L^2(t_0, T;H^{\frac 12}_{1}) + N^1(t_0, T)} \\
   \leq& 2 C_T \ve + C_T D'_T \sum\limits_{k=0}^j  C(k) \ve^{(\frac{4}{d-2})^{k+1}}
   \leq \delta.
\end{align*}
Thus, Proposition \ref{Pro-ShortP-H1} yields that the claim holds on $I_{j+1}$.

Therefore,
using the inductive arguments we prove the claim on any $I_j$, $0\leq j\leq l$.
This yields \eqref{Sta-H1.2}, \eqref{Sta-H1.3} and that for some
$c=c(C_T,D'_T,E,L)>0$,
\begin{align*}
   \|w-\wt{w}\|_{\bbx(I)}
   \leq C'(C_T, D'_T,E) \ve^c.
\end{align*}
Finally, taking into account \eqref{Lp-bbXS1} and \eqref{Sta-H1.2},
we obtain for some $0<c'\leq 1$,
\begin{align*}
    \|w-\wt{w}\|_{L^{\frac{2(d+2)}{d-2}}(I)}
    \leq \|w-\wt{w}\|_{\bbx(I)}^{c'} \|w-\wt{w}\|_{S^1(I)}^{1-c'}
    \leq C''(C_T, D'_T,E) \ve^{c},
\end{align*}
thereby yielding \eqref{Sta-H1.1}.
The proof is complete.
\hfill $\square$

\section{Global well-posedness} \label{Sec-GWP}

This section is mainly devoted to the global well-posedness of \eqref{equa-x} in the
mass and energy critical cases.

\subsection{Mass-critical case} \label{Subsec-GWP-L2}

We first recall the global well-posedness and scattering results
in the deterministic defocusing mass-critical case,
based on the work of Dodson \cite{D12, D16.1,D16.2}
\footnote{\cite{D12,D16.1, D16.2} study the equation $i\p_t u= -\Delta u + |u|^{\frac 4 d} u$,
which can be easily transformed into \eqref{equa-u-L2}
by reversing the time.
Hence, the results in \cite{D12,D16.1, D16.2} also hold for \eqref{equa-u-L2}.}.

\begin{theorem} (\cite{D12,D16.1, D16.2}) \label{Thm-L2GWP-Det}
For any $u_0\in L^2$,
there exists a unique global $L^2$-solution $u$ to the equation
\begin{align} \label{equa-u-L2}
   i\p_t u =&   \Delta u -  |u|^{\frac{4}{d}} u,  \\
   u(0)=& u_0  \nonumber
\end{align}
with $d\geq 1$. Moreover,
\begin{align} \label{globdd-u-L2-Det}
   \|u\|_{V(\bbr) \cap L^2(\bbr; H^\frac 12_{-1})} \leq B_0(|u_0|_2) <\9,
\end{align}
where $B_0(|u_0|_2)$ depends continuously on $u_0$ in $L^2$,
and $u$ scatters at infinity,
i.e., there exist  $u_{\pm} \in L^2$ such that
\begin{align} \label{sca-det-L2}
  | e^{it\Delta} u(t) - u_{\pm}|_2 \to 0,\ as\ t\to \pm\9.
\end{align}
\end{theorem}

We remark that the bound of $\|u\|_{L^2(\bbr; H^\frac 12_{-1})}$
in \eqref{globdd-u-L2-Det}   follows standardly from
Strichartz estimates and the bound of $\|u\|_{V(\bbr)}$,
and the continuity of $B_0(|u_0|_2)$ on $u_0$
is a consequence of the mass-critical stability result Lemma $3.6$ of \cite{TVZ07}.

We also need the following boundedness of $X$ in the space $L^2$.
\begin{lemma} \label{Lem-bdd-L2}
Assume the conditions of Theorem \ref{Thm-GWP} $(i)$ to hold.
Then, for each $X_0\in L^2$,
we have $\bbp$-a.s.
\begin{align} \label{Ito-L2}
  |X(t)|_2^2 = |X_0|_2^2 + 2 \sum\limits_{k=1}^N\int_0^t \int \Re G_k(s) |X(s)|^2 dx d\beta_k(s),
  \ \ 0\leq t<\tau^*,
\end{align}
where $\tau^*$ is the maximal existing time as in Theorem \ref{Thm-LWP}.
Moreover, for any $0<T<\9$, $p\geq 1$,
\begin{align}    \label{bdd-X-L2}
   & \bbe \|X\|^p_{C([0,\tau^*\wedge T);L^2)} \leq C(p,T) <\9.
\end{align}
In particular,
\begin{align} \label{globdd-M}
   M_T:= \sup\limits_{0\leq t<\tau^*\wedge T} |X(t)|_2 \leq C(T) <\9,\ \ a.s..
\end{align}
\end{lemma}

{\bf Proof.} The It\^o formula \eqref{Ito-L2} was obtained in \cite[(6.1)]{HRZ18}.
The proof of \eqref{bdd-X-L2}
is similar to that of   \cite[Lemma $3.6$]{BRZ16},
based on the
Burkholder-Davis-Gundy inequality and the Gronwall inequality.
We omit the details here for simplicity. \hfill $\square$  \\

{\bf Proof of Theorem \ref{Thm-GWP} $(i)$.}
(Mass-Critical Case).
Let $X$ be the unique $L^2$-solution to \eqref{equa-x}
on the maximal time interval $[0,\tau^*)$.
In order to prove the global existence of $X$,
we only need to prove the global bound \eqref{gloexist-L2} for any $0<T<\9$.

For this purpose, we proceed as follows:
we first consider a small (random) time interval $I_1$,
determined by the smallness condition \eqref{Sta-L2-ve} of Theorem \ref{Thm-Sta-L2},
and we apply the rescaling transformation
and   the stability result Theorem \ref{Thm-Sta-L2}
to obtain the bound of $V(I_1)$-norm of the resulting random solution
and so of the stochastic solution $X$.
Then, we construct consecutively small (random) subintervals $I_j$, $2\leq j\leq l$,
on which we obtain the bound of $\|X\|_{V(I_j)}$ by using
Theorem \ref{Thm-Rescale-sigma}
as well as  Theorem \ref{Thm-Sta-L2}.
At last, by virtue of the global $L^2$ bound in Lemma \ref{Lem-bdd-L2},
we are able to show that the total number $l$  is finite almost surely,
thus we obtain the desirable global bound \eqref{gloexist-L2}.

Let us start with the first step.

{\bf Step $1$.}
Set
$g(t):= \sum_{k=1}^N|\int_0^t g_k(s) d\beta_k(s)| + \int_0^t g_k^2(s) ds$,
$0\leq t\leq T$.
Let $\tau_1:= \inf\{ 0<t<T\wedge \tau^*: g(t) \geq \ve_1(t)\} \wedge (T\wedge \tau^*)$
with
\begin{align} \label{ve1-L2-GWP}
   \ve_1(t) = \frac{\ve_*(C_t,D_t, |X_0|_2)}{D_0(|X_0|_2)},
\end{align}
where
$\ve_*(C_t,D_t, |X_0|_2) := \ve_*(C_t,D_t,|X_0|_2, 0,  B_0(|X_0|_2)$
is as in Theorem \ref{Thm-Sta-L2},
and $D_0(|X_0|_2) = C_0(B_0(|X_0|_2) + (B_0(|X_0|_2))^{1+\frac 4d})$
with $C_0$ specified in  \eqref{esti-e1-L2} below.

Let $\vf$ be as in \eqref{vf} with $\sigma\equiv0$.
By Theorem \ref{Thm-Rescale-sigma},
$v_1:= e^{-\vf} X$
satisfies the random equation \eqref{equa-RNLS}
with $\a= 1+\frac 4d$.
In order to obtain the bound of $\|v_1\|_{V(0,\tau_1)}$,
we compare $v_1$ with the solution $\wt{v}_1$ to \eqref{equa-NLS}
(or, equivalently, \eqref{equa-NLS*}) with
$\a=1+\frac 4d$
and with the same initial datum, i.e.,
$\wt{v}_1(0) = v_1(0) = X_0$.

Then, Theorem \ref{Thm-L2GWP-Det} implies that
$\wt{v}_1$ exists globally and satisfies that
\begin{align} \label{esti-u1-L2}
   \|\wt{v}_1\|_{V(0,\tau_1 ) \cap L^2(0,\tau_1; H^\frac 12_{-1})}
   \leq B_0(|\wt{v}_1(0)|_2) = B_0(|X_0|_2) <\9.
\end{align}

Moreover,
in order to estimate the error term \eqref{Error-NLS}, i.e.,
\begin{align*}
   e_1 :=  -(b \cdot \na + c)\wt{v}_1
       - (1- e^{\frac 4d \Re \vf} ) F(\wt{v}),
\end{align*}
where  $b$, $c$ are as in \eqref{b} and \eqref{c}
with $\sigma\equiv 0$, respectively,
we note that
\begin{align} \label{esti-e1-L2.0}
   &\|e_1\|_{ N^0(0,\tau_1) + L^2(0,\tau_1; H^{-\frac 12}_1)}  \nonumber \\
   \leq& \|(b\cdot \na+ c) \wt{v}_1\|_{L^2(0,\tau_1; H^{-\frac 12}_1)}
         + \| (1- e^{\frac 4d \Re \vf} ) F(\wt{v}_1) \|_{L^{\frac{2(d+2)}{d+4}}((0,\tau_1)\times \bbr^d)}.
\end{align}
To estimate the first term on the right-hand side above,
we see that, by \eqref{asymflat},
\begin{align*}
  \<x\> \<\na\>^{-\frac 12} (b\cdot \na + c) \wt{v}_1
  = \Psi_p \<x\>^{-1} \<\na\>^\frac 12 \wt{v}_1,
\end{align*}
where
$\Psi_p := \<x\> \<\na\>^{-\frac 12} (b\cdot \na + c) \<\na\>^{-\frac 12} \<x\>$
is a pseudo-differential operator of zero order.
Then, using Lemma \ref{Lem-L2-Bdd} we get for some $m\geq 1$,
\begin{align} \label{esti-e1-L2.1}
   \|(b\cdot \na+ c) \wt{v}_1\|_{L^2(0,\tau_1; H^{-\frac 12}_1)}
    \leq& C \sup\limits_{0\leq t\leq \tau_1} |p(t)|_{S^0}^{(m)}
         \|\wt{v}_1 \|_{L^2(0,\tau_1; H^{\frac 12}_{-1})} \nonumber \\
   \leq& C'  \sup\limits_{0\leq t\leq\tau_1 }  g(t)
             \|\wt{v}_1 \|_{L^2(0,\tau_1; H^{\frac 12}_{-1})}.
\end{align}
Moreover, using  \eqref{ineq-V}
and the inequality  $|1-e^x|\leq e|x|$ for $|x|\leq 1$,
we have
\begin{align} \label{esti-e1-L2.2}
    \|(1-e^{\frac 4d \Re \vf}) F(\wt{v}_1)\|_{L^{\frac{2(d+2)}{d-2}}((0,\tau_1) \times \bbr^d)}
    \leq C''  \sup\limits_{0\leq t\leq \tau_1} g(t)
            \|\wt{v}_1\|^{1+ \frac 4d}_{V(0,\tau_1)}.
\end{align}
Hence, plugging \eqref{esti-e1-L2.1} and \eqref{esti-e1-L2.2} into \eqref{esti-e1-L2.0}
and using \eqref{globdd-u-L2-Det} we  arrive at
\begin{align} \label{esti-e1-L2}
        \|e_1\|_{L^2(0,\tau_1; H^{-\frac 12}_1) + N^0(0,\tau_1)}
   \leq&  C  \sup\limits_{0\leq t\leq \tau_1} g(t)
           (\|\wt{v}_1\|_{L^2(0,\tau_1; H^{\frac 12}_{-1})}
        +   \|\wt{v}_1\|^{1+\frac 4 d}_{V(0,\tau_1)}) \nonumber \\
   \leq& C( B_0(|X_0|_2) + (B_0(|X_0|_2))^{1+\frac 4d}) \ve_1 (\tau_1)  \nonumber \\
   \leq& \ve_*(C_{\tau_1}, D_{\tau_1}, |X_0|_2),
\end{align}
where $C_0 := \max\{C', C''\}$.

Thus, in view of Theorem \ref{Thm-Sta-L2},
we obtain
\begin{align} \label{esti-v1-L2*}
    \|v_1\|_{V(0,\tau_1)} \leq  C(C_{\tau_1}, D_{\tau_1}, |X_0|_2, 0, B_0(|X_0|_2))
    =: C(C_{\tau_1}, D_{\tau_1}, |X_0|_2) ,
\end{align}
which implies that
\begin{align}
    \|X\|_{V(0,\tau_1)} \leq \|e^{\vf}\|_{C([0,\tau_1];L^\9)}  C(C_{\tau_1}, D_{\tau_1}, |X_0|_2).
\end{align}

Thus, \eqref{esti-v1-L2*} implies \eqref{gloexist-L2} if $\tau_1 = T\wedge \tau^*$.
Otherwise, we come to the next step.

{\bf Step $2$.}
Set $\sigma_0:=0$,
$\sigma_1 := \tau_1$.
For $j\geq 1$,
we define random times inductively:
\begin{align*}
   & \tau_{j+1} := \inf\{t\in (0,(T\wedge \tau^*)-\sigma_j): g_{\sigma_{j}}(t) \geq \ve_{j+1}(t) \} \wedge (T\wedge \tau^* - \sigma_j), \\
   & \sigma_{j}:= \sum_{k=1}^{j} \tau_k (\leq T\wedge \tau^*), \ \ l:=\inf\{j\geq 1:\sigma_j = T\wedge \tau^*\}.
\end{align*}
Here,
Let $g_{\sigma_j}(t) := \sum_{k=1}^N |\int_{\sigma_j}^{\sigma_{j}+t} g_k(s) d\beta_k(s)| + \int_{\sigma_j}^{\sigma_{j}+t} g^2_k(s) ds$
and
\begin{align} \label{vej-L2-GWP}
   \ve_{j+1}(t):= \frac{\ve_*(C_{\sigma_j+t}, D_{\sigma_j+t}, |X(\sigma_j)|_2)}{D_0(|X(\sigma_j)|_2)}
\end{align}
with
$\ve_*(C_{\sigma_j+t}, D_{\sigma_j+t}, |X(\sigma_j)|_2)
:= \ve_*(C_{\sigma_j+t}, D_{\sigma_j+t}, |X(\sigma_j)|_2, 0, B_0(|X(\sigma_j)|_2))$
as in Theorem \ref{Thm-Sta-L2},
and $D_0(|X(\sigma_j)|_2)$ is defined similarly to $D_0(|X_0|_2)$.
We see that
$\tau_{j+1}$ (resp. $\sigma_{j}$) are $\mathscr{G}_j(t):=\mathscr{F}(\sigma_j+t)$ (resp. $\mathscr{F}(t))$-stopping times, $0\leq t\leq T$.

We use the inductive arguments to obtain the bound of $\|X\|_{V(0,\sigma_j)}$ for any $1\leq j\leq l$.
Suppose that for each $1\leq k\leq j<l$,
\begin{align} \label{esti-v-Xj}
   \|X\|_{V(0,\sigma_k)}
   \leq \sum\limits_{i=0}^{k-1} \|e^{\vf_{\sigma_i}}\|_{C([0,\tau_{i+1}];L^\9)}
         C(C_{\sigma_{i+1}}, D_{\sigma_{i+1}}, |X(\sigma_j)|_2),
\end{align}
where $\vf_{\sigma_i}$ is as in  \eqref{vf} with $\sigma_i$ replacing $\sigma$,
and
$C(C_{\sigma_{i+1}}, D_{\sigma_{i+1}}, |X(\sigma_j)|_2) :=
C(C_{\sigma_{i+1}}, D_{\sigma_{i+1}}, |X(\sigma_j)|_2, 0, B_0(|X(\sigma_j)|_2))$
is as in Theorem \ref{Thm-Sta-L2}.
Below we  show that \eqref{esti-v-Xj} also holds when $k$ is replaced by $j+1$.

For this purpose,
we apply Theorem \ref{Thm-Rescale-sigma} to obtain that
\begin{align} \label{res-zj1}
     v_{j+1}(t) := e^{-\vf_{\sigma_j}(t)} X(\sigma_j +t), \ \ 0\leq t<(T\wedge \tau^*)-\sigma_j.
\end{align}
satisfies the equation
\begin{align} \label{equa-vj1-L2}
   i\p_t v_{j+1}  &=   e^{-\vf_{\sigma_j}}\Delta (e^{\vf_{\sigma_j}}v_{j+1})
                 - e^{\frac 4d \Re \vf_{\sigma_j}} F(v_{j+1}),  \\
    v_{j+1}(0)&= X(\sigma_j).               \nonumber
\end{align}

Similarly to Step $1$,
we  compare \eqref{equa-vj1-L2} with the equation
\begin{align} \label{equa-wtu-j1.0}
    i \p_t \wt{v}_{j+1}
    =   \Delta \wt{v}_{j+1} -  F(\wt{v}_{j+1}),
\end{align}
or equivalently,
\begin{align} \label{equa-wtu-j1}
    i \p_t \wt{v}_{j+1}
    =  e^{-\vf_{\sigma_j}}\Delta (e^{\vf_{\sigma_j}}\wt{v}_{j+1})
                - e^{\frac 4d \Re \vf_{\sigma_j}}  F(\wt{v}_{j+1}) + e_{j+1},
\end{align}
with $\wt{v}_{j+1}(0) = X(\sigma_j)$
and
\begin{align} \label{err-ej1-L2}
   e_{j+1}:= -(b_{\sigma_j}(t) \cdot \na  + c_{\sigma_j}(t))\wt{v}_{j+1}
              -(1- e^{\frac 4d \Re \vf_{\sigma_j}(t)}) F(\wt{v}_{j+1}).
\end{align}
where  $b_{\sigma_j}$ and $c_{\sigma_j}$ are as in \eqref{b} and \eqref{c}
with $\sigma_j$ replacing $\sigma$, respectively.

Again, Theorem \ref{Thm-L2GWP-Det} yields that
$\wt{v}_{j+1}$ exists globally and
\begin{align} \label{esti-wtu-Vj1}
    \|\wt{v}_{j+1}\|_{V(0,\tau_{j+1})} \leq B_0 (|\wt{v}_{j+1}(0)|_2) = B_0 (|X(\sigma_j)|_2).
\end{align}
This implies that,
similarly to \eqref{esti-e1-L2},
\begin{align} \label{esti-ej1}
        \|e_{j+1}\|_{N^0(0,\tau_{j+1} ) + L^2(0,\tau_{j+1}; H^{-\frac 12}_{1})}
   \leq& C \sup\limits_{0\leq t\leq \tau_{j+1}} g_{\sigma_j}(t)
          (\|\wt{v}_{j+1}\|_{L^2(0,\tau_{j+1} ; H^\frac 12_{-1})}
       + \|\wt{v}_{j+1}\|^{1+\frac 4 d}_{V(0,\tau_{j+1} )}) \nonumber \\
   \leq& C(B_0(|X_{\sigma_j}|_2) + (B_0(|X_{\sigma_j}|_2))^{1+\frac 4d }) \ve_{j+1}(\tau_{j+1})   \nonumber \\
   \leq& \ve_*(C_{\sigma_{j+1}}, D_{\sigma_{j+1}}, |X_{\sigma_j}|_2).
\end{align}

Thus, by virtue of Theorem \ref{Thm-Sta-L2},
we obtain
\begin{align}
    \|v_{j+1}\|_{V(0,\tau_{j+1})}
    \leq C(C_{\sigma_{j+1}}, D_{\sigma_{j+1}}, |X(\sigma_j)|_2),
\end{align}
and so
\begin{align} \label{esti-Xj1-L2}
   \|X\|_{V(\sigma_j,\sigma_{j+1})}
   \leq \|e^{\vf_{\sigma_j}}\|_{C([0,\tau_{j+1}];L^\9)} C(C_{\sigma_{j+1}}, D_{\sigma_{j+1}}, |X(\sigma_j)|_2).
\end{align}
This along with the inductive assumptions yields \eqref{esti-v-Xj} with $j+1$ replacing $k$.

Thus, using the inductive arguments
we conclude that \eqref{esti-v-Xj} holds for all $1\leq j\leq l$.
This yields that
\begin{align} \label{esti-X9-L2}
   \|X \|_{V(0,\sigma_l)}
    \leq \sum\limits_{k=0}^{l-1}
         \|e^{\vf_{\sigma_k}}\|_{C([0,T];L^\9)}
            C(C_T,D_T, M_T),
\end{align}
where
$C(C_T,D_T, M_T):=C(C_T,D_T, M_T,0,\sup_{0\leq x\leq M_T} B_0(x))$ is as in Theorem \ref{Thm-Sta-L2},
and $M_T$ is  as in Lemma \ref{Lem-bdd-L2}.

{\bf Step $3$.}
We  claim that
\begin{align} \label{sigma9-L2}
   \bbp (l<\9) =1.
\end{align}

To this end, we use the contradiction argument.
Suppose that \eqref{sigma9-L2} is not true.
We consider $\omega\in \{l=\9\}$.
For simplicity, we omit the argument $\omega$ below.

On one hand, by the definition of $\tau_{j+1}$,
\begin{align*}
  g_{\sigma_j}(\tau_{j+1})
  = \frac{\ve_*(C_{\sigma_{j+1}}, D_{\sigma_{j+1}}, |X(\sigma_j)|_2) } {D_0(|X(\sigma_{j})|_2)}
  \geq \frac{\ve_*(C_T, D_T, M_T) } {D_0(M_T)} >0,
\end{align*}
where
$\ve_*(C_T, D_T, M_T) := \ve_*(C_T, D_T, M_T, 0,\sup_{0\leq x\leq M_T} B_0(x))$, and
$D_0(M_T):= C_0 \sup_{0\leq x\leq M_T} ((B_0(x)) + (B_0(x))^{1+\frac 4d})$.

On the other hand,
Since the processes $t \mapsto \int_0^t g_k d\beta_k(s)$
and $t\mapsto \int_0^t g_k^2 ds$
are $(\frac 12 -\kappa)$-H\"older continuous
for any $\kappa<\frac 12$ and $1\leq k\leq N$,
we have for some positive $C(T)$ (depending on $\omega$)
$$  g_{\sigma_j} (\tau_{j+1})
   \leq C(T) (\tau_{j+1})^{\frac 12 - \kappa},\ \ \forall j\geq 1.
$$
Thus, we conclude that
$$
   \tau_{j+1} \geq \(\frac{\ve_*(C_T, D_T, M_T) } {C(T)D_0(M_T)}\)^{\frac{2}{1-2\kappa}} > 0,\ \ \forall j\geq 1.
$$
Thus, for $\omega\in \{l=\9\}$,
$$
   \sigma_l(\omega)= \sum_{j=1}^\9 \tau_{j}(\omega) = \9,
$$
which contracts the fact that $\sigma_l(\omega)\leq (T\wedge \tau^*)(\omega) \leq T<\9$.
This yields  \eqref{sigma9-L2}, as claimed.

Now, since $\{l<\9\} \subseteq \{\sigma_l = T\wedge \tau^*\}$,
combining \eqref{esti-X9-L2} and \eqref{sigma9-L2} together
we conclude that
\begin{align} \label{bdd-X-V}
   \|X\|_{V(0,\tau^*\wedge T)}
   \leq \sum\limits_{k=0}^{l-1}
         \|e^{\vf_{\sigma_k}}\|_{C([0,T];L^\9)}
        C(C_T,D_T, M_T) <\9,\ \ a.s..
\end{align}
This yields \eqref{gloexist-L2}
and so the global well-posedness of \eqref{equa-x}.

Finally,
the estimate \eqref{thm-L2-L2} follows from \eqref{bdd-X-L2},
and the estimate \eqref{thm-L2-Lpq} can be proved standardly by \eqref{bdd-X-V} and
Strichartz estimates.

Therefore, the proof of Theorem \ref{Thm-GWP} $(i)$ is complete. \hfill $\square$

\subsection{Energy-critical case} \label{Subsec-GWP-H1}

We start with the global well-posedness and scattering results
in the deterministic defocusing energy-critical case,
mainly based on the work of I-team \cite{CKSTT08},
Ryckman and Visan \cite{RV07} and Visan \cite{V07}
\footnote{As in the mass-critical case,
although the equation studied in \cite{CKSTT08,RV07,V07}
is $i\p_t u = -\Delta u + |u|^{\frac {4}{d-2}}u$,
the results there are also valid for \eqref{equa-u-H1}
by reversing the time.}

\begin{theorem} \label{Thm-H1GWP-Det} (\cite{CKSTT08,RV07,V07})
For every $u_0 \in H^1$,
there exists a unique global $H^1$-solution $u$ to the equation
\begin{align} \label{equa-u-H1}
   i\p_t u =&   \Delta u - |u|^{\frac{4}{d-2}} u, \\
   u(0) =& u_0  \nonumber
\end{align}
with $d\geq 3$.
Moreover,
\begin{align} \label{globdd-u-H1-Det}
   \|u\|_{S^1(\bbr) \cap L^2(\bbr; H^\frac 32_{-1})} \leq B_1(|u_0|_{H^1}) <\9,
\end{align}
where $ B_1(|u_0|_{H^1})$ depends continuously on $u_0$ in $H^1$,
and $u$ scatters at infinity, i.e.,
there exist $u_\pm \in H^1$ such that
\begin{align} \label{sca-det-H1}
   |e^{it\Delta} u(t) - u_\pm|_{H^1} \to 0,\ \ as\ t\to \pm\9.
\end{align}
\end{theorem}

As is the mass-critical case,
the bound of $\|u\|_{L^2(\bbr; H^\frac 32_{-1})}$
follows standardly from Strichartz estimates
and the bound of $\|u\|_{S^1(\bbr^+)}$,
and the continuous dependence on $u_0$
follows from the energy-critical stability result Lemma $3.8$ of \cite{TVZ07}.

We also need the global energy estimates below.
\begin{lemma}  \label{Lem-bdd-H1}
Assume the conditions of Theorem \ref{Thm-GWP} $(ii)$ to hold,  $3\leq d\leq 6$.
Define the Hamiltonian of $X$ by
$H(X) :=\frac 12 |\na X|_2^2 - \frac{\lbb}{\a+1} |X|^{\a+1}_{L^{\a+1}}$.
Then, for each $X_0\in H^1$, we have $\bbp$-a.s., for any $t\in (0,\tau^*)$,
\begin{align} \label{Ito-H}
   H(X(t))
    &= H(X_0) - \int_0^t \Re \int \na \ol{X} \na(\mu X) dx ds
              + \frac 12 \sum\limits_{k=1}^N \int_0^t |\na (G_k X)|^2 dx ds  \nonumber \\
            &   - \frac{\lbb(\a-1)}{2}\sum\limits_{k=1}^N \int_0^t \int (\Re G_k)^2 |X|^{\a+1} dx ds  \\
            &+ \sum\limits_{k=1}^N\int_0^t \Re \int \na \ol{X} \na (G_k X) dx d\beta_k(s)
              -\lbb \sum\limits_{k=1}^N\int_0^t \int \Re G_k |X|^{\a+1} dx d\beta_k(s). \nonumber
\end{align}
Moreover,
in the defocusing case where $\lbb = -1$,
for any
$0<T<\9$, $p\geq 1$,
\begin{align} \label{bdd-X-H1}
   \bbe \|X\|^p_{C([0;\tau^*\wedge T);H^1)} + \bbe \|X\|^p_{L^{\frac{2d}{d-2}}(0,\tau^* \wedge T;L^{\frac{2d}{d-2}})} \leq C(p,T) <\9,
\end{align}
where $\tau^*$ is the maximal existing time as in Theorem \ref{Thm-LWP}.
In particular,
\begin{align} \label{globdd-E}
     E_T:= \sup\limits_{0\leq t\leq \tau^*\wedge T} |X(t)|_{H^1} \leq C(T) <\9,\ \ a.s..
\end{align}
\end{lemma}
The proof is postponed to the Appendix.

\begin{remark} \label{Rem-globdd-E-d6-proof}
Similar formula was proved  in \cite{BRZ16} in the energy-subcritical case
(i.e., $\a\in (1,1+\frac{4}{d-2})$, $d\geq 3$),
where an approximating procedure was used to derive the It\^o formula of
corresponding potential energy $|X|^{\a+1}_{L^{\a+1}}$.
The arguments there
apply also to the energy-critical case where $3\leq d\leq 6$,
except that we need to apply the stability result Theorem  \ref{Thm-Sta-H1-dlow} instead
in the approximating procedure.
\end{remark}

\begin{remark}
In the high dimensional case where $d>6$,
the It\^o formula \eqref{Ito-H} also can be obtained
by a formal computation,
however, the rigorous derivation is technically unclear.
Actually,
one can not use the Stability result Theorem \ref{Thm-Sta-H1} to derive \eqref{Ito-H}
as in the case $3\leq d\leq 6$.
The main reason is that in \eqref{Sta-H1-ve} the time function
$g$ is imposed to satisfy a small constant,
which is even smaller than that of
the difference of two solutions \eqref{Sta-H1.1}.
Hence, the approximating procedure as in the case $3\leq d\leq 6$ in Appendix
does not work.
\end{remark}

We are now ready to prove Theorem \ref{Thm-GWP} $(ii)$.

{\bf Proof of Theorem \ref{Thm-GWP} $(ii)$.} (Energy-Critical Case).
The arguments below are similar to those in the mass-critical case,
however, based on the more delicate stability results Theorems \ref{Thm-Sta-H1-dlow} and \ref{Thm-Sta-H1}.
Below we mainly treat the case $d>6$,
the case $3\leq d\leq 6$ is easier and can be proved similarly by using Theorem \ref{Thm-Sta-H1-dlow}.

Let $X$ be the unique $H^1$-solution to \eqref{equa-x} with $\a= 1+ \frac{4}{d-2}$
on the maximal time interval $[0,\tau^*)$.
In view of Theorem \ref{Thm-LWP},
we only need to prove \eqref{gloexist-H1} for
any $0<T<\9$.

For this purpose,
we define $g$   as in Step $1$ in the proof of Theorem \ref{Thm-GWP} $(i)$
and let
$\tau_1:= \inf\{t\in (0,T\wedge \tau^*): g(t) \geq \ve_1(t)\} \wedge (T\wedge \tau^*)$
with
\begin{align}
   \ve_1(t):= \frac{\ve_*(C_t,D_t, |X_0|_{H^1})}{D_1(|X_0|_{H^1})},
\end{align}
where
$\ve_*(C_t,D'_t,|X_0|_{H^1}) := \ve_*(C_t,D'_t,\sqrt{2H(X_0)}, B_1(|X_0|_{H^1}))$ is as in Theorem \ref{Thm-Sta-H1},
and $D_1(|X_0|_{H^1}) =1+ C_1(B_1(|X_0|_{H^1}) + (B_1(|X_0|_{H^1}))^{1+\frac{4}{d-2}})$
with $C_1$ as in \eqref{esti-e1-H1} below.

We can take $g$ as the time function in Theorem \ref{Thm-Sta-H1}.
Hence,
$\sup_{0\leq t\leq \tau_1} g(t) \leq \ve_*(C_{\tau_1},D_{\tau_1}, |X_0|_{H^1})$,
and so the smallness condition on $g$ in \eqref{Sta-H1-ve} is satisfied on $[0,\tau_1]$.

Let $\vf$ be as in \eqref{vf} with $\sigma \equiv0$.
By Theorem \ref{Thm-Rescale-sigma},
$w_1:= e^{-\vf} X$ satisfies \eqref{equa-RNLS} with
$\sigma\equiv 0$ and $\a= 1+\frac{4}{d-2}$.

Moreover,
let $\wt{w}_1$ be the solution to \eqref{equa-NLS}
with $\a= 1+\frac{4}{d-2}$ and $\wt{w}(0)= w_1(0) = X_0$.
Then, Theorem \ref{Thm-L2GWP-Det} implies that
\begin{align} \label{esti-u1-L2}
   \|\wt{w}_1\|_{\bbw(0,\tau_1) \cap L^2(0,\tau_1; H^\frac 32_{-1})}
   \leq B_1 (|\wt{w}_1(0)|_{H^1}) = B_1 (|X_0|_{H^1}) <\9,
\end{align}
and by the conservation law of Hamiltonian (i.e., $H(\wt{w}(t)))= H(\wt{w}(0))$),
\begin{align} \label{wtw-H1-tau1-GWP}
    \|\wt{w}\|_{C([0,t];H^1)}
    \leq \sqrt{2\sup\limits_{0\leq s\leq t} H(\wt{w}(s))}
    = \sqrt{2H(X_0)},\ \ t\in [0,\tau_1].
\end{align}
For the error term
\begin{align*}
   e_{1}
   := - (b \cdot \na  + c)\wt{w}_{1}
      -  (1- e^{\frac{4}{d-2} \Re \vf})  F(\wt{w}_{1}),
\end{align*}
where $b,c$ are as in \eqref{b} and \eqref{c} with $\sigma \equiv 0$, respectively,
using \eqref{ineq-W.2}
and similar arguments as in the proof of \eqref{esti-e1-L2.1}
we have
\begin{align} \label{esti-e1-H1}
   &\|e_1\|_{N^1(0,\tau_1) + L^2(0,\tau_1; H^{\frac 12}_1)} \nonumber \\
   \leq& \|(b \cdot \na  + c)\wt{w}_{1}\|_{L^2(0,\tau_1; H^\frac 12_1)}
         + \|(1- e^{\frac{4}{d-2}\Re \vf}) F(\wt{w}_{1}) \|_{L^2(0,\tau_1; W^{1,\frac{2d}{d+2}})} \nonumber \\
   \leq& C \sup\limits_{0\leq t\leq \tau_1} g(t)
           (\|\wt{w}_1\|_{L^2(0,\tau_1; H^\frac 32_{-1})}
        + \|\wt{w}_1\|^{\frac {d+2} {d-2}}_{\bbw(0,\tau_1)}) \nonumber \\
   \leq& C_1 (B_1(|X_0|_{H^1}) + (B_1(|X_0|_{H^1}))^{1 + \frac{4}{d-2}}) \ve_1(\tau_1) \nonumber \\
   \leq& \ve_*(C_{\tau_1}, D'_{\tau_1}, |X_0|_{H^1}).
\end{align}
Then,
applying Theorem \ref{Thm-Sta-H1}
we obtain
\begin{align} \label{esti-v1-H1}
    \|w_1\|_{\bbw(0,\tau_1)} \leq C(C_{\tau_1}, D'_{\tau_1}, \sqrt{2H(X_0)}, B_1(|X_0|_{H^1})) =:C(C_{\tau_1}, D'_{\tau_1}, |X_0|_{H^1}),
\end{align}
which implies that
\begin{align} \label{bdd-X-bbw-tau1}
    \|X\|_{\bbw(0,\tau_1)}
    \leq \|e^\vf\|_{C([0,\tau_1];W^{1,\9})} C(C_{\tau_1}, D'_{\tau_1}, |X_0|_{H^1}).
\end{align}
Thus, \eqref{gloexist-H1} follows from \eqref{bdd-X-bbw-tau1} if $\tau_1\geq \tau^*$.
Otherwise,  we turn to the next step.

Let $\tau_j, \sigma_j, g_{\sigma_{j}}$ and $l$
be as in Step 2 in the proof of Theorem \ref{Thm-GWP} $(i)$,
but with
$$\ve_{j+1}(t) = \frac{\ve_*(C_{\sigma_j+t},D'_{\sigma_j+t}, |X(\sigma_j)|_{H^1})}{D_1(|X(\sigma_j)|_{H^1})}, $$
where $D_1(|X(\sigma_j)|_{H^1})$ is defined similarly to $D_1(|X_0|_{H^1})$,
$\ve_*(C_{\sigma_j+t},D'_{\sigma_j+t}, |X(\sigma_j)|_{H^1})
   := \ve_*(C_{\sigma_j+t},D'_{\sigma_j+t}, \sqrt{2|X(\sigma_j)|_{H^1}}, B_1(|X(\sigma_j)|_{H^1}))$
is as in Theorem \ref{Thm-Sta-H1}.

We use the inductive arguments to prove that for any $1\leq j\leq l$,
\begin{align} \label{bdd-X-bbw-H1}
     \|X\|_{\bbw(0,\sigma_j)}
     \leq \sum\limits_{k=0}^{j-1}
          \|e^{\vf_{\sigma_k}}\|_{C([0,\tau_{k+1}];W^{1,\9})}
          C(C_{\sigma_{k+1}},D'_{\sigma_{k+1}}, |X(\sigma_k)|_{H^1})<\9,
\end{align}
where $\vf_{\sigma_k}$ are as in Theorem \ref{Thm-Rescale-sigma}
with $\sigma_k$ replacing $\sigma$,
and
$C(C_{\sigma_{k+1}},D'_{\sigma_{k+1}}, |X(\sigma_k)|_{H^1})$
$:= C(C_{\sigma_{k+1}},D'_{\sigma_{k+1}}, \sqrt{2H(X(\sigma_k))}, B_1(|X(\sigma_k)|_{H^1}))$
are as in Theorem \ref{Thm-Sta-H1}.

We see from \eqref{bdd-X-bbw-tau1} that \eqref{bdd-X-bbw-H1} holds for $j=1$.
Suppose that \eqref{bdd-X-bbw-H1} holds for each $1\leq k\leq j<l$.

In order to prove \eqref{bdd-X-bbw-H1} with $j+1$ replacing $j$,
we consider the rescaling transformation
$w_{j+1}(t) := e^{-\vf_{\sigma_j}(t)} X(\sigma_j+t)$, $0\leq t < (T\wedge \tau^*)-\sigma_j$,
and apply Theorem \ref{Thm-Rescale-sigma} to obtain
\begin{align} \label{equa-wj1}
   i \p_t w_{j+1} &=   e^{-\vf_{\sigma_j}}\Delta (e^{\vf_{\sigma_j}}w_{j+1})
                 - e^{\frac{4}{d-2} \Re \vf_{\sigma_j}} F(w_{j+1}) ,  \\
    w_{j+1}(0)&= X(\sigma_j).               \nonumber
\end{align}
Then, we compare \eqref{equa-wj1} with the equation
\begin{align} \label{equa-wtw-j1.0-H1}
    i\p_t \wt{w}_{j+1}
    =  \Delta \wt{w}_{j+1} - F(\wt{w}_{j+1}),
\end{align}
or equivalently,
\begin{align} \label{equa-wtw-j1-H1}
    i\p_t \wt{w}_{j+1}
    =   e^{-\vf_{\sigma_j}}\Delta (e^{\vf_{\sigma_j}}\wt{w}_{j+1})
                -  e^{\frac{4}{d-2} \Re \vf_{\sigma_j}} F(\wt{w}_{j+1})   + e_{j+1},
\end{align}
with $\wt{w}_{j+1}(0)  =w_{j+1}(0) = X(\sigma_j)$
and the error term
\begin{align*}
   e_{j+1}
   = - (b_{\sigma_j}\cdot \na + c_{\sigma_j}) \wt{w}_{j+1}
     - (1-e^{\frac{4}{d-2}\Re \vf_{\sigma_j}}) F(\wt{w}_{j+1}).
\end{align*}
Theorem \ref{Thm-H1GWP-Det} yields that
$\wt{w}_{j+1}$ exists globally and satisfies
\begin{align}
    \|\wt{w}_{j+1}\|_{\bbw(0,\tau_{j+1})} \leq B_1 (|\wt{w}_{j+1}(0)|_{H^1})
    = B_1(|X(\sigma_j)|_{H^1}),
\end{align}
and similarly to  \eqref{wtw-H1-tau1-GWP},
\begin{align}
    \|\wt{w}_{j+1}\|_{C([0,\tau_{j+1}];H^1)}
    \leq \sqrt{2\sup\limits_{0\leq s\leq \tau_{j+1}} H(\wt{w}_{j+1}(s))}
    = \sqrt{2H(X(\sigma_j))}.
\end{align}

Note that,
we can take $g_{\sigma_j}$ the time function in Theorem \ref{Thm-Sta-H1}.
Then, by the definition of $\tau_{j+1}$,
$\sup_{0\leq t\leq \tau_{j+1}} g_{\sigma_j}(t) \leq \ve_*(C_{\sigma_{j+1}},D'_{\sigma_{j+1}},|X({\sigma_{j+1}})|_{H^1})$,
and so the smallness condition on $g_{\sigma_j}$ in Theorem \ref{Thm-Sta-H1} is satisfied on $[0,\tau_{j+1}]$.

Moreover, similarly to \eqref{esti-e1-H1},
\begin{align} \label{esti-ej1}
   &\|e_{j+1}\|_{L^{2}(0,\tau_{j+1}; H^{\frac 12}_1) + N^1(0,\tau_{j+1})} \nonumber \\
   \leq& C_1 \sup\limits_{0\leq t\leq \tau_{j+1}} g_{\sigma_j} (t)
         (\|\wt{w}_{j+1}\|_{L^2(0,\tau_{j+1}; H^\frac 32 _1)}
         +  \|\wt{w}_{j+1}\|^{\frac {d+2}{d-2}}_{\bbw(0,\tau_{j+1})}) \nonumber \\
   \leq& C_1(B_1(|X(\sigma_j)|_{H^1}) + (B_1(|X(\sigma_j)|_{H^1}))^{\frac{d+2}{d-2}}) \ve_{j+1}(\tau_{j+1}) \nonumber \\
   \leq& \ve_*(C_{\sigma_{j+1}},D'_{\sigma_{j+1}},|X(\sigma_j)|_{H^1}).
\end{align}
Thus, by virtue of Theorem \ref{Thm-Sta-H1}, we obtain
\begin{align}
    \|w_{j+1}\|_{\bbw(0,\tau_{j+1} )}
    \leq C(C_{\sigma_{j+1}},D'_{\sigma_{j+1}},|X(\sigma_j)|_{H^1}),
\end{align}
and so
\begin{align*}
    \|X\|_{\bbw(\sigma_j, \sigma_{j+1})}
    \leq  \|e^{ \vf_{\sigma_j}}\|_{C([0,\tau_{j+1}];W^{1,\9})}
          C(C_{\sigma_{j+1}},D'_{\sigma_{j+1}},|X(\sigma_j)|_{H^1})<\9,
\end{align*}
thereby yielding  \eqref{bdd-X-bbw-H1} with $j+1$ replacing $j$.

Therefore, the inductive arguments yield
\eqref{bdd-X-bbw-H1}  for all $1\leq j\leq l$
and so
\begin{align}   \label{bdd-X-W1}
    \|X\|_{\bbw(0,\sigma_l)}
   \leq& \sum\limits_{k=0}^{l-1}
        \|e^{\vf_{\sigma_k}}\|_{C([0,\tau_{k+1}];W^{1,\9})}
         C(C_{\sigma_{j+1}},D'_{\sigma_{j+1}}, |X(\sigma_j)|_{H^1}) \nonumber \\
   \leq& \sum\limits_{k=0}^{l-1}
         \|e^{\vf_{\sigma_k}}\|_{C([0,\tau_{k+1}];W^{1,\9})}
          C(C_T,D'_T, C' E_T, \sup\limits_{0\leq x\leq E_T} B_1(x)),
\end{align}
where $E_T$ is as in \eqref{bdd-E-assum-d6},
and we also used the inequality
$\sup_{0\leq t<\tau^*} \sqrt{2H(X(t))} \leq C'E_T$ in the last step,
implied by the Sobolev embedding.

Since by \eqref{bdd-E-assum-d6},
$E_T<\9$, a.s.,
and for any $0\leq j\leq l-1$,
\begin{align*}
   g_{\sigma_j}(\tau_{j+1})
   =\frac{\ve_*(C_{\sigma_{j+1}},D'_{\sigma_{j+1}},|X(\sigma_j)|_{H^1})}{D_1(|X(\sigma_j)|_{H^1})}
   \geq \frac{\ve_*(C_T,D_T',C'E_T, \sup\limits_{0\leq x\leq E_T}B_1(x))}{ \sup\limits_{0\leq x\leq E_T} D_1(x)}
   >0,
\end{align*}
we can use similar arguments as in Step $3$
in the proof of Theorem \ref{Thm-GWP} $(i)$ to
deduce that $l<\9$, a.s.,
which together with \eqref{bdd-X-W1}
yields the global bound \eqref{gloexist-H1},
thereby implying the global existence of $X$ to \eqref{equa-x}.

Finally, the estimate \eqref{thm-H1-H1} follows from Lemma \ref{Lem-bdd-H1},
and \eqref{thm-H1-LpWq} follows from  \eqref{bdd-X-W1} and Strichartz estimates.
The proof of Theorem \ref{Thm-GWP}  is complete.
\hfill $\square$

\section{Scattering} \label{Sec-Sca}

In this section we prove the scattering behavior of global solutions to \eqref{equa-x}.

The idea here is based on the very recent work \cite{HRZ18}.
More precisely, we use a new rescaling transformation \eqref{z*},
i.e., $z_*= e^{-\vf_*}X$
with $\vf_*$ as in \eqref{vf*},
and compare the resulting random equation \eqref{equa-RNLS-Sca}
with \eqref{equa-NLS}
but after some large time $T$, i.e.,
\begin{align} \label{equa-NLS-Sca}
   i\p_t u =& \Delta u- |u|^{\a-1}u, \\
   u(T) =& z_*(T).   \nonumber
\end{align}
Let us start with the mass-critical case.

\subsection{Mass-critical case}

First we enhance the bounds \eqref{bdd-X-L2} and \eqref{globdd-M}
to the whole time regime,
under the condition that $g_k\in L^2(\bbr^+)$, a.s..

\begin{lemma} \label{Lem-globdd-X-L2}
Consider the situations in Theorem \ref{Thm-Sca} $(i)$.
Then, for any $p\geq 1$,
\begin{align} \label{globdd-EX-L2}
   \bbe \sup\limits_{0\leq t<\9}  |X(t)|_2^p \leq C(p) <\9.
\end{align}
In particular, we have the global pathwise  bound
\begin{align} \label{globdd-X-L2}
   M_\9:= \sup\limits_{0\leq t<\9} |X(t)|_2  \leq C  <\9,\ \ a.s..
\end{align}
\end{lemma}

{\bf Proof.} Estimate \eqref{globdd-EX-L2} can be proved by using
the It\^o formula \eqref{Ito-L2}
and similar arguments as in the proof of \cite[(1.7)]{HRZ18}.
We omit the details here for simplicity. \hfill $\square$ \\

Below we also have the important uniform boundedness (independent of $T$).

\begin{lemma} \label{Lem-globdd-u-L2}
For $z_*(T) \in L^2$,
there exists a unique global $L^2$-solution $u$ (depending on $T$)
to \eqref{equa-NLS-Sca} with $\a=1+\frac 4d$, $d\geq 1$,
which scatters at infinity and satisfies
\begin{align} \label{globdd-u-Sca-L2}
     \|u\|_{S^0(T,\9) \cap L^2(T,\9;H^\frac 12_{-1})}  \leq C <\9,\ \ a.s.,
\end{align}
where $C$ is independent of $T$.
\end{lemma}

{\bf Proof.}
For each $z_*(T)\in L^2$ fixed,
the global well-posedness and scattering follow from Theorem \ref{Thm-L2GWP-Det}.

Regarding \eqref{globdd-u-Sca-L2},
applying the global-in-time Strichartz estimates  in Theorem \ref{Thm-Stri*} to \eqref{equa-NLS*}
and using the H\"older inequality \eqref{ineq-V}
we have that for any $t>T$,
\begin{align*}
  \|u\|_{L^2(T,t;H^\frac 12_{-1})} + \|u\|_{S^0(T,t)}
  \leq C|u(T)|_2 + C \|u\|_{V(T,t)}^{1+\frac 4d}
\end{align*}
with $C$ independent of $T$ and $t$.
using \eqref{globdd-u-L2-Det}
and that $u(T)=z_*(T)$, we get
\begin{align*}
   \|u\|_{L^2(T,t;H^\frac 12_{-1})} + \|u\|_{S^0(T,t)}
  \leq C|z_*(T)|_2
       + C  (B_0(|z_*(T)|_2))^{1+\frac 4d}.
\end{align*}
Since $g_k\in L^2(\bbr^+)$, a.s.,
we have $\vf_* \in C(\bbr^+; L^\9)$,
and so
$|z_*(T)|_2 \leq C|X(T)|_2$ with
$C$ independent of $T$.
In view of the global bound \eqref{globdd-X-L2},
we obtain
\begin{align*}
    \|u\|_{L^2(T,t;H^\frac 12_{-1})} + \|u\|_{S^0(T,t)}
    \leq CM_\9 + C \sup\limits_{0\leq x\leq M_\9} (B_0(x))^{1+\frac 4d}
    <\9
\end{align*}
with $C, M_\9$  independent of $T$ and $t$.
Thus, letting $t\to \9$, we obtain  \eqref{globdd-u-Sca-L2}.
\hfill $\square$

The following result
is crucial for the scattering in the mass-critical case.

\begin{lemma} \label{Lem-z*u-0-L2}
Consider the situations in Theorem \ref{Thm-Sca} $(i)$.
Let $u$ be the solution to \eqref{equa-NLS-Sca}
with $u(T) = z_*(T)$.
Then,  $\bbp$-a.s.
as  $T \to \9$,
\begin{align} \label{z*u-0-scaL2}
     \|z_*- u\|_{S^0(T,\9) \cap L^2(T,\9; H^{\frac 12}_{-1})} \to 0.
\end{align}
\end{lemma}

{\bf Proof.}
We use the idea of comparison as in the proof of Theorem \ref{Thm-GWP} $(i)$.
Precisely,
we compare the solution $z_*$ to \eqref{equa-RNLS-Sca}
with the solution $u$ to \eqref{equa-NLS-Sca}.

For this purpose,
we rewrite \eqref{equa-NLS-Sca} with $\a=1+\frac 4d$ as follows
\begin{align*}
   i\p_t u =   e^{-\vf_*}\Delta (e^{\vf_*} u)
                 - e^{\frac 4d \Re \vf_*} F(u) + e,
\end{align*}
with the error term
\begin{align*}
    e = -(b_*\cdot \na + c_*) u
        - (1-e^{\frac 4d \Re \vf_*}) F(u).
\end{align*}

Since $g_k\in L^2(\bbr^+)$, a.s.,
$1\leq k\leq N$,
for any multi-index $\g$, as $T\to \9$,
\begin{align} \label{part-vf*}
  \sup\limits_{t\geq T} |\p_x^\g \vf_*(t,x)|
  \leq C \<x\>^{-2} \sup\limits_{t\geq T}  \sum\limits_{k=1}^N
       \( \bigg|\int_t^\9 g_k d\beta_k \bigg| + \int_t^\9 g_k^2 ds \) \to 0.
\end{align}
Hence, for $T$ large enough,
Theorem \ref{Thm-Stri*} yields that
global-in-time Strichartz and local smoothing estimates hold for the operator
$e^{-\vf_*}\Delta(e^{\vf_*} \cdot)$.

Note that, for any $t\geq T$,
\begin{align} \label{esti-u-scaL2.0}
  \|e\|_{N^0(T,t) + L^2(T,t;H^{-\frac 12}_{1})}
    \leq& \|(b_*\cdot \na + c_*) u\|_{L^2(T,t;H^{-\frac 12}_{1})} \nonumber \\
        & +  \|(e^{\frac 4d \Re \vf_*}-1) F(u)\|_{L^{\frac{2d+4}{d+4}}((T,t)\times \bbr^d)}.
\end{align}

Since $g_k\in L^2(\bbr^+$, a.s.,
we have
\begin{align} \label{ve1-u-scaL2}
   \ve_1(T)
   := \sup\limits_{t\geq T}  \sum\limits_{k=1}^N
       \( \bigg|\int_t^\9 g_k d\beta_k \bigg| + \int_t^\9 g_k^2 ds \)  \to 0, \ \ as\ T\to \9, \ a.s..
\end{align}
Then, estimating as in \eqref{esti-e1-L2.1}, we get
\begin{align} \label{esti-u-scaL2.1}
    \|(b_*\cdot \na + c_*)u\|_{L^2(T,t;H^{-\frac 12}_{1})}
    \leq C \ve_1(T) \|u\|_{L^2(T,t;H^{\frac 12}_{-1})},
\end{align}
where $C$ is independent of $T$ and $t$.

Moreover,
using again $g_k\in L^2(\bbr^+)$, a.s.,
we deduce that
\begin{align} \label{ve2-u-scaL2}
     \ve_2 (T): = \sup\limits_{t\geq T} \|\vf_*(t)\|_{W^{1,\9}} \to 0,\ \ as\ T\to \9,\ a.s..
\end{align}
Then, using the inequality $|e^x-1|\leq e|x|$ for $|x|\leq 1$ and \eqref{ineq-V},
we get
\begin{align}   \label{esti-u-scaL2.2}
   \|(e^{\frac 4d \Re \vf_*}-1)  F(u) \|_{N^0(T,t)}
   \leq C \ve_2(T) \|u\|_{V(T,t)}^{1+\frac 4d}.
\end{align}

Plugging \eqref{esti-u-scaL2.1} and \eqref{esti-u-scaL2.2} into \eqref{esti-u-scaL2.0}
and using \eqref{globdd-u-Sca-L2} we obtain that,
\begin{align} \label{esti-e-Sca-L2}
    \|e\|_{N^0(T,t) + L^2(T,t;H^{-\frac 12}_{1}) }
     \leq C \ve(T) (\|u\|_{ L^2(T,t;H^{\frac 12}_{-1})} + \|u\|_{V(T,t)}^{1+\frac 4d})
    \leq C \ve(T),
\end{align}
where $\ve(T): = \max\{\ve_1(T), \ve_2(T)\}$,
and $C$ is independent of $T, t$,
due to  \eqref{globdd-u-Sca-L2}.

Thus,
in view of Remark \ref{Rem-glob-Sta-L2}
and \eqref{globdd-u-Sca-L2},
we obtain that for $T$ large enough,
\begin{align} \label{esti-z*-u-Sca-L2}
     \|z_* - u\|_{S^0(T,t) \cap  L^2(T,t;H^{\frac 12}_{-1})} \leq C \ve(T),
\end{align}
where   $C$  is independent of $T$ and $t$.
(Note that, since $|z_*(T)-u(T)|_2 =0$,
we can choose $M'=\ve(T)$ when applying
the stability result.)

Therefore, letting $t\to \9$ in \eqref{esti-z*-u-Sca-L2}
and using \eqref{ve1-u-scaL2}, \eqref{ve2-u-scaL2}
we obtain \eqref{z*u-0-scaL2}.
\hfill $\square$ \\

{\bf Proof of Theorem \ref{Thm-Sca} }
$(i)$ { (Mass-Critical Case).}
Let $u$ be as in Lemma \ref{Lem-globdd-u-L2}.
We have $\bbp$-a.s. for any $t_1, t_2 \geq T$,
\begin{align*}
   |e^{it_1\Delta} z_*(t_1) - e^{it_2\Delta} z_*(t_2) |_2
   \leq& |e^{it_1\Delta} (z_*-u)(t_1) - e^{it_2\Delta} (z_*-u)(t_2) |_2 \\
       & + |e^{it_1\Delta} u(t_1) - e^{it_2\Delta} u(t_2) |_2.
\end{align*}
By lemma \ref{Lem-globdd-u-L2}, the scattering of $u$ yields
\begin{align*}
    |e^{it_1\Delta} u(t_1) - e^{it_2\Delta} u(t_2) |_2 \to 0,\ \ as\ t_1, t_2 \to \9.
\end{align*}
Hence, taking into account \eqref{z*u-0-scaL2} we obtain
\begin{align*}
   \limsup\limits_{t_1,t_2\to \9}
   |e^{it_1\Delta} z_*(t_1) - e^{it_2\Delta} z_*(t_2) |_2
   \leq& \limsup\limits_{t_1,t_2\to \9}|e^{it_1\Delta} (z_*-u)(t_1) - e^{it_2\Delta} (z_*-u)(t_2) |_2  \\
   \leq& 2\|z_*-u\|_{C([T,\9);L^2)} \to 0, \ \ as\ T\to \9, \ a.s..
\end{align*}
This implies that $\{e^{it\Delta} z_*(t)\}$  is a Cauchy sequence in $L^2$,
thereby yielding \eqref{Sca-L2.1}.

Next we prove \eqref{Sca-L2.2}.
Recall that $U_*(t,s)$, $s,t\geq 0$,
are the evolution operators related to the operators $e^{-\vf_*}\Delta(e^{\vf_*} \cdot)$, $t\geq 0$.
Then, by Equation \eqref{equa-RNLS-Sca},
\begin{align*}
   z_*(t) = U_*(t,0)X_0
            + i \int_0^t U_*(t,s)e^{\frac 4d \Re\vf_*}  F(z_*) ds,\ \ t\geq 0.
\end{align*}
Since $U_*(0,t) U_*(t,s) = U_*(0,s)$ for $s\geq 0$,
applying $U_*(0,t)$ to both sides we get
\begin{align*}
    U_*(0,t)z_*(t)
    =X_0 + i \int_0^t U_*(0,s) e^{\frac 4d \Re\vf_*}   F(z_*)ds,
\end{align*}
which implies that for any $0<t_1<t_2<\9$,
\begin{align} \label{equa-z*-Sca-L2}
   U_*(0,t_2) z_*(t_2) - U_*(0,t_1)z_*(t_1)
   =&  i \int_{t_1}^{t_2}  U_*(0,s) e^{\frac 4d \Re\vf_*}   F(z_*)ds \nonumber  \\
   =& U_*(0,t_2) \(i  \int_{t_1}^{t_2}  U_*(t_2,s) e^{\frac 4d \Re\vf_*}   F(z_*) ds\) \nonumber  \\
   =:& U_*(0,t_2) w(t_2),
\end{align}
Thus,
applying homogeneous Strichartz estimates in Theorem \ref{Thm-Stri*}
leads to
\begin{align} \label{V-t12.1}
  | U_*(0,t_2) z_*(t_2) - U_*(0,t_1)z_*(t_1) |_2
  \leq \|U_*(\cdot,t_2) w(t_2)\|_{C([0,t_2];L^2)}
  \leq C|w(t_2)|_2,
\end{align}
where $C$ is independent of $t_1, t_2$.
Moreover, since
$w(\cdot)$ satisfies \eqref{equa-RNLS}
with the initial datum $w(t_1)=0$,
applying Theorem \ref{Thm-Stri*} again and using \eqref{ineq-V} we obtain
\begin{align}  \label{V-t12.2}
   |w(t_2)|_2
   \leq \|w\|_{C([t_1,t_2];L^2)}
   \leq  C \|e^{\frac 4d \Re\vf_*}   F(z_*) \|_{L^{\frac{2d+4}{d+4}}((t_1, t_2)\times \bbr^d)}
   \leq C \|z_*\|^{1+\frac 4d}_{V(t_1,t_2)},
\end{align}
where $C$ is independent of $t_1,t_2$,
due to the global-in-time Strichartz estimates
and that $e^{\frac 4d \Re \vf_*} \in C(\bbr^+; L^\9)$, a.s..

Moreover, taking into account \eqref{globdd-u-Sca-L2} and \eqref{z*u-0-scaL2} we have
\begin{align} \label{globdd-z*-sca-L2}
   \|z_*\|_{V(T,\9)}
   \leq \|u\|_{V(T,\9)} +  \|z_*-u\|_{V(T,\9)} <\9,\ \ a.s..
\end{align}

Thus, plugging \eqref{V-t12.2} into \eqref{V-t12.1}
and using the global bound \eqref{globdd-z*-sca-L2}
we obtain
\begin{align*}
   |U_*(0,t_2) z_*(t_2) - U_*(0,t_1)z_*(t_1)|_2 \leq C\|z_*\|_{V(t_1,t_2)}^{1+\frac 4d} \to 0,\ \ as\ t_1,t_2\to \9,\ a.s..
\end{align*}
This implies that
$\{U_*(0,t)z_*(t)\}$ is a Cauchy sequence in $L^2$,
and so \eqref{Sca-L2.2} follows.
The proof of Theorem \ref{Thm-Sca} $(i)$ is complete.

\subsection{Energy-critical and pseudo-conformal cases}

As in the mass-critical case, we have the global bound of
solutions  $X$  and $u$ below.

\begin{lemma} (\cite{HRZ18}) \label{Lem-globdd-X-H1}
Consider the situations in Theorem \ref{Thm-Sca} $(ii)$
with $3\leq d\leq 6$.
Then, for any $p\geq 1$,
\begin{align} \label{globdd-EX-H1}
   \bbe \sup\limits_{0\leq t<\9}  |X(t)|_{H^1}^p + |X(t)|^p_{L^\frac{2d}{d-2}} \leq C(p) <\9.
\end{align}
In particular,
\begin{align} \label{globdd-X-H1}
   E_{\9}:= \sup\limits_{0\leq t<\9} |X(t)|_{H^1} \leq C <\9,\ \ a.s..
\end{align}
\end{lemma}

\begin{lemma} \label{Lem-globdd-u-H1}
For  $z_*(T) \in H^1$,
there exists a unique global $H^1$-solution $u$ (depending on $T$)
to \eqref{equa-NLS-Sca} with $\a= 1+\frac{4}{d-2}$, $d\geq 3$,
which scatters at infinity and satisfies
\begin{align} \label{globdd-u-Sca-H1}
     \|u\|_{S^1(T, \9) \cap L^2(T,\9; H^{\frac 32}_{-1})} \leq C <\9,\ \ a.s.,
\end{align}
where $C$ is independent of $T$.
\end{lemma}

{\bf Proof.}
The proof is analogous to that of Lemma \ref{Lem-globdd-u-L2}.
First,
the global well-posedness and scattering in the space $H^1$
follows from Theorem \ref{Thm-H1GWP-Det}.

In order to prove \eqref{globdd-u-Sca-H1},
applying Theorem \ref{Thm-Stri*} to \eqref{equa-NLS-Sca}
and using \eqref{ineq-W.2} and \eqref{globdd-u-H1-Det} we obtain
that for any $t>T$,
\begin{align*}
   \|u\|_{L^2(T,t; H^\frac 32_{-1})}
   + \|u\|_{S^1(T,t)}
   \leq& C|u(T)|_{H^1} + C\|u\|_{\bbw(T,t)}^{\frac{d+2}{d-2}} \\
   \leq&  C|z_*(T)|_{H^1} + C(B_1(\|z_*(T)\|_{H^1}))^{\frac{d+2}{d-2}},
\end{align*}
where $C$ is independent of $T$ and $t$.
Since
$|z_*(T)|_{H^1} \leq \|e^{\vf_*}\|_{C(\bbr^+; W^{1,\9})} |X(T)|_{H^1}$,
using \eqref{globdd-X-H1} and letting $t\to \9$ we prove \eqref{globdd-u-Sca-H1}.
\hfill $\square$

We have the crucial asymptotics of difference between the solutions $z_*$ and $u$.

\begin{lemma} \label{Lem-z*u-0-H1}
Consider the situations in Theorem \ref{Thm-Sca} $(ii)$.
Let $u$ be the solution to \eqref{equa-NLS-Sca} with $\a= 1+\frac {4}{d-2}$,
$u(T) = z_*(T)$, $d\geq 3$.
Then,
\begin{align} \label{z*u-0-scaH1}
     \|z_*- u\|_{S^1(T,\9) \cap L^2(T,\9; H^\frac 32_{-1})} \to 0, \ \ as\ T\to \9,\ a.s..
\end{align}
\end{lemma}

{\bf Proof.}
The case $3\leq d\leq 6$ was proved in the recent work \cite{HRZ18}
under weaker condition on $\phi_k$.
Below we mainly consider the high dimensional case  $d>6$.

We shall apply Theorem \ref{Thm-Sta-H1} to compare the solutions $z_*$ and $u$.
For this purpose, we reformulate Equation \eqref{equa-NLS-Sca} as follows
\begin{align*}
   i \p_t u =   e^{-\vf_*}\Delta (e^{\vf_*} u)
                 -  e^{\frac{4}{d-2}\Re \vf_*}   F(u) + e,
\end{align*}
with the error term
\begin{align*}
    e = -(b_*\cdot \na + c_*) u
        - (1-e^{\frac{4}{d-2}\Re \vf_*}) F(u).
\end{align*}

We see that,
\eqref{part-vf*} implies that for $T$ large enough,
\begin{align*}
   \sup\limits_{t\geq T} g(t)
   := \sup\limits_{t\geq T} \sum\limits_{k=1}^N
      \(\bigg| \int_t^\9 g_k d\beta_k \bigg| + \int_t^\9 g_k^2 ds\)
   \leq \ve,
\end{align*}
so the smallness condition on $g$ in Theorem \ref{Thm-Sta-H1}
is satisfied on $[T,t]$ for any $t>T$.

Regarding the error term,
similarly to \eqref{esti-e-Sca-L2},
\begin{align*}
   \|e\|_{ N^1(T,t) + L^2(T,t;H^{\frac 32}_{-1})}
   \leq  C \ve(T)
         ( \|u\|_{L^2(T,t;H^{\frac 32}_{-1})}
           + \|u\|^{\frac{d+2}{d-2}}_{\bbw(T,t)})
   \leq C \ve(T) \to 0,\ as\ T\to \9.
\end{align*}
where $C$ is independent of $T$ and $t$ due to \eqref{globdd-u-Sca-H1},
and $\ve(T)$ is as in \eqref{esti-e-Sca-L2}.

Thus, by virtue of  Theorem \ref{Thm-Sta-H1}
we deduce that
there exist  $c, C>0$, independent of $T$ and $t$, such that
\begin{align*}
  \|z_* - u \|_{S^1(T,t) \cap L^2(T,t;H^{\frac 32}_{-1})}
  \leq C (\ve(T))^c.
\end{align*}
Therefore,
letting $t\to \9$ and then taking $T \to \9$ we obtain \eqref{z*u-0-scaH1}.
\hfill $\square$  \\

Now, we are ready to prove Theorem \ref{Thm-Sca} $(ii)$.

{\bf Proof of Theorem \ref{Thm-Sca} $(ii)$.}
In the case where $3\leq d\leq 6$,
because of the global well-posedness of \eqref{equa-x}   and the estimate \eqref{thm-H1-LpWq} in
Theorem \ref{Thm-GWP} $(ii)$,
Assumption $(H0')$ in \cite{HRZ18} in the energy-critical case is satisfied.
Thus, the asymptotics \eqref{Sca-L2.1} and \eqref{Sca-L2.2} with $H^1$ replacing $L^2$
follow from  Theorem $1.4$ in the recent work \cite{HRZ18}.

Below we consider the case where $d>6$.
The proof is similar to that of mass-critical case,
thanks to Lemmas \ref{Lem-globdd-X-H1}, \ref{Lem-globdd-u-H1} and \ref{Lem-z*u-0-H1}.

Actually, let $u$ be as in Lemma \ref{Lem-z*u-0-H1}. We have
for any $t_1,t_2 \geq T$,
\begin{align*}
   |e^{it_1\Delta} z_*(t_1) - e^{it_2\Delta} z_*(t_2) |_{H^1}
   \leq  |e^{it_1\Delta} u(t_1) - e^{it_2\Delta} u(t_2) |_{H^1}
        + 2 \| z_*-u \|_{C([T,\9);H^1)}.
\end{align*}
Then, we first use the scattering of $u$ in Lemma \ref{Lem-globdd-u-H1} to
pass to the limits $t_1,t_2 \to \9$,
and then we use Lemma \ref{Lem-z*u-0-H1} to take the limit $T\to \9$.
It follows that
\begin{align*}
   \limsup\limits_{t_1,t_2\to \9}
   |e^{it_1\Delta} z_*(t_1) - e^{it_2\Delta} z_*(t_2) |_{H^1}
   \leq  2 \| z_*-u \|_{C([T,\9);H^1)} \to 0,\ as\ T\to \9,\ a.s..
\end{align*}
This yields that $\{e^{it\Delta} z_*(t)\}$ is a Cauchy sequence in $H^1$,
thereby implying \eqref{Sca-L2.1} with $H^1$ replacing $L^2$.

Moreover, applying Theorem \ref{Thm-Stri*} to \eqref{equa-z*-Sca-L2}
with $\frac{4}{d-2}$ replacing $\frac 4d$ we get
\begin{align*}
   |U_*(0,t_2)z_*(t_2) - U_*(0,t_1)z_*(t_1)|_{H^1}
   \leq C \|e^{\frac{4}{d-2} \Re \vf_*}F(z_*)\|_{N^1(t_1,t_2)}
   \leq C \| z_*\|^{\frac{d+2}{d-2}}_{\bbw(t_1,t_2)} ,
\end{align*}
where $C$ is independent of $t_1$ and $t_2$.

Then, taking into account the global pathwise bound of $z_*$
implied by \eqref{globdd-u-Sca-H1} and \eqref{z*u-0-scaH1},
we pass to the limits $t_1,t_2\to\9$ to
obtain that right-hand side above tends to $0$ almost surely.
This yields that
$\{U_*(0,t)z_*(t)\}$ is a Cauchy sequence in $H^1$,
thereby implying \eqref{Sca-L2.2}  with $H^1$ replacing $L^2$.

Therefore, the proof of Theorem \ref{Thm-Sca} $(ii)$  is complete. \hfill $\square$\\

{\bf Proof of Theorem \ref{Thm-Sca} $(iii)$.}
In view of the  global well-posedness of \eqref{equa-x},
we see that Assumption $(H0')$ in \cite{HRZ18} in the energy-critical case is satisfied.
Thus, the asymptotic  \eqref{Sca-L2.1}  with $\Sigma$ replacing $L^2$
follows  from  Theorem $1.3$ in the recent work \cite{HRZ18}. \hfill $\square$

\begin{remark}
In $(H0')$ of \cite{HRZ18},
the assumption on the boundedness of $\||\cdot|X\|_{L^\g(0,T; L^\rho)}$
is redundant,
which can be deduced from Theorem $1.2$ and Lemma $2.1$ of \cite{HRZ18}.
\end{remark}

We close this section with the proof of Theorem \ref{Thm-S0S1-Global}.

{\bf Proof of Theorem \ref{Thm-S0S1-Global}.}
In the mass-critical case,
the global bound \eqref{globbd-L2-Lpq} follows from
\eqref{thm-L2-Lpq} and Lemma \ref{Lem-globdd-u-L2}.
In the energy-critical case,
the global bound \eqref{globdd-H1-LpWq}
is a consequence of \eqref{thm-H1-LpWq} and
Lemma \ref{Lem-globdd-u-H1}.
\hfill $\square$

\section{Support theorem}  \label{Sec-Supp}

In this section we prove Theorem \ref{Thm-Supp}
concerning the support theorem for \eqref{equa-x}.
We combine the idea of \cite{MS94}
with the stability results in Section \ref{Sec-Sta}.

Recall that, for any $h=(h_1,\cdots,h_N)\in \mathscr{H}$
(i.e., the Cameron-Martin space),
$X(\beta+h)$ denotes the solution to \eqref{equa-x}
with the driving processes $\beta_k+h_k$
replacing the Brownian motions $\beta_k$,
$1\leq k\leq N$,
$S(h)$ denotes the controlled solution to \eqref{equa-Sh}.
The global existence and uniqueness of $X(\beta+h)$ and $S(h)$ can be proved similarly as in Section \ref{Sec-GWP}.

In view of Proposition $2.2$ in \cite{MS94},
we only need to prove that, for any $\ve >0$,
\begin{align}
    & \lim\limits_{n\to \9} \bbp (\|S(\beta^n) - X(\beta)\|_{\calx(0,T)}\geq \ve) = 0,   \label{0-Sn-X-E}\\
    & \lim\limits_{n\to \9} \bbp (\|X(\beta^n-\beta+h) - S(h)\|_{\calx(0,T)}\geq \ve)=0, \label{0-Xnh-Sh-E}
\end{align}
where
$\calx(0,T) = S^0(0,T)\cap L^2(0,T; H^\frac 12_{-1})$
or
$\calx(0,T) = S^1(0,T)\cap L^2(0,T; H^\frac 32_{-1})$
in the mass-critical or energy-critical case, respectively,
and $\beta^n$ are the adapted linear interpolation of   Brownian motions in \cite{MS94}, defined by
\begin{align*}
   \beta^n(t) = \beta(\ol{t}_n) + 2^n(t-\tilde{t}_n)(\beta_{\tilde{t}_n} - \beta_{\ol{t}_n}),
\end{align*}
$\wt{t}_n = \frac{k}{2^n}$ and $\ol{t}_n = \frac{k-1}{2^n} \vee 0$ if $\frac{k}{2^n} \leq t< \frac{k+1}{2^n}$.

For this purpose,
we first prove the asymptotic result below.

\begin{lemma} \label{Lem-betan-beta}
Assume  $g_k$ are deterministic and continuous, $1\leq k\leq N$.
Then,
\begin{align} \label{E-betan-beta-0}
   \bbe \(\sup\limits_{t\in [0,T]}
   \bigg|\int_0^t g_k(s) \dot\beta_k^n(s) ds - \int_0^t g_k(s) d\beta_k(s) \bigg|^2 \)
   \to 0, \ \ as\ n\to \9.
\end{align}
\end{lemma}

The proof is quite technical and is postponed to the Appendix. \\

{\bf Proof of Theorem \ref{Thm-Sca}.}
We mainly prove Theorem \ref{Thm-Supp} in the energy-critical case when $3\leq d\leq 6$.
The mass-critical case can be proved similarly,
based on the stability result Theorem \ref{Thm-Sta-L2}.

It order to obtain \eqref{0-Sn-X-E} and \eqref{0-Xnh-Sh-E},
it is  equivalent to prove that  for any subsequence $\{n_j\}$,
there exists some subsequence $\{n_{j_k}\}$ of $\{n_j\}$
such that as $k\to \9$,
\begin{align}
    & \lim\limits_{k\to \9} \|S(\beta^{n_{j_k}}) - X(\beta)\|_{S^1(0,T) \cap L^2(0,T; H^\frac 32_{-1})} = 0,\ \ a.s., \label{0-Sn-X}\\
    & \lim\limits_{k\to \9} \|X(\beta^{n_{j_k}}-\beta+h) - S(h)\|_{S^1(0,T)\cap L^2(0,T; H^\frac 32_{-1})} =0,\ \ a.s..  \label{0-Xnh-Sh}
\end{align}
Below, for simplicity,
we still denote the subsequence $\{n_j\}$ by $\{n\}$.

{\bf Proof of \eqref{0-Sn-X}}.
For this purpose,
for each $h\in H^1(0,T; \bbr^N)$,
we set
$$\psi(\beta+ h)(t)
   := \sum\limits_{k=1}^N \int_0^t G_k(s,x) d \beta_k(s)
      + \sum\limits_{k=1}^N \int_0^t G_k(s,x) \dot{h}_k(s) ds - \int_0^t \wh{\mu}(s,x) ds, $$
and $\psi(h)$ is defined similarly.
Using
\begin{align} \label{res-zn-Sbetan}
   z_n:=e^{-\psi(\beta^n)} S(\beta^n),
\end{align}
we have
\begin{align} \label{equa-zn-Supp}
   i\p_t z_n =&   e^{-\psi(\beta^n)} \Delta(e^{\psi(\beta^n)} z_n)
          -e^{\frac{4}{d-2} \Re \psi(\beta^n)} F(z_n), \\
     z_n(0)=& X_0. \nonumber
\end{align}
Similarly,
$\wt{z}:= e^{-\psi(\beta)} X(\beta)$ satisfies the random  equation
\begin{align}
   i\p_t \wt{z} =&   e^{-\psi(\beta)} \Delta(e^{\psi(\beta)} \wt{z}) dt
             - e^{\frac {4}{d-2} \Re \psi(\beta)} F(\wt{z}), \label{equa-wtz-Supp} \\
   \wt{z}(0)=& X_0, \nonumber
\end{align}
or equivalently,
\begin{align*}
   i \p_t \wt{z} =&   e^{-\psi(\beta^n)} \Delta(e^{\psi(\beta^n)} \wt{z}) dt
             - e^{\frac {4}{d-2} \Re \psi(\beta^n)} F(\wt{z}) + e_n
\end{align*}
with the error term
\begin{align} \label{en-betan}
   e_n =&  ( (b(\psi(\beta))-b(\psi(\beta^n)))\cdot \na + (c(\psi(\beta))-c(\psi(\beta^n))) ) \wt{z}  \nonumber \\
        & - (e^{\frac{4}{d-2} Re \psi(\beta)} - e^{\frac{4}{d-2} \Re \psi(\beta^n)}) F(\wt{z}),
\end{align}
where
$b(\psi(\beta)) = 2 \na \psi(\beta)$,
$c(\psi(\beta)) = \Delta\psi(\beta) + \sum_{j=1}^d (\p_j \psi(\beta))^2$,
and $b(\psi(\beta^n))$, $c(\psi(\beta^n))$
are defined similarly.

Note that,
the global well-posedness and
the bound of $S^1(0,T)\cap L^2(0,T;H^{\frac 32}_{-1})$-norm of $\wt{z}$
can be proved similarly as in Section \ref{Sec-GWP}
by using the energy-critical stability result Theorem $1.4$ of \cite{TV05}.
Similar assertions also hold for $z_n$, $n\geq 1$.

Let
\begin{align*}
    & g^n(t): =\sum_{k=1}^N  \int_0^t g_k(s)\dot{\beta}^n_k(s) ds   + \int_0^t g_k^2(s) ds,
\end{align*}
and define $g(t)$  similarly with $d\beta_k(s)$ replacing $\dot{\beta}^n_k(s) ds$.
Set
\begin{align*}
      \ve_n(t) := \sum_{k=1}^N  \bigg| \int_0^t g_k \dot{\beta}_k^n(s) ds - \int_0^t g_k d\beta_k(s) \bigg|.
\end{align*}

Lemma \ref{Lem-betan-beta} implies that
for some subsequence of $\{n\}$ (still dented by $\{n\}$),
$\bbp$-a.s., as $n\to\9$,
\begin{align} \label{gn-g-0}
    \sup\limits_{0\leq t\leq T} \|\psi(\beta^n) - \psi(\beta)\|_{W^{1,\9}}
    \leq C\sup\limits_{0\leq t\leq T} |g^n(t) - g(t)|
    \leq C \sup\limits_{0\leq t\leq T} \ve_n(t)    \to 0.
\end{align}
In particular, for some positive constant $C$ uniformly bounded of $n$,
\begin{align} \label{bdd-psin-gn}
    \sup\limits_{0\leq t\leq T} \|\psi(\beta^n)\|_{W^{1,\9}} \leq C, \ \
   \sup\limits_{n\geq 1} \sup\limits_{0\leq t\leq T} |g^n(t)|\leq C.
\end{align}
This along with  Assumption $(H0)$ yields that for any multi-index $\g$,
\begin{align} \label{part-psibeta}
  \sup\limits_{n\geq 1} \sup\limits_{0\leq t\leq T}  |\p_x^\g \psi(\beta^n)| \leq
            C\<x\>^{-2} \sup\limits_{n\geq 1} \sup\limits_{0\leq t\leq T}  g^n(t)
  \leq C\<x\>^{-2}.
\end{align}

Then, Theorem \ref{Thm-Stri} yields that Strichartz and local smoothing estimates hold
for the operator $e^{-\psi(\beta^n)} \Delta (e^{\psi(\beta^n)} \cdot )$,
and the corresponding Strichartz constants $C_T$
are uniformly bounded for all $n$.

Estimating as in \eqref{esti-e1-H1}
and using the global bounds of the
$L^2(0,T;H^{\frac 32}_{-1})$- and $S^1(0,T)$-norms of $\wt{z}$,
we obtain
\begin{align*}
   \|e_n\|_{N^1(0,T) + L^2(0,T;H^{\frac 12}_{1})}
   \leq  C(T) \sup\limits_{0\leq t\leq T} \ve_n(t)
        ( \|\wt{z}\|_{L^2(0,T;H^{\frac 32}_{-1})}  + \|\wt{z}\|_{\bbw(0,T)}^{1+\frac{4}{d-2}})
   \leq  C'(T) \ve_n(t) .
\end{align*}
We note that $C'(T)$ is independent of $n$,
due to the uniform bound \eqref{bdd-psin-gn}.
This along with \eqref{gn-g-0} yields that
\begin{align} \label{e-b-bn-0}
   \|e_n\|_{N^1(0,T) + L^2(0,T;H^{\frac 12}_{1})}
   \to 0,\ \ as\ n\to\9,\ a.s..
\end{align}

Then,
by virtue of Theorem \ref{Thm-Sta-H1-dlow}, we obtain
that $\bbp$-a.s. as $n\to \9$,
\begin{align} \label{zn-z-0}
      \| z_n - z \|_{S^1(0,T) \cap L^2(0,T;H^{\frac 32}_{-1})}
      \leq C(T) \sup\limits_{0\leq t\leq T} \ve_n(t)
      \to 0.
\end{align}
In particular, this yields the uniform bound
\begin{align*}
   \sup\limits_{n\geq 1} \|z_n\|_{S^1(0,T) \cap L^2(0,T; H^{\frac 32}_{-1})}
   \leq C(T) <\9,\ \ a.s..
\end{align*}

We claim that \eqref{zn-z-0} implies \eqref{0-Sn-X}.
Actually,  we have
\begin{align} \label{Sbetan-Xbeta}
  &\|S(\beta^n) - X(\beta)\|_{S^1(0,T) \cap L^2(0,T;H^{\frac 32}_{-1})}  \nonumber  \\
  =& \|e^{\psi(\beta^n)} z_n - e^{\psi(\beta)} \wt{z}\|_{S^1(0,T) \cap L^2(0,T;H^{\frac 32}_{-1})}\\
  \leq& \|e^{\psi(\beta^n)} (z_n-\wt{z})\|_{S^1(0,T) \cap L^2(0,T;H^{\frac 32}_{-1})}
        + \|(e^{\psi(\beta^n)} - e^{\psi(\beta)}) \wt{z}\|_{S^1(0,T) \cap L^2(0,T;H^{\frac 32}_{-1})} .   \nonumber
\end{align}
In order to pass to the limit $n\to \9$,
using  \eqref{bdd-psin-gn} and the inequality
$|e^x-e^y| \leq C|e^y||x-y|$ for $|x|,|y|\leq \frac 12$,
we have,
as $n\to \9$,
\begin{align}   \label{zn-z-S1-0}
  & \|e^{\psi(\beta^n)} (z_n-\wt{z})\|_{S^1(0,T) }
        + \|(e^{\psi(\beta^n)} - e^{\psi(\beta)}) \wt{z}\|_{S^1(0,T)} \nonumber   \\
  \leq& C(T) (\| (z_n-\wt{z})\|_{S^1(0,T) }
        + \|\psi(\beta^n) - \psi(\beta)\|_{C([0,T]; W^{1,\9})} )  \nonumber  \\
  \leq& C(T) \(\| (z_n-\wt{z})\|_{S^1(0,T) }
               + \sup\limits_{0\leq t\leq T} \ve_n(t)  \)
  \to 0,
\end{align}
where in the last step we used \eqref{gn-g-0} and \eqref{zn-z-0}.

Regarding the $L^2(0,T; H^{\frac 32}_{-1})$-norm in \eqref{Sbetan-Xbeta},
we deduce from \eqref{part-psibeta} that
$e^{\psi(\beta^n)} \in S^0$,
and so $\Psi_{p}:= \<x\>^{-1}\<\na\>^{\frac 32} e^{\psi(\beta^n)}\<\na\>^{-\frac 32}  \<x\>$
is a pseudo-differential operator of order zero.
This along with Lemma \ref{Lem-L2-Bdd} yields
\begin{align*}
   \|e^{\psi(\beta^n)} (z_n-\wt{z})\|_{L^2(0,T;H^{\frac 32}_{-1})}
   =&  \| \Psi_p  \<x\>^{-1} \<\na\>^{\frac 32} (z_n-\wt{z})\|_{L^2(0,T;L^2)}   \\
   \leq&  C \|(z_n-\wt{z})\|_{L^2(0,T;H^{\frac 32}_{-1})}
\end{align*}
Moreover,
using Assumption $(H0)$
and the inequality
$|e^x-e^y| \leq C|e^y||x-y|$ for $|x|,|y|\leq \frac 12$
we have that
for any multi-index $\g$,
\begin{align*}
   \p_x^\g |(e^{\psi(\beta^n)} - e^{\psi(\beta)})(t,x)|
   \leq C(T) \<x\>^{-2}   \ve_n(t)
\end{align*}
where $C(T)$ is independent of $n$. Then, estimating as above we get
\begin{align*}
   \|(e^{\psi(\beta^n)} - e^{\psi(\beta)}) \wt{z}\|_{L^2(0,T;H^{\frac 32}_{-1})}
   \leq C(T)
         \sup\limits_{0\leq t\leq T} \ve_n(t)
          \| \wt{z}\|_{L^2(0,T;H^{\frac 32}_{-1})}
   \leq C(T)
         \sup\limits_{0\leq t\leq T} \ve_n(t) .
\end{align*}
Combing the  estimates above together we conclude that,
$\bbp$-a.s. as $n \to \9$,
\begin{align} \label{zn-z-LS-0}
  &\|e^{\psi(\beta^n)} (z_n-\wt{z})\|_{L^2(0,T;H^{\frac 32}_{-1})}
        + \|(e^{\psi(\beta^n)} - e^{\psi(\beta)}) \wt{z}\|_{L^2(0,T;H^{\frac 32}_{-1})}  \nonumber \\
  \leq&  C  (\|(z_n-\wt{z})\|_{L^2(0,T;H^{\frac 32}_{-1})}
         + \sup\limits_{0\leq t\leq T} \ve_n(t)  )
  \to 0.
\end{align}
Therefore,
plugging \eqref{zn-z-S1-0} and \eqref{zn-z-LS-0} into \eqref{Sbetan-Xbeta}
we prove \eqref{0-Sn-X}, as claimed.  \\

{\bf Proof of \eqref{0-Xnh-Sh}.}
The proof is similar as above.
Now we use a new transformation
\begin{align} \label{res-yn-Xbetan}
   y_n= e^{-\psi(\beta^n-\beta+h)} X(\beta^n-\beta+h)
\end{align}
to obtain
\begin{align} \label{equa-zn-bbnh}
   i\p_t y_n =  e^{-\psi(\beta^n-\beta+h)} \Delta (e^{\psi(\beta^n-\beta+h)}y_n)
           - e^{\frac{4}{d-2} \Re \psi(\beta^n-\beta+h)} F(y_n)
\end{align}
with $y_n(0)=X_0$.
Arguing as  above,
we have the Strichartz and local smoothing estimates for
the operator $e^{-\psi(\beta^n-\beta+h)} \Delta (e^{\psi(\beta^n-\beta+h)} \cdot)$,
the related Strichartz constants $C_T$
are uniformly bounded of $n$.

Moreover, letting $\wt{y}:= e^{-\psi(h)} S(h)$ we have
\begin{align}
   i \p_t \wt{y} =&   e^{-\psi(h)} \Delta (e^{-\psi(h)} \wt{y}) - e^{\frac{4}{d-2}\Re \psi(h)}  F(\wt{y}) \nonumber \\
       =& e^{-\psi(\beta^n-\beta+h)} \Delta (e^{\psi(\beta^n-\beta+h)}\wt{y})
           - e^{\frac{4}{d-2}\Re \psi(\beta-\beta^n+h)}  F(\wt{y}) + e'_n
\end{align}
with $\wt{y}(0) = X_0$ and the error term
\begin{align} \label{en-supp.1}
   e'_n =& (( b(\psi(h)) - b(\psi(\beta^n-\beta+h)) ) \cdot \na
         + ( c(\psi(h)) - c(\psi(\beta^n-\beta+h)) )) \wt{y} \nonumber \\
        & -(e^{\frac{4}{d-2} \Re \psi(h)} - e^{\frac{4}{d-2} \Re \psi(\beta^n - \beta + h)})  F(\wt{y}),
\end{align}
where
$b(\psi(h))= 2\na \psi(h)$,
$c(\psi(h)) = \Delta \psi(h) + \sum_{j=1}^d(\p_j \psi(h))^2$,
and $b(\psi(\beta^n-\beta +h))$, $c(\psi(\beta^n-\beta+h))$ are defined similarly.

Similarly to \eqref{gn-g-0},
for some subsequence  of $\{n\}$ (still denoted by $\{n\}$),
for any multi-index $\g$,
\begin{align*}
    \sup\limits_{0\leq t\leq T} |\p_x^\g (\psi(\beta^n-\beta+h)(t) - \psi(h)(t))|
   \leq C \<x\>^{-2} \sup\limits_{0\leq t\leq T} \ve_n(t),
\end{align*}
and
\begin{align*}
    \sup\limits_{0\leq t\leq T} \|\psi(\beta^n-\beta+h)(t) - \psi(h)(t)\|_{W^{1,\9}}
    \leq C \sup\limits_{0\leq t\leq T} \ve_n(t).
\end{align*}
Then,
similarly to \eqref{en-supp.1},
we have that
\begin{align*}
   \|e'_n\|_{N^1(0,T)+ L^2(0,T;H^{\frac 12}_{1})} \to 0,\ \  as\ n \to \9,\ a.s..
\end{align*}
which along with Theorem \ref{Thm-Sta-H1-dlow} implies  that $\bbp$-a.s.
as $n\to \9$,
\begin{align} \label{yn-wty-0}
   \|y_n -\wt{y}\|_{S^1(0,T) \cap L^2(0,T;H^{\frac 32}_{-1})}
   \to 0.
\end{align}
In particular,
\begin{align*}
   \sup\limits_{n\geq 1} \|y_n\|_{S^1(0,T) \cap L^2(0,T; H^{\frac 32}_{-1})} \leq C(T) <\9,\ \ a.s..
\end{align*}

Thus, estimating as those below \eqref{zn-z-0}
and using \eqref{yn-wty-0} we obtain \eqref{0-Xnh-Sh}.
Therefore, the proof of Theorem \ref{Thm-Supp} is complete.
\hfill $\square$

\section{Appendix} \label{Sec-App}

{\bf Proof of Theorem \ref{Thm-Rescale-sigma}.}
The case where $\sigma \equiv 0$ can be proved similarly as in
\cite[Lemma 6.1]{BRZ14} and \cite[Lemma 2.4]{BRZ16}
in the $L^2$ and $H^1$ space, respectively.
For the  general case,
we  prove the $L^2$ case below,
the $H^1$ case can be proved similarly.

Set $\vf(t) := \vf_0(t) =  \int_0^{t} G_k d\beta_k(s) - \int_0^t \wh{\mu}ds$
and $v(t) := v_0(t)=e^{-\vf(t)}X(t)$, $t\in[0,\tau^*)$.
For any $0\leq t<\tau^*-\sigma$, we have
\begin{align} \label{X-v}
     X(t) = e^{\vf(t)} v(t),\ \
     X(\sigma+t) = e^{\vf_\sigma(t)} v_\sigma(t),
\end{align}
and
\begin{align} \label{vfsigma-vf}
            \vf(\sigma+t) - \vf(\sigma) = \vf_\sigma(t).
\end{align}
It follows that
\begin{align} \label{vsigma-v}
    v_\sigma(t) = e^{-\vf_\sigma(t)} X(\sigma+t)
                = e^{-(\vf(\sigma+t) - \vf(\sigma) )} X(\sigma+t)
                = e^{\vf(\sigma)} v(\sigma+t).
\end{align}
Then, similar arguments as in the proof of \cite[Lemma 6.1]{BRZ14} show that,
$v$ satisfies pathwisely the equation \eqref{equa-RNLS} on $[0,\tau^*)$
in the space $H^{-2}$,
with $0$ replacing $\sigma$.
Hence, $\bbp$-a.s. for any $t\in[0,\tau^*-\sigma)$,
\begin{align*}
    iv(\sigma+t)
    = iv(\sigma) + \int_\sigma^{\sigma+t} e^{-\vf(s)} \Delta (e^{\vf(s)}v(s)) ds
                 - \int_\sigma^{\sigma+t} e^{(\a-1)\Re \vf(s)} F(v(s)) ds,
\end{align*}
where the equation is taken in $H^{-2}$.
Plugging this into  \eqref{vsigma-v} yields that
\begin{align} \label{vsigma-v.0}
    iv_\sigma(t)
    =& ie^{\vf(\sigma)}v(\sigma)
       +\int_\sigma^{\sigma+t} e^{\vf(\sigma)-\vf(s)} \Delta (e^{\vf(s)}v(s)) ds  \nonumber \\
     &  -  \int_\sigma^{\sigma+t}  e^{\vf(\sigma)} e^{(\a-1)\Re \vf(s)} F(v(s))ds  \nonumber  \\
    =& iX(\sigma)
       +  \int_0^{t}   e^{\vf(\sigma)-\vf(\sigma+s)} \Delta (e^{\vf(\sigma+s)}v(\sigma+s)) ds  \nonumber \\
    &   -  \int_0^{t}   e^{\vf(\sigma)} e^{(\a-1)\Re \vf(\sigma+s)} F(v(\sigma+s)) ds.
\end{align}

Note that, by \eqref{X-v} and \eqref{vfsigma-vf},
\begin{align*}
   e^{\vf(\sigma)-\vf(\sigma+s)} \Delta (e^{\vf(\sigma+s)}v(\sigma+s))
   = e^{-\vf_\sigma(s)} \Delta (X(\sigma+s))
   = e^{-\vf_\sigma(s)} \Delta (e^{\vf_\sigma(s)}v_\sigma(s) ).
\end{align*}
Moreover,
\begin{align*}
     e^{\vf(\sigma)} e^{(\a-1)\Re \vf(\sigma+s)}  F(v(\sigma+s))
     =&  e^{\vf(\sigma) - \vf(\sigma+s)}  F(X(\sigma+s)) \nonumber \\
     =& e^{-\vf_\sigma(s)} F(e^{\vf_\sigma(s)}v_\sigma(s)) \nonumber \\
     =& e^{(\a-1)\Re \vf_\sigma(s)}F(v_{\sigma}(s)).
\end{align*}

Thus, plugging the two identities above into \eqref{vsigma-v.0},
we obtain that $\bbp$-a.s. for any $0\leq t<\tau^*- \sigma$,
\begin{align*}
    i v_\sigma(t)
    = i X(\sigma) + \int_0^t e^{-\vf_\sigma(s)} \Delta (e^{\vf_{\sigma(s)}} v_\sigma(s)) ds
      - \int_0^t e^{(\a-1)\Re \vf_\sigma(s)}  F(v_{\sigma}(s)) ds
\end{align*}
as an equation in $H^{-2}$,
which implies \eqref{equa-RNLS},
thereby finishing the proof.
\hfill $\square$ \\

{\bf Proof of Lemma \ref{Lem-bdd-H1}.}
Below, we mainly prove the It\^o formula \eqref{Ito-H}.
The estimate \eqref{bdd-X-H1} can be obtained from \eqref{Ito-H}
by using similar arguments as in the proof of \cite[(2.4)]{BRZ18},
involving the Burkholder-Davis-Gundy inequality and the Gronwall inequality.

In order to prove \eqref{Ito-H},
we use the stability result to pass to the limit
in the approximating procedure
as in the proof of $(2.4)$ in \cite{BRZ16}.

More precisely,
we consider the solution $X_m$ to \eqref{equa-x},
with the nonlinearity $\Theta_m (|X|^{4/(d-2)}X)$ replacing $|X|^{4/(d-2)}X$,
where $\Theta_m f:= \mathscr{F}^{-1}(\theta(\frac{|\cdot|}{m})) \ast f$,
$\theta \in C_c^\9$ is real-valed, nonnegative,
and $\theta(x)=1$ for $|x|\leq 1$,
$\theta(x)=0$ for $|x|\geq 2$.

Since the operators $\{\Theta_m\}$ are uniformly bounded in $L^p$ for any $1<p<\9$
(see \cite[(3.2)]{BRZ16}),
arguing as in the proof of \cite[Theorem 3.1]{BRZ16},
we deduce that  $X_m$, $m\geq 1$, exist on the common
time regime $[0,\tau^*)$,
and
\begin{align*}
   \sup_{m\geq 1}\|X_m\|_{S^1(0,t) \cap L^2(0,t; H^{\frac 32}_{-1})} \leq C(t)<\9, \ \ t\in (0,\tau^*),\  a.s..
\end{align*}
Moreover, similar arguments as in the proof of \cite[(3.9)]{BRZ16} yield that
\begin{align*}
    &H(X_m(t)) \nonumber \\
    =& H(X_0) - \int_0^t \Re \int \na \ol{X}_m \na(\mu X_m) dx ds
              + \frac 12 \sum\limits_{k=1}^N \int_0^t |\na (G_k X_m)|^2 dx ds  \nonumber \\
            &   - \frac{\lbb(\a-1)}{2} \sum\limits_{k=1}^N \int_0^t \int (\Re G_k)^2 |X_m|^{\a+1} dx ds  \\
            & - \lbb \int_0^t \Re \int i \na ((\Theta_m -1) F(X_m) ) \na \ol{X}_m dx ds \\
            & + \sum\limits_{k=1}^N\int_0^t \Re \int \na \ol{X}_m \na (G_k X_m) dx d\beta_k(s)
              -\lbb \sum\limits_{k=1}^N \int_0^t \int \Re G_k |X_m|^{\a+1} dx d\beta_k(s). \nonumber
\end{align*}

Then, in order to pass to the limit $m\to \9$,
we only need to show that, for $w_m :=e^{-\vf} X_m$ and $w :=e^{-\vf} X$
with $\vf$ as in \eqref{vf} with $\sigma \equiv 0$,
\begin{align} \label{asym-ym-y}
    w_m \to w,\ \ in\ S^1(0,t), \ \  as\ m\to \9,\ t\in (0,\tau^*),\ a.s..
\end{align}

For this purpose,
we use the stability result in Section \ref{Sec-Sta}
to replace  the subcritical arguments in \cite{BRZ16}.
Note that, $w_m$ satisfies
\begin{align} \label{equa-ym}
   i\p_t w_m =  e^{-\vf}\Delta (e^{\vf}w_m)  - e^{-\vf} \Theta_m(F(e^\vf z_m))
\end{align}
with $w_m(0) = X_0$.
Moreover,  $w$ satisfies \eqref{equa-w-p}
with $\vf$ replacing $\psi$, i.e.,
\begin{align} \label{equa-z}
   i\p_t w =   e^{-\vf}\Delta (e^{\vf}w)  -  e^{-\vf} \Theta_m(F(e^\vf z)) + e_m
\end{align}
with the error
\begin{align*}
    e_m =   e^{-\vf} ( \Theta_m(F(e^\vf w))-F(e^\vf w)).
\end{align*}
Since for   $p\in (1,\9)$,
$\Theta_m f \to f$ in $L^p$ (see \cite[(3.3)]{BRZ16}),
we have for  $t\in (0,\tau^*)$,
\begin{align*}
    \|e_m\|_{L^{q2}(0,t; W^{1,\frac{2d}{d+2}})}
    \leq C(t) \| \Theta_m(F(e^\vf w))-F(e^\vf w)\|_{L^{2}(0,t; W^{1,\frac{2d}{d+2}})}
    \to 0,\ m\to \9,
\end{align*}
where $C(t)$ is independent of $m$.

Therefore,
we deduce that
the asymptotic \eqref{asym-ym-y} holds
by using the stability result similar to Theorem \ref{Thm-Sta-H1-dlow}
with the nonlinearity $e^{-\vf}\Theta_m (F(e^{\vf}w_m))$
replacing $e^{\frac{4}{d-2}\Re \Phi}F(w)$,
which can be proved similarly as in the proof of Theorem \ref{Thm-Sta-H1-dlow}.
Then,
we use \eqref{asym-ym-y} to pass to the limit $m\to \9$
in the It\^o formula of $H(X_m)$ to obtain \eqref{Ito-H}.
The proof is complete.
\hfill $\square$\\

{\bf Proof of Lemma \ref{Lem-betan-beta}.}
Note that, for each $1\leq k\leq N$ fixed,
\begin{align} \label{gk-betanbeta}
   \bigg|\int_0^t g_k(s) \dot\beta_k^n(s) ds - \int_0^t g_k(s) d\beta_k(s) \bigg|
   \leq& \bigg|\int_{\frac{[2^nt]}{2^n}}^t g_k(s) d\beta_k(s)\bigg|
      + \bigg|\int_{\frac{[2^nt]}{2^n}}^t g_k(s) \dot{\beta}^n_k(s)ds \bigg| \nonumber \\
      +& \bigg|\int_0^{\frac{[2^nt]}{2^n}} g_k(s) \dot\beta^n_k(s) ds
     - \int_0^{\frac{[2^nt]}{2^n}} g_k(s) d\beta_k(s) \bigg| \nonumber \\
   =:& J'_{n,1}(t) + J_{n,2}'(t) + J'_{n,3}(t).
\end{align}
Below we estimate $J'_{n,1}, J'_{n,2}, J'_{n,3}$ respectively.

First we prove that
\begin{align} \label{Jn1-0}
   \bbe \sup\limits_{0\leq t\leq T} (J'_{n,1}(t))^2
   \to 0,\ \ as\ n\to \9.
\end{align}
To this end, we set
$M_k(t):= \int_0^t g_k(s) d\beta_k(s)$.
Since $g_k\in C(0,T)$,
using the Burkholder-Davis-Gundy inequality
we have that for any $p\geq 1$,
\begin{align*}
   \bbe|M_k(t) - M_k(s)|^{2p} \leq C(p) |t-s|^p.
\end{align*}
Then, in view of
Kolmogorov's continuity criterion
(see, e.g., \cite[Proposition 2.1]{MS94}), we get that for any
$\lbb>0$, $\g<\frac{2p}{2p-1}$,
\begin{align*}
     \bbp \( \sup\limits_{t\not =s} \frac{|M_k(t) - M_k(s)|}{|t-s|^\g} >\lbb \) \leq C \lbb^{-2p}.
\end{align*}
In particular, taking $p=3$ and $\g=\frac 14$, we  arrive at
\begin{align*}
    \bbp(\sup\limits_{0\leq t\leq T} \bigg|M_k(t) - M_k(\frac{[2^n t]}{2^n}) \bigg|  >\lbb )
    \leq C \lbb^{-6} 2^{-\frac 32 n}.
\end{align*}
This yields that
\begin{align*}
   \bbe \sup\limits_{0\leq t\leq T} (J'_{n,1}(t))^2
   =& 2 \int_0^\9 \lbb \bbp \(\sup\limits_{0\leq t\leq T} |M_k(t) - M_k(\frac{[2^n t]}{2^n})| >\lbb \)  d\lbb \\
   \leq& \frac 2n + 2 \int_{\frac 1n}^\9\lbb^{-5} 2^{-\frac 32 n}
   =  \frac 2n + \frac 12 n^4 2^{-\frac 32n} \to 0,\ \ as\ n \to \9,
\end{align*}
which implies \eqref{Jn1-0}, as claimed.

Similarly, since $g_k\in C(0,T)$, $0\leq t- \frac{[2^nt]}{2^n} \leq \frac{1}{2^n}$,
\begin{align*}
   J'_{n,2}(t)
   =& \bigg| 2^n \int_{\frac{[2^nt]}{2^n}}^t g_k(s) ds \(\beta_k(\frac{[2^n]t}{2^n}) - \beta_k(\frac{[2^n]t-1}{2^n})\)\bigg| \\
   \leq& C \bigg|\beta_k(\frac{[2^n]t}{2^n}) - \beta_k(\frac{[2^n]t-1}{2^n})\bigg|.
\end{align*}
Arguing as above we have
\begin{align} \label{Jn2-0}
    \bbe \sup\limits_{0\leq t\leq T} (J'_{n,2}(t))^2
   \to 0,\ \ as\ n\to 0.
\end{align}

It remains to prove that
\begin{align} \label{Jn3-0}
    \bbe \sup\limits_{0\leq t\leq T} (J'_{n,3}(t))^2
   \to 0,\ \ as\ n\to 0.
\end{align}

For this purpose,
since
\begin{align*}
   \int_{\frac{j-1}{2^n}}^{\frac{j}{2^n}} g_k(s)  \dot{\beta}^n_k(s) ds
   = \int_{\frac{j-2}{2^n}}^{\frac{j-1}{2^n}} \(\int_{\frac{j-1}{2^n}}^{\frac{j}{2^n}} g_k(r)2^n dr \) d\beta_k(s),
\end{align*}
we have
\begin{align*}
  J'_{n,3}(t)
  =& \bigg| \sum\limits_{j=1}^{[2^nt]} \int_{\frac{j-2}{2^n}}^{\frac{j-1}{2^n}}
     \(\int_{\frac{j-1}{2^n}}^{\frac{j}{2^n}} g_k(r) 2^n dr\) d\beta_k(s)
    - \sum\limits_{j=1}^{[2^nt]}  \int_{\frac{j-1}{2^n}}^{\frac{j}{2^n}} g_k(s) d\beta_k(s)\bigg| \nonumber \\
  =& \bigg|\sum\limits_{j=1}^{[2^nt]-1} \int_{\frac{j-1}{2^n}}^{\frac{j}{2^n}}
      \( \(\int_{\frac{j}{2^n}}^{\frac{j+1}{2^n}} g_k(r) 2^n dr\) -g_k(s)\) d\beta_k(s)
      -  \int_{\frac{[2^nt]-1}{2^n}}^{\frac{[2^nt]}{2^n}} g_k(s) d\beta_k(s)\bigg|.
\end{align*}
Let
$ M_n(t):= \sum_{j=1}^{[2^nt]-1} \int_{\frac{j-1}{2^n}}^{\frac{j}{2^n}}
      ( (\int_{\frac{j}{2^n}}^{\frac{j+1}{2^n}} g_k(r) 2^n dr) -g_k(s)) d\beta_k(s)$.
We get
\begin{align} \label{Jn3.1}
   J'_{n,3}(t) \leq |M_t(t)| + \bigg|  \int_{\frac{[2^nt]-1}{2^n}}^{\frac{[2^nt]}{2^n}} g_k(s) d\beta_k(s) \bigg|.
\end{align}

Estimating as in the proof of \eqref{Jn1-0} we see that
\begin{align} \label{Jn3.2}
     \bbe \sup\limits_{0\leq t\leq T} \bigg|\int_{\frac{[2^nt]-1}{2^n}}^{\frac{[2^nt]}{2^n}}g_k(s) d\beta_j(s) \bigg|^2
     \to 0,\ \ as\ n\to \9.
\end{align}
Moreover, since $g_k(s)$ is deterministic,
using the independence of increments of Brownian motions
we have that for each $n\geq 1$,
$t \mapsto M_n(t)$ is a right-continuous martingale.
Then, using the maximal inequality and
the  Burkholder-Davis-Gundy inequality we get
\begin{align} \label{esti-Mn}
    \bbe \sup\limits_{0\leq t\leq T} |M_n(t)|^2
    \leq& 4 \bbe |M_n(T)|^2   \nonumber \\
    \leq&  C  \sum\limits_{j=1}^{[2^nT]-1} \int_{\frac{j-1}{2^n}}^{\frac{j}{2^n}}
      \( \(\int_{\frac{j}{2^n}}^{\frac{j+1}{2^n}} g_k(r) 2^n dr\) -g_k(s)\)^2 ds.
\end{align}
For any $\ve>0$,
by the uniform continuity of $g_k$ on $[0,T]$,
we have that for $n$ large enough,
$|g_k(r_1) - g_k(r_2)| \leq \ve$
for any $|r_1-r_2| \leq 2^{1-n}$.
Then, by the mean-value theorem for integrals,
we get that for any $1\leq j\leq [2^nT]-1$,
\begin{align*}
   \bigg| \(\int_{\frac{j}{2^n}}^{\frac{j+1}{2^n}} g_k(r) 2^n dr\) -g_k(s) \bigg|
   \leq |g_k(s_{n,j}) - g_k(s)| \leq \ve,
\end{align*}
where $s_{n,j} \in (\frac{j}{2^n}, \frac{j+1}{2^n})$.
Thus the right hand-side of \eqref{esti-Mn} is bounded by
\begin{align*}
   C  \sum\limits_{j=1}^{[2^nT]-1} \int_{\frac{j-1}{2^n}}^{\frac{j}{2^n}} \ve^2 ds
   \leq C T \ve^2.
\end{align*}
This implies that
\begin{align} \label{Jn3.3}
   \bbe \sup\limits_{0\leq t\leq T} |M_n(t)|^2 \to 0,\ \ as\ n\to \9.
\end{align}
Thus, we obtain \eqref{Jn3-0} from \eqref{Jn3.2} and \eqref{Jn3.3}.

Therefore, collecting \eqref{gk-betanbeta}, \eqref{Jn1-0}, \eqref{Jn2-0} and \eqref{Jn3-0} together
we prove   \eqref{E-betan-beta-0}.
\hfill $\square$\\

\subsection*{Acknowledgements}
The author would like to thank Prof. Michael R\"ockner for warm hospitality
during the visit to the University of Bielefeld
during August in 2018,
when part of this work is done.
The author also thanks Prof. Viorel Barbu, Jun Cao and Jiqiang Zheng for helpful discussions.
This work is supported by NSFC (No. 11501362, No. 11871337).
Financial support by the DFG through CRC 1283 is also gratefully acknowledged.

\end{document}